\pdfoutput=1

\documentclass[12pt,reqno]{amsart}
\usepackage{lmodern} 
\usepackage[scale=0.89]{tgheros} 
\usepackage[osf]{Baskervaldx}
\usepackage[baskervaldx,cmintegrals,bigdelims,vvarbb]{newtxmath} 
\usepackage[cal=boondoxo]{mathalfa}
         
  \usepackage{amsmath}          
  \usepackage{amsfonts}           
  \usepackage{amsthm}
  \usepackage{amscd}
  \usepackage{enumerate}
  \usepackage[active]{srcltx}
  \usepackage[percent]{overpic}
  \usepackage[mathscr]{eucal}
  \usepackage{graphicx,epsfig} 
  \usepackage{tikz-cd}        
  \usepackage{color}
  \usepackage{marginnote}
  \usepackage{hyperref}
  \usepackage{bm,bbm}
  \usepackage{enumitem} 
  \usepackage{tikz-cd}

\newtheorem{theorem}{Theorem}[section]
\newtheorem{corollary}[theorem]{Corollary}
\newtheorem{lemma}[theorem]{Lemma}
\newtheorem{proposition}[theorem]{Proposition}
\theoremstyle{definition}
\newtheorem{definition}[theorem]{Definition}
\newtheorem{example}[theorem]{Example}
\theoremstyle{remark}
\newtheorem{remark}[theorem]{Remark}
\newtheorem{question}[theorem]{Question}
\newtheorem{conjecture}{Conjecture}
\theoremstyle{thmx}
\newtheorem{thmx}{Theorem}

\newtheorem*{theorem*}{Theorem}

\numberwithin{equation}{section}
  \def\C{\mathbb{C}}
  \def\D{\mathbb{D}}
  \def\N{\mathbb{N}}
  \def\R{\mathbb{R}}
  \def\Q{\mathbb{Q}}
  \def\Z{\mathbb{Z}}    
  
  \def\H{\mathbb{H}}

  \def\XX{\mathcal{X}}
  \def\YY{\mathcal{Y}}
  
  \def\OO{\mathcal{O}}
  \def\MM{\mathcal{M}}
  \def\RR{\mathcal{R}}
  \def\PP{\mathcal{P}}
  \def\CC{\mathcal{C}}
  \def\HH{\mathcal{H}}
  \def\WW{\mathcal{W}}
  \def\LL{\mathcal{L}}
  \def\KK{\mathcal{K}}
  \def\NN{\mathcal{L}_{\operatorname{ren}}}

  \def\m{\mathbf{m}} 
  \def\e{\mathrm{e}}
  \def\ii{\mathrm{i}}

  \def\Mbrot{{\mathcal{M}}}
  \def\odd{{\textrm{odd}}}

  \def\Chat{\hat{\mathbb{C}}}

  \def\De{\Delta}
  \def\eps{\epsilon}
  \def\ga{\gamma}
  \def\Ga{\Gamma}
  \def\om{\omega}
  \def\Om{\Omega}
  \def\la{\lambda}
  \def\La{\Lambda}
  \def\si{\sigma}
  \def\Si{\Sigma}

	\newcommand{\vs}{\vspace{6pt}}
	\newcommand{\bit}{\it \bfseries}
	\newcommand{\es}{\emptyset}
	\newcommand{\sm}{\smallsetminus}
	\newcommand{\ov}{\overline}
	\newcommand{\bd}{\partial}
	
	\newcommand{\modd}{\ (\operatorname{mod} {\mathbb Z})}
	
	\newcommand{\dist}{\operatorname{dist}}
	\newcommand{\area}{\operatorname{area}}
  
	\newcommand{\per}{\operatorname{Per}}
 \newcommand{\re}{\operatorname{Re}}
 \newcommand{\im}{\operatorname{Im}}
 \newcommand{\Mod}{\operatorname{mod}}
 \newcommand{\kcen}{K_{\operatorname{cen}}} 
 \newcommand{\kper}{K_{\operatorname{per}}}
 \newcommand{\lcen}{{\mathcal L}_{\operatorname{cen}}} 
 \newcommand{\lper}{{\mathcal L}_{\operatorname{per}}}
 \newcommand{\chicen}{\chi_{\operatorname{cen}}} 
 \newcommand{\chiper}{\chi_{\operatorname{per}}}
 \newcommand{\Chi}{\boldsymbol{\chi}}
 \newcommand{\rin}{r_{\operatorname{\, in}}}
 \newcommand{\rout}{r_{\operatorname{out}}}

\newcommand{\ALIGN}{\begin{align*}}
\newcommand{\ENDALIGN}{\end{align*}}
\newcommand{\ENUM}{\begin{enumerate}}
\newcommand{\ENUMa}{\begin{enumerate}[a.]}
\newcommand{\ENUMA}{\begin{enumerate}[A.]}
\newcommand{\ENUMAF}{\begin{enumerate}[\bf A.]}
\newcommand{\ENUMi}{\begin{enumerate}[i)]}
\newcommand{\ENDENUM}{\end{enumerate}}
\newcommand{\ITMZ}{\begin{itemize}}
\newcommand{\ENDITMZ}{\end{itemize}}
\newcommand{\REFEQN}[1] { \begin{equation}\label{#1} }
\newcommand{\ENDEQN}{\end{equation}}
\newcommand{\THM}{\begin{theorem}}
\newcommand{\EX}{ \begin{example}}
\newcommand{\REFEX}[1] { \begin{example}\label{#1} }
\newcommand{\ENDEX}{\end{example}}
\newcommand{\MMTX}{ \begin{matrix}}
\newcommand{\ENDMTX}{ \end{matrix}}
\newcommand{\REM}{ \begin{remark}}
\newcommand{\ENDREM}{\end{remark}}
\newcommand{\REFTHM}[1] { \begin{theorem}\label{#1} }
\newcommand{\RREFTHM}[2] { \begin{theorem}[#1]\label{#2} }
\newcommand{\ENDTHM}{\end{theorem}}
\newcommand{\REFNTH}[1] { \begin{newthm}\label{#1} }
\newcommand{\ENDNTH}{\end{newthm}}
\newcommand{\REFPROP}[1]{\begin{proposition}\label{#1} }
\newcommand{\RREFPROP}[2]{\begin{proposition}[#1]\label{#2} }
\newcommand{\PROP}{\begin{proposition}}
\newcommand{\ENDPROP}{\end{proposition} }
\newcommand{\REFDEF}[1]{\begin{definition}\label{#1} }
\newcommand{\DEF}{\begin{definition}}
\newcommand{\ENDEF}{\end{definition} }
\newcommand{\REFLEM}[1]{\begin{lemma}\label{#1} }
\newcommand{\RREFLEM}[2]{\begin{lemma}[#1]\label{#2} }
\newcommand{\LEM}{\begin{lemma}}
\newcommand{\ENDLEM}{\end{lemma} }
\newcommand{\REFCOR}[1]{\begin{corollary}\label{#1} }
\newcommand{\COR}{\begin{corollary}}
\newcommand{\ENDCOR}{\end{corollary} }
\newcommand{\CONJ}{\begin{conjecture}}
\newcommand{\REFCONJ}[1]{\begin{conjecture}\label{#1}}
\newcommand{\RREFCONJ}[2]{\begin{conjecture}{#1}\label{#2}}
\newcommand{\ENDCONJ}{\end{conjecture} }
\newcommand{\QS}{\begin{question}}
\newcommand{\ENDQS}{\end{question} }
\newcommand{\REFDEFTHM}[1] { \begin{defthm}\label{#1} }
\newcommand{\ENDEFTHM}{\end{defthm}}
\newcommand{\corref}[1]{Corollary~\ref{#1}}

\newcommand{\exaref}[1]{Example~\ref{#1}}
\newcommand{\figref}[1]{Fig.~\ref{#1}}
\newcommand{\lemref}[1]{Lemma~\ref{#1}}
\newcommand{\remref}[1]{Remark~\ref{#1}}
\newcommand{\thmref}[1]{Theorem~\ref{#1}}

\newcommand{\PROOF}{\begin{proof}}
\newcommand{\ENDPROOF}{\end{proof}}

\newcommand{\Grotszch}{Gr{\"o}tszch }
\newcommand{\Bottcher}{B{\"o}ttcher }

\definecolor{uibgreen}{RGB}{119, 175, 0}

\begin{document}

\title[]{Lemon limbs of the cubic connectedness locus} \vs

\author[]{Carsten Lunde Petersen and Saeed Zakeri}

\address{Department of Mathematical Sciences, University of Copenhagen, Universitetsparken 5, DK-2100 Copenhagen {{\O}}, Denmark}

\email{lunde@math.ku.dk} 

\address{Department of Mathematics, Queens College of CUNY, 65-30 Kissena Blvd., Queens, NY 11367, USA \\ 
and The Graduate Center of CUNY, 365 Fifth Ave., New York, NY 10016, USA}

\email{saeed.zakeri@qc.cuny.edu}

\date{\today}

\begin{abstract}
We describe a limb structure in the connectedness locus of complex cubic polynomials, where the limbs are indexed by the periodic points of the doubling map $t \mapsto 2t \modd$. The main renormalization locus in each limb is parametrized by the product of a pair of (punctured) Mandelbrot sets. This parametrization is the inverse of the straightening map and can be thought of as a tuning operation that manufactures a unique cubic of a given combinatorics from a pair of quadratic hybrid classes. Our results build a new understanding of the structure of the cubic connectedness locus and its hyperbolic components.            
\end{abstract}

\maketitle

\tableofcontents

\section{A bird's-eye view}\label{sec:intro}

In this paper we define and investigate a limb structure in the connectedness locus of complex cubic polynomials that is inspired by the well-known case of the Mandelbrot set $\MM$ described by Douady and Hubbard (see \cite{DH1} and \cite{M3}). In their description $\MM$ is essentially partitioned into a main hyperbolic component---the ``main cardioid''---and a countable collection of disjoint limbs attached to it. These limbs are uniquely indexed by the non-trivial rotation cycles of the doubling map $\m_2: t \mapsto 2t \modd$ or, equivalently, by their non-zero rotation numbers. A parameter $c \in \MM$ belongs to the limb associated with a given rotation cycle precisely when the dynamic rays of the quadratic polynomial $Q_c: z \mapsto z^2+c$ with angles in that rotation cycle co-land at a fixed point. This description allowed for a divide-and-conquer strategy for understanding the structure of $\MM$. For example, it was the starting point of Yoccoz's celebrated work on the problem of local connectivity of the Mandelbrot set. \vs

For cubic maps the situation is far more complicated. Consider the family $\PP(3) \cong \C^2$ of complex cubic polynomials in the normal form $P_{a,b}: z \mapsto z^3 + 3a z^2+ b$, with $a,b \in \C$, and its connectedness locus $\CC(3)$, the cubic analogue of the Mandelbrot set. Branner and Hubbard \cite{BH} showed that $\CC(3)$ is compact and cellular, hence connected. But soon after Lavaurs \cite{La} proved that $\CC(3)$ is not locally connected and therefore its topology must be quite intricate. Up until now, there seems to have been no attempt at understanding $\CC(3)$ through a suitable notion of limbs or open wakes that should contain limbs except for their root points (compare the recent paper \cite{BOTW} which constructs a locally connected combinatorial model of $\CC(3)$ analogous to Thurston's model of $\MM$).  \vs
  
\begin{figure}[t]
\centering
\begin{overpic}[width=\textwidth]{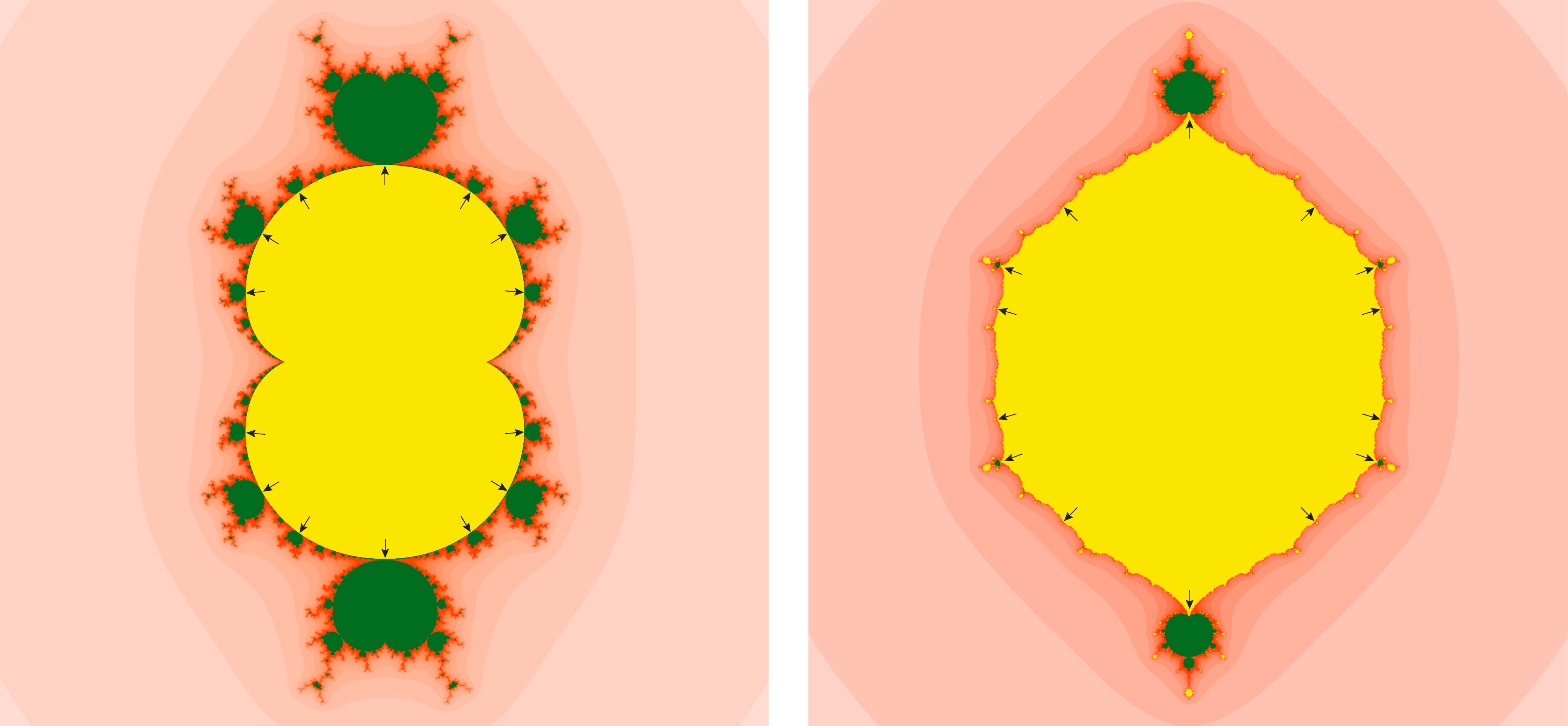}
\put (35.6,27) {\scalebox{0.5}{$L(1/4)$}}
\put (35.5,33) {\scalebox{0.5}{$L(1/3)$}}
\put (30.5,36) {\scalebox{0.5}{$L(2/5)$}}
\put (22.6,44) {\scalebox{0.5}{$L(1/2)$}}
\put (14.8,36) {\scalebox{0.5}{$L(3/5)$}}
\put (9.8,33) {\scalebox{0.5}{$L(2/3)$}}
\put (9.7,27) {\scalebox{0.5}{$L(3/4)$}} 
\put (35.6,18.2) {\scalebox{0.5}{$L^*(3/4)$}}
\put (35.5,12.7) {\scalebox{0.5}{$L^*(2/3)$}}
\put (30.5,9.8) {\scalebox{0.5}{$L^*(3/5)$}}
\put (22.6,1.5) {\scalebox{0.5}{$L^*(1/2)$}}
\put (14.2,9.8) {\scalebox{0.5}{$L^*(2/5)$}}
\put (9.2,12.7) {\scalebox{0.5}{$L^*(1/3)$}}
\put (9.1,18.2) {\scalebox{0.5}{$L^*(1/4)$}}
\put (77.3,2) {\scalebox{0.5}{$\LL_0(0)$}}
\put (84.3,12) {\scalebox{0.5}{$\LL_0(1/5)$}}
\put (90,16) {\scalebox{0.5}{$\LL_0(1/3)$}}
\put (89,19.1) {\scalebox{0.5}{$\LL_0(2/5)$}}
\put (89,26.3) {\scalebox{0.5}{$\LL_0(3/5)$}}
\put (90,29.7) {\scalebox{0.5}{$\LL_0(2/3)$}}
\put (84.3,33.2) {\scalebox{0.5}{$\LL_0(4/5)$}}
\put (77.3,43) {\scalebox{0.5}{$\LL_0^*(0)$}}
\put (62.7,33.2) {\scalebox{0.5}{$\LL_0^*(1/5)$}}
\put (57.2,29.7) {\scalebox{0.5}{$\LL_0^*(1/3)$}}
\put (58.2,26.3) {\scalebox{0.5}{$\LL_0^*(2/5)$}}
\put (58.2,19.1) {\scalebox{0.5}{$\LL_0^*(3/5)$}}
\put (57.2,16) {\scalebox{0.5}{$\LL_0^*(2/3)$}}
\put (62.7,12) {\scalebox{0.5}{$\LL_0^*(4/5)$}}
\put (74.5,22.5) {\small $\HH_0$}
\end{overpic}
\caption{\small Two transversal one-dimensional slices of the cubic connectedness locus $\CC(3) \subset \C^2$. Left: The unicritical family $\{ P_{0,b}: z \mapsto z^3+b \}_{b \in \C}$, where the limbs attached to the main hyperbolic component are indexed by non-zero rotation numbers. Here $b$ belongs to the limb $L(p/q)$ iff the $q$ dynamic rays of $P_{0,b}$ that co-land at a fixed point are permuted with rotation number $p/q$ and have their angles in the interval $]0,1/2[$. Right: The lemon family $\{ P_{a,0} : z \mapsto z^3+3az^2 \}_{a \in \C}$, where the limbs are indexed by the periodic points of the doubling map $\m_2: t \mapsto 2t \modd$. Here $a$ belongs to the limb $\LL_0(t)$ iff the dynamic rays of $P_{a,0}$ with the ``simulating angles'' of $t$ co-land (see \S \ref{subsec:sp}). Both pictures have a $180^{\circ}$ rotational symmetry since $P_{a,b}$ and $P_{-a,-b}$ are affinely conjugate. This gives rise to a countable collection of rotated limbs designated by asterisks.}  
\label{mlimb}
\end{figure}

To illustrate the varied nature of limbs in the cubic case, consider the ``unicritical family'' $\{ P_{0,b}: z \mapsto z^3 +b \}_{b \in \C}$ where the two critical points coincide, and the ``lemon family'' $\per_1(0)=\{ P_{a,0}: z \mapsto z^3+3az^2 \}_{a \in \C}$ where one critical point is fixed at $0$ and the other is free. These families can be identified with complex lines in $\C^2$ that meet transversally at the center of the main hyperbolic component $\HH$ of $\CC(3)$ corresponding to the map $z \mapsto z^3$. Both of these lines intersect $\CC(3)$ along a closed topological disk (the closure of their intersection with $\HH$) decorated with a countable collection of disjoint limbs. But the limbs in these one-dimensional slices of $\CC(3)$ have completely different dynamical and topological natures. In the unicritical case, the limbs are naturally indexed by non-zero rational rotation numbers $p/q$ plus an element of $\Z/2\Z$, where the family has a period-$q$ bifurcation at the root parameter, a situation that bears strong resemblance to the case of the Mandelbrot set (\figref{mlimb} left). For the lemon family, the limbs are naturally indexed by the periodic {\it points} of $\m_2$ plus an element of $\Z/2\Z$ and the family has a saddle-node bifurcation at the root parameter (\figref{mlimb} right). Combinatorially speaking, the unicritical limbs are associated with certain cycles of the tripling map $\m_3: t \mapsto 3t \modd$ whose combinatorics have degree $1$, while the limbs in $\per_1(0)$ are associated with pairs of adjacent cycles of $\m_3$ that can have any combinatorics of degree $1$ or $2$ (see \S \ref{copo} for details). \vs

%%%%%%%%%%%%%%%%%%%%%%%%%%%%%%%%%%%%%%%%%%%
%%%%%%%%%%%%%%%%%%%%%%%%%%%%%%%%%%%%%%%%%%% 

We use the combinatorial characterization of limbs in $\per_1(0)$ to define the ``lemon limb'' $\LL(t) \subset \CC(3)$ for every periodic point $t \in \Q/\Z$ of $\m_2$. The combinatorial pattern of co-landing rays   of $P=P_{a,b}$ in $\LL(t)$ implies the existence of two quadratic-like restrictions, a ``central'' $P:U \to U'$ about the first critical point $0$ and a ``peripheral'' $P^{\circ q}: V \to V'$ about the second critical point $-2a$, where $q \geq 1$ is the period of $t$ under $\m_2$. Under favorable conditions these restrictions have connected filled Julia sets and therefore they are hybrid conjugate to unique quadratic polynomials $Q_c$ and $Q_\nu$, respectively. This defines a straightening map $\Chi: P \mapsto (c,\nu)$ from the main renormalization locus $\NN(t) \subset \LL(t)$ to the product of two copies of $\MM$. We prove that this map is a homeomorphism between $\NN(t)$ and the product $\check{\MM}_t \times \check{\MM}_0$. Here $\check{\MM}_t = \MM \sm \{ c(t_1), \ldots , c(t_q) \}$, where $\{ t_1, \ldots, t_q \}$ is the $\m_2$-orbit of $t$ and $c(t_i)$ is the parabolic landing point of the external ray of $\MM$ at angle $t_i$ (in particular, $\check{\MM}_0=\MM \sm \{ 1/4 \}$). Thus, the inverse $\Chi^{-1}$ gives a parametrization of $\NN(t)$ by the product of two finitely punctured Mandelbrot sets. The assignment $\Chi^{-1}:(c,\nu) \mapsto P$ can be thought of as a ``tuning'' operation that manufactures a unique cubic $P$ with a peripheral quadratic hybrid class $Q_\nu$ having the combinatorics of $t$, tuned by a central quadratic hybrid class $Q_c$. \vs

Our results extend and tie together various pieces of prior knowledge. For example, the one-dimensional limbs $\LL_0(t):=\LL(t) \cap \per_1(0)$ have been studied in the works of Milnor \cite{M4} and Roesch \cite{R2}, and our introduction of the lemon limbs $\LL(t)$ can be thought of as a direct generalization from dimension $1$ to $2$. Inou and Kiwi \cite{IK} have studied properties of the straightening map in the general framework of mapping schemata of postcritically finite polynomials and, among other things, shown its injectivity. For cubic polynomials, they proved the existence of homeomorphic copies of $\MM \times \MM$ in $\CC(3)$ associated with the hyperbolic components of primitive disjoint type. Our work on lemon limbs provides a complementary result for {\it non-primitive} disjoint type renormalizations, where one has to deal with the issue of non-compactness of the domain of the straightening map. The paper \cite{IK} also showed the existence of injected copies of $\KK := \{ (c,z) : c \in \MM \ \text{and} \ z \in K_c \}$ in $\CC(3)$ associated with the hyperbolic components of primitive capture type (see also \cite{BuHe}). Such copies of $\KK$ can be observed in our full lemon limb $\LL(t)$ as the extensions of the capture components in $\LL_0(t)$. Another related result is the intertwining surgery studied by Epstein and Yampolsky \cite{EY}, and by Ha\"{i}ssinsky \cite{Ha}, which can be interpreted as a special case of our construction of $P=\Chi^{-1}(c,\nu)$ when the combinatorics of $t$ under $\m_2$ is of rotation type $p/q$ and the central hybrid class $Q_c$ is chosen from the $p/q$-limb of $\MM$. Finally, lemon limb structures show up in a recent paper of Bonifant, Estabrooks and Sharland \cite{BES} which establishes a correspondence between escape loci in $\per_1(0)$ and in the curve $\per_2(0) \subset \PP(3)$ consisting of cubics with a critical point of period $2$.% 
\vs

The ideas presented in this paper open up several new directions and avenues for further exploration. We refer the reader to a list of problems at the end of our expanded introduction in the next section.

\section{Outline of results}\label{sec:oor}

We now give a more precise account of the main results. We represent cubic polynomials in the normal form 
$$
P_{a,b}: z \mapsto z^3 + 3a z^2+ b \qquad (a,b \in \C)
$$
with marked critical points $\om_1=0, \om_2=-2a$. The space $\PP(3)$ of all such critically marked cubics is isomorphic to $\C^2$ with coordinates $(a,b)$. The space of all affine conjugacy classes of these cubics (the so-called ``moduli space'') is also isomorphic to $\C^2$ with coordinates $(a^2,b^2)$. The {\bit cubic connectedness locus} $\CC(3)$ is the set of all cubics in $\PP(3)$ with connected Julia set, which is known to be compact and connected (in fact, cellular). Consider the complex line $\per_1(0) \subset \PP(3)$ consisting of all cubics $P_{a,b}$ for which $\om_1$ is fixed, with the normal form 
$$
P_a=P_{a,0}: z \mapsto z^3+3az^2 \qquad (a \in \C).
$$
The connectedness locus $\CC_0=\CC(3) \cap \per_1(0)$ of this family is shown in Figs. \ref{mlimb} and \ref{lemonfig}. The prominent central region $\HH_0$ in these figures, resembling the shape of a lemon, consists of all $a \in \C$ for which $\om_2$ belongs to the immediate basin of attraction of $\om_1$ (equivalently, the Julia set of $P_a$ is a Jordan curve). For this reason, we informally refer to $\{ P_a \}_{a \in \C}$ as the {\bit lemon family} of cubics. The limbs of $\CC_0$ attached to $\bd \HH_0$ have been the subject of several investigations (see \cite{F}, \cite{M4}, and \cite{R2}). Although structurally different, they are reminiscent of the primary limbs of the Mandelbrot set attached to the boundary of the main cardioid. \vs

The ``lemon limbs'' of the connectedness locus $\CC(3)$ that we study in this paper are the natural extensions of the one-dimensional limbs of $\CC_0$. Their description depends on a basic combinatorial result on the existence and uniqueness of pairs of periodic orbits of the tripling map $\m_3: t \mapsto 3t$ (mod $\Z$) that simulate the combinatorics of the periodic points of the doubling map $\m_2$. By the {\bit combinatorics} of a periodic orbit $\{ t_1, \ldots, t_q \}$ under a covering map $f: \R/\Z \to \R/\Z$, labeled so that $0 \leq t_1<\cdots<t_q<1$, we mean the $q$-cycle $\sigma$ in the permutation group $S_q$ determined by the condition $f(t_i)=t_{\sigma(i)}$ for all $i$. \vs

In \S \ref{subsec:sp} we use the general machinery developed in \cite{PZ1} to prove the following 

\begin{thmx} [Existence and uniqueness of simulating orbits] \label{A}
Let $t \in \Q/\Z$ be a period $q$ point of $\m_2$ with the orbit $\{ t_1,\ldots, t_k=t,\ldots,t_q \}$ and combinatorics $\sigma \in S_q$. Then $t$ has a unique pair of ``simulating orbits,'' i.e., period $q$ orbits $\{ x_1, \ldots, x_q \}$ and $\{ y_1, \ldots, y_q \}$ under $\m_3$, having the same combinatorics $\sigma$, that interlace in the sense  
$$
0<x_1<y_1<\cdots<x_k<\frac{1}{2}<y_k<\cdots<x_q<y_q<1.
$$
\end{thmx}

The theorem may be interpreted as saying that the $\m_2$-orbit $\{ t_1, \ldots, t_q \}$ can be ``blown up'' to a pair of $\m_3$-orbits of the same combinatorics so that the marked point $t_k$ blows up to a pair $x_k,y_k$ surrounding the fixed point $1/2$ (compare Figures \ref{one} and \ref{two} and the associated examples). \vs 

In fact, there are precisely $q+1$ period $q$ orbits $\OO_0, \ldots, \OO_q$ of $\m_3$ that have the combinatorics $\sigma$; they are uniquely characterized by the property $\#(\OO_i \cap [0,1/2[)=i$. We prove \thmref{A} by showing that the neighboring pairs $\OO_{i-1}, \OO_i$ interlace, so $\{ x_1, \ldots, x_q \}=\OO_k$ and $\{ y_1, \ldots, y_q \}=\OO_{k-1}$ will be the simulating orbits of $t$ (\thmref{interlace}). \vs 

We note that the idea of constructing orbits of $\m_3$ from those of $\m_2$ (and more generally $\m_{d+1}$ from $\m_d$) for a combinatorics of rotation type has appeared in \cite{BMMOP}. \vs 

With $t=t_k$ and $x_k,y_k$ as above, the limb $\LL_0(t) \subset \per_1(0)$ attached to $\bd \HH_0$ can be described as the set of cubics $P_a \in \CC_0$ for which the dynamic rays $R_a(x_k), R_a(y_k)$ co-land (see \S \ref{lfdyn} and \S \ref{lfpar}). Inspired by this characterization, we define the {\bit lemon limb} $\LL(t)$ by
$$
\LL(t) = \big\{ P \in \CC(3) : \text{the dynamic rays} \ R_P(x_k), R_P(y_k) \ \text{co-land} \, \big\}.
$$
Unlike the slices $\LL_0(t)$ or the primary limbs of the Mandelbrot set, the lemon limbs $\LL(t)$ are not compact (see \exaref{noncompact}). In fact, for $P \in \ov{\LL(t)} \sm \LL(t)$ the rays $R_P(x_k), R_P(y_k)$ land at distinct points, one of which being parabolic. For an extended discussion of this phenomenon we refer to \S \ref{propchi}. \vs
 
When $P \in \LL(t)$, the dynamic rays $R_P(x_i),R_P(y_i)$ co-land at a point $z_i$ for every $1 \leq i \leq q$ and the {\bit co-landing orbit} $\{ z_1, \ldots, z_q \}$ has period dividing $q$, with each $z_i$ either repelling or parabolic with multiplier $(P^{\circ q})'(z_i)=1$. Define the dynamical wake $W_i$ as the connected component of $\C \sm (R_P(x_i) \cup R_P(y_i) \cup \{ z_i \})$ not containing the first critical point $\om_1$. The fact that the interval $(x_k,y_k) \subset \R/\Z$ contains the fixed point $1/2$ of $\m_3$ implies that the wake $W_k$ contains the second critical point $\om_2$. Under $P$ each wake $W_i$ maps conformally to $W_{\sigma(i)}$ if $i \neq k$ and to the entire plane if $i=k$. These wakes can be used to define two renormalizations of $P$ of periods $1$ and $q$ as follows (see \S \ref{cpren}). If $\{ z_1, \ldots, z_q \}$ is ``peripherally repelling,'' i.e., repelling or parabolic with no immediate basin in $\bigcup_{i=1}^q W_i$, there is a quadratic-like restriction $P:U \to U'$ around $\om_1$, where $U'$ is a slightly thickened copy of $\C \sm \bigcup_{i=1}^q \ov{W_i}$ cut off by some equipotential. We call $P$ {\bit centrally renormalizable} if the filled Julia set $\kcen$ of $P:U \to U'$ is connected. Similarly, if $\{ z_1, \ldots, z_q \}$ is ``centrally repelling,'' i.e., repelling or parabolic with no immediate basin outside $\bigcup_{i=1}^q W_i$, there is a quadratic-like restriction $P^{\circ q}: V \to V'$ around $\om_2$, where $V'$ is a slightly thickened copy of $W_k$ cut off by some equipotential. We call $P$ {\bit peripherally renormalizable} if the filled Julia set $\kper$ of $P^{\circ q}: V \to V'$ is connected. The {\bit main renormalization locus} $\NN(t)$ is the set of all cubics in $\LL(t)$ that are both centrally and peripherally renormalizable. \vs

According to the straightening theorem of Douady and Hubbard \cite{DH2}, for each $P \in \NN(t)$ there is a unique pair $(c,\nu) \in \MM \times \MM$ such that the central renormalization $P: U \to U'$ is hybrid equivalent to $Q_c: z \mapsto z^2+c$ and the peripheral renormalization $P^{\circ q}: V \to V'$ is hybrid equivalent to $Q_\nu: z \mapsto z^2+\nu$. This gives two straightening maps
$$
\chicen, \chiper : \NN(t) \to \MM
$$
defined by $\chicen(P):=c, \chiper(P):=\nu$. If $\varphi$ is a hybrid equivalence between $P: U \to U'$ and $Q_c$, then for every $1 \leq i \leq q$ the dynamic ray $R_c(t_i)$ of $Q_c$ lands at the repelling point $w_i:=\varphi(z_i)$ of $Q_c$, and if $\psi$ is a hybrid equivalence between $P^{\circ q}: V \to V'$ and $Q_\nu$, then $\psi(z_k)$ is the $\beta$-fixed point of $Q_\nu$ (\lemref{renorm}). Since $\{ z_1, \ldots, z_q \}$ is repelling, we conclude that the periodic orbit $\{ w_i \}$ of $Q_c$ and the $\beta$-fixed point of $Q_\nu$ must be repelling. It follows from the general theory of the Mandelbrot set that $\nu \neq 1/4$ and $c \neq c(t_i)$ for all $1 \leq i \leq q$, where $c(t_i)$ is the parabolic landing point of the external ray $\RR_\MM(t_i)$ of the Mandelbrot set. Let  
$$
\check{\MM}_t := \MM \sm \{ c(t_1), \ldots, c(t_q) \}.
$$   
In particular, $\check{\MM}_0=\MM \sm \{ 1/4 \}$.   

\begin{thmx}[Parametrization of the main renormalization locus] \label{B}
The straightening map
$$
\Chi: \NN(t) \to \check{\MM}_t \times \check{\MM}_0
$$ 
defined by $\Chi(P)=(\chicen(P), \chiper(P))$ is a homeomorphism.  
\end{thmx}

The inverse map $\Chi^{-1}$ gives a homeomorphic embedding of the product $\check{\MM}_t \times \check{\MM}_0$ into the lemon limb $\LL(t)$. The construction of $P=\Chi^{-1}(c,\nu)$ amounts to the aforementioned tuning operation: it is the unique cubic with the peripheral hybrid class of $Q_\nu$ with the combinatorics of $t$, tuned by the central hybrid class of $Q_c$. \vs 

Our proof of \thmref{B} is structured as follows:  \vs 

\begin{enumerate}[leftmargin=*]
\item[$\bullet$]
In \S \ref{injchi} we prove injectivity of $\Chi$ by a pull-back argument. The key ingredient is to show that for each $P \in \NN(t)$ the residual Julia set 
$$
K_\infty := K_P \sm \bigcup_{n \geq 0} P^{-n}(\kcen \cup \kper)
$$
has Lebesgue measure zero (\lemref{K=null}). The required bounded distortion argument is handled using hyperbolic geometry estimates for the inclusion map $i : X \hookrightarrow Y$ between Riemann surfaces (Lemmas \ref{lip} and \ref{lipcor}). \vs 

\item[$\bullet$] 
In \S \ref{propchi} we show that the inverse $\Chi^{-1}$ is continuous, that is, if $P_n,P_\infty \in \NN(t)$ and $\Chi(P_n) \to \Chi(P_\infty)$, then $P_n \to P_\infty$. Since $\Chi$ is a continuous injection, it suffices to prove that if $P_n$ accumulates on some $P \in \CC(3)$, then $P \in \NN(t)$. The non-trivial step here is to ensure that the co-landing orbit of $P_n$ does not degenerate into a parabolic orbit of $P$ and that the rays $R_P(x_k), R_P(y_k)$ continue to co-land at the limit of the co-landing point of $R_{P_n}(x_k), R_{P_n}(y_k)$ (Lemmas \ref{onelands} and \ref{bothland}). The latter amounts to ruling out the creation of heteroclinic arcs in the Hausdorff limit of the pair $R_{P_n}(x_k), R_{P_n}(y_k)$ as $n \to \infty$.\footnote{This issue seems to have been overlooked in the special case of intertwining cubics in \cite{EY}.} \vs 

\item[$\bullet$] 
In \S \ref{surg} we prove surjectivity of $\Chi$. Given $(c,\nu) \in \check{\MM}_t \times \check{\MM}_0$, we first construct a holomorphic map $f: \C \sm \ov{\La} \to \C$ with a central renormalization of period $1$ hybrid equivalent to $Q_c$ and a peripheral renormalization of period $q$, with the combinatorics of $t$, hybrid equivalent to $Q_\nu$ (see \S \ref{mec}). Here the ``hole'' $\La$ is a Jordan domain which, synthetically speaking, plays the role of the preimage of $\kcen$ other than $\kcen$. We then use an auxiliary dynamics near $\bd \La$ and a careful quasiconformal interpolation to extend $f$ across a sector containing $\La$ (see \S \ref{bcon} and \S \ref{extf}). This gives a quasiregular polynomial-like map of degree $3$ whose straightening, after suitable normalization, will be the desired $P \in \NN(t)$ with $\Chi(P)=(c,\nu)$.   \vs     
\end{enumerate}

Every periodic point $t \in \Q/\Z$ of $\m_2$ determines an $\m_3$-invariant lamination defined by joining the corresponding pairs $x_i,y_i$ given by \thmref{A}. This lamination can be thought of as the combinatorial signature of the lemon limb $\LL(t)$. The combinatorics arising here is intrinsically quadratic in the sense that the permutation in $S_{2q}$ induced by the union $\{ x_1, \ldots, x_q \} \cup \{ y_1, \ldots, y_q \}$ has degree $\leq 2$ (see \S \ref{copo} for the notion of degree). Other pairs of periodic orbits of $\m_3$ can still define co-landing patterns but the combinatorics of their union necessarily has degree $3$ (compare the Appendix). It would be desirable to understand which pairs are compatible and to investigate the resulting laminations. \vs

\noindent
{\bf Some open problems.} We conclude this introduction with a list of open problems that suggest possible directions for further exploration and a path toward a divide-and-conquer approach to $\CC(3)$. 
 
\begin{enumerate}[leftmargin=*]
\item
Study the topology of lemon limbs. In particular, is $\ov {\LL(t)}$ connected? Is $\LL(t)$ connected? Is it contractible? A relevant notion here is that of the ``lemon wake'' $\WW(t)$ consisting of all cubics $P \in \PP(3)$ for which the rays $R_P(x_k), R_P(y_k)$ are {\it smooth} and co-land at a {\it repelling} periodic point. By definition $\WW(t)$ is open. It can be shown that the following relations hold: 
$$
\WW(t) \cap \CC(3) \subsetneq \LL(t) \subsetneq \ov{\LL(t)} = \ov{\WW(t)} \cap \CC(3).
$$
Moreover, for $P \in \bd \WW(t) \cap \CC(3)$ one of the following occurs: \vs
\begin{enumerate}[leftmargin=*]
\item[(i)]
$R_P(x_k), R_P(y_k)$ co-land at a parabolic point. In this case $P \in \LL(t)$; \vs
\item[(ii)]
$R_P(x_k), R_P(y_k)$ land at distinct points, one of which being  parabolic. In this case $P \notin \LL(t)$. \vs 
\end{enumerate}     
For $P$ in the complementary boundary part $\bd \WW(t) \sm \CC(3)$ there are more possibilities. Both $R_P(x_k), R_P(y_k)$ can be smooth and co-landing at a parabolic point, or they can land at distinct points one of which being parabolic. Also one or both of these rays can be broken. \vs

\item
Give a dynamical characterization for the cubics in $\ov{\LL(t)} \sm \LL(t)$. This set is contained in $\bd \WW(t) \cap \CC(3)$ and therefore every cubic in it must satisfy the condition (ii) above, but there are actually more restrictions. See \remref{hetconj} for a discussion of these restrictions and a proposed characterization for this difference set. \vs

\item
Investigate the extension of the straightening map $\Chi$ to cubics in $\LL(t)$ with parabolic co-landing orbits. This appears to be a straightforward application of the parabolic-like mappings of \cite{Lo}. Along the same lines, it seems plausible that the inverse $\Chi^{-1}$ extends continuously to all parameters $(c,\nu)$ where at least one of the conditions $c \in \check{\MM}_t$ or $\nu \in \check{\MM}_0$ holds. \vs

\item
Each one-dimensional limb $\LL_0(t)=\LL(t) \cap \per_1(0)$ has a unique root point which is also the root of a copy $\MM_0(t)$ of the Mandelbrot set (see Figs. \ref{lemonfig} and \ref{onethirdlimb}). This ``main Mandelbrot copy'' can be interpreted as the closure of the fiber $\Chi^{-1}(\{ 0 \} \times \check{\MM}_0)$ of the straightening map. It follows from our \thmref{B} that there is a continuous motion $c \mapsto \Chi^{-1}(\{ c \} \times \check{\MM}_0)$ as $c$ varies over $\check{\MM}_t$, that is, the derooted main Mandelbrot copy persists when the central hybrid class moves from $c=0$ to anywhere in the Mandelbrot set with finitely many exceptional parabolic parameters. As noted above, we expect $\Chi^{-1}(c,\check{\MM}_0)$ to have a well-defined root, giving rise to a full copy of the Mandelbrot set. When is this copy quasiconformally homeomorphic to $\MM$ (in a suitable sense)? With parabolic-like mappings in mind, we speculate that when $c$ is an exceptional parameter, a corresponding Mandelbrot copy still exists but it is quasiconformally homeomorphic to the ``parabolic Mandelbrot set'' $\MM_1$ studied in \cite{PR}. \vs 

\item \label{holomotions}
The central renormalization locus in $\LL(t)$ can be naturally decomposed into one-dimensional slices $\{ \LL_c(t) \}_{c\in \MM}$ as follows. For $c \in \MM$, consider the central renormalization curve $\XX_c$ consisting of all cubics in $\PP(3)$ which have a quadratic-like restriction around $\om_1=0$ hybrid equivalent to $Q_c$. This is a one-dimensional complex analytic set in $\PP(3)$ whose escape region $\XX_c \sm \CC(3)$ is isomorphic to a punctured disk (when $c=0$, we recover the subset $\XX_0=\per_1(0) \sm \ov{\HH_0}$ of the complex line shown in Figs. \ref{mlimb} and \ref{lemonfig}). We define $\LL_c(t)$ as the intersection $\XX_c \cap \LL(t)$. For $c \in \check{\MM}_t$ the slice $\LL_c(t)$ can be viewed as an analytic extension of the embedded copy $\Chi^{-1}(c,\check{\MM}_0)$ of the derooted Mandelbrot set. \vs

By the work of Petersen and Tan \cite{PT} there is a holomorphic motion $\LL_0(t) \to \LL_c(t)$ as $c$ varies over the main cardioid of the Mandelbrot set $\MM$. More generally, one would expect $\LL_c(t)$ to move holomorphically as long as $c$ varies in an interior component of $\MM$. What is the structure of $\LL_c(t)$ when $c$ is on the boundary of $\MM$? Is it always a quotient of $\LL_0(t)$? \vs

\item \label{puz}
Give a description of the cubics in the limb $\LL(t)$ that are centrally renormalizable. We expect that in some sense these form a continuous motion of the parameters in the one-dimensional limb $\LL_0(t)=\LL(t) \cap \per_1(0)$ and often a holomorphic motion of them. However, as there will be collisions due to the change in the topology of the central filled Julia set, the holomorphic motions of quotients of $\LL_0(t)$ only exist when the central renormalization vary over an interior component of $\MM$ (compare (\ref{holomotions}) above). Similarly, give a description of the cubics in $\LL(t)$ that are peripherally renormalizable. There is an asymmetry between the two types of renormalizations, but the latter problem can be viewed as a companion of the former. As a further step, one could attempt to give a ``complete description'' of the whole limb $\LL(t)$ by analyzing what happens when the forward orbit of $\om_1$ or $\om_2$ exits the domain of the corresponding renormalization.\vs

A key tool in studying these problems is the machinery of {\it Yoccoz puzzles} (compare Roesch's work in \cite{R2} where she uses puzzles to study the slice $\per_1(0)$, although her construction is different from what we propose below). Each lemon limb $\LL(t)$ gives a separation of the two critical points by a pair of co-landing rays and thus the basis for constructing puzzle pieces, where the critical points are already separated by the depth-zero puzzle. More precisely, let $P \in \LL(t)$ and $x_k,y_k$ be the simulating angles of $t=t_k$. There are unique angles $x'_k,y'_k$ with $x_k<y'_k<1/2<x'_k< y_k$ such that $\m_3(x_k’)=x_{\sigma(k)}$ and $\m_3(y_k’)=y_{\sigma(k)}$ (see \S \ref{subsec:sp}). For each $j \neq k$, take the unique angles $x'_j,y'_j$ with $x_j<y_j’<x_j’<y_j$ such that $\m_3(x'_j)=x'_{\sigma(j)}$ and $ \m_3(y'_j)=y'_{\sigma(j)}$. Take the union of the $2q$ co-landing pairs $R_{x_j}, R_{y_j}$ and $R_{x'_j}, R_{y'_j}$ for $1 \leq j \leq q$ together with their landing points and some Green's equipotential curve. The depth-zero puzzle pieces are defined as the closure of the bounded connected components of the complement of this union. The depth-$n$ puzzle pieces are obtained by taking pull backs under $P^{\circ n}$. Denoting by $Y_i^P$ the depth-zero puzzle piece containing the critical point $\omega_i$, it follows that $P$ is centrally (resp. peripherally) renormalizable precisely when $\{ z_1, \ldots z_q \}$ is peripherally (resp. centrally) repelling and the forward orbit of $\om_1$ under $P$ (resp. $\om_2$ under $P^{\circ q}$) is contained in $Y_1^P$ (resp. $Y_2^P$). These dynamical puzzles induce parameter space puzzles, often called ``para-puzzles,'' in each wake $\WW(t)$ which promise to be a useful tool in studying the renormalization loci in $\LL(t)$.
\end{enumerate}

\section{Cubic orbits with quadratic combinatorics}\label{coqc}

The following notation will be adopted in this paper: \vs

\begin{enumerate}[leftmargin=*]
\item[$\bullet$]
$[t_1,t_2]$ is the closed interval from $t_1$ to $t_2$ if $t_1,t_2 \in \R$ and is the closed arc of the circle $\R/\Z$ from $t_1$ to $t_2$ in the counter-clockwise direction if $t_1,t_2 \in \R/\Z$. \vs 

\item[$\bullet$]
$A \Subset B$ means $\ov{A} \subset B$. \vs 

\item[$\bullet$]
$\m_k: t \mapsto k t \ \modd$ is the multiplication-by-$k$ map of the circle $\R/\Z$. Here $k$ is an integer $\geq 2$ \vs

\item[$\bullet$]
$Q_c$ is the normalized quadratic polynomial $z \mapsto z^2+c$, with the filled Julia set $K_c$, the Julia set $J_c=\bd K_c$, the \Bottcher coordinate $\phi_c$, and the external ray $R_c(\theta)$ at angle $\theta \in \R/\Z$. \vs

\item[$\bullet$]
$\MM = \{ c \in \C: 0 \in K_c \}$ is the Mandelbrot set, with the external rays $\RR_\MM(\theta)$. \vs

\item[$\bullet$]
$P=P_{a,b}$ is the normalized cubic polynomial $z \mapsto z^3+3az^2+b$ with marked critical points at $\om_1=0$ and $\om_2=-2a$, the filled Julia set $K_P$, the Julia set $J_P=\bd K_P$, the \Bottcher coordinate $\phi_P$, and the external rays $R_P(\theta)$. \vs 

\item[$\bullet$]
$\PP(3)$ and $\CC(3)$ are the cubic polynomial parameter space and its connectedness locus, respectively. \vs

\item[$\bullet$]
For a hyperbolic Riemann surface $X$, $\lambda_X=\lambda_X(z) \, |dz|$ is the hyperbolic metric with constant curvature $-1$, $\dist_X$ and $\area_X$ are the hyperbolic distance and area, and $B_X(w,r)$ is the hyperbolic $r$-ball centered at $w$. 
\end{enumerate}

\subsection{Combinatorics of periodic orbits}\label{copo}

Let $f: \R/\Z \to \R/\Z$ be an orientation-preserving covering map of the circle with a period $q$ orbit $\OO= \{ t_1, \ldots, t_q \}$. Unless otherwise stated, we label periodic orbits in positive cyclic order, that is, choosing unique representatives in $[0,1)$, we assume $0 \leq t_1 < \cdots < t_q <1$. By the {\bit combinatorics} of $\OO$ under $f$ we mean the cyclic permutation $\sigma$ in the symmetric group $S_q$ defined by 
$$
f(t_i)=t_{\sigma(i)} \qquad \text{for all} \ 1 \leq i \leq q.
$$
We often describe this situation by saying that the orbit $\OO$ {\bit realizes} $\sigma \in S_q$.  
      
\DEF
$\sigma \in S_q$ is said to have a {\bit descent} at $i$ if $\sigma(i) > \sigma(i+1)$. Here $1 \leq i \leq q$ and we always adopt the convention $q+1=1$). The number of descents of $\sigma$ is called its {\bit degree}:\footnote{In the combinatorics literature this invariant is called the {\it cyclic descent number} of $\sigma$.}
$$
\deg(\sigma) = \# \{ 1 \leq i \leq q: \sigma \ \text{has a descent at} \ i \}.  
$$ 
\ENDEF

For example, $\sigma=(1243) \in S_4$, which can be written in the classical notation as $\sigma=\begin{pmatrix} 1 & 2 & 3 & 4 \\ 2 & 4 & 1 & 3 \end{pmatrix}$, has descents at $i=2,4$ so $\deg(\sigma)=2$. \vs 

It is easy to see that $\deg(\sigma)=1$ if and only if $\sigma$ is a {\bit rotation cycle}, i.e., a cyclic permutation of the form $\sigma: i \mapsto i+p \ (\operatorname{mod} q)$ for some integer $1 \leq p <q$ relatively prime to $q$. The reduced fraction $p/q$ is called the {\bit rotation number} of $\sigma$. \vs 

The question of possible combinatorics of periodic orbits of the standard covering maps $\m_k: t \mapsto kt \ \modd$ and the frequency with which they occur is answered in \cite{PZ1}. The combinatorial invariant $\deg(\sigma)$ can be characterized topologically as the minimum degree of an orientation-preserving covering map of the circle having a periodic orbit which realizes $\sigma$. Thus, if $\sigma \in S_q$ is realized under $\m_k$, then necessarily $k \geq \deg(\sigma)$. The following theorem shows sufficiency of this condition and provides a precise count: 

\THM[\cite{PZ1}]\label{realize}
Let $\sigma \in S_q$ be a cyclic permutation with $\deg(\sigma)=d$. Then, for every integer $k \geq \max \{ d, 2 \}$ the number of period $q$ orbits of $\m_k$ having combinatorics $\sigma$ is 
$$
\begin{cases} 
\ \ \displaystyle{\binom{q+k-d}{q}} & \qquad \text{if} \ \sigma(q)>\sigma(1), \vspace{3mm} \\
\ \ \displaystyle{\binom{q+k-d-1}{q}} & \qquad \text{otherwise.}
\end{cases}
$$
\ENDTHM

The periodic orbits of $\m_k$ that realize a given $\sigma \in S_q$ are uniquely determined by their ``fixed point distribution,'' i.e., the knowledge of how the $k-1$ fixed points of $\m_k$ are distributed among the $q$ intervals in between the orbit points, or equivalently, how the orbit points are deployed in the $k-1$ intervals between the fixed points. \vs

The following special case of the above theorem plays a central role in this paper. Suppose a cyclic permutation $\sigma \in S_q$ can be realized under the doubling map $\m_2$. Then either $\deg(\sigma)=1$ so $\sigma$ is a rotation cycle and $\sigma(q)<\sigma(1)$, or $\deg(\sigma)=2$ and $\sigma(q)>\sigma(1)$. Whichever the case, \thmref{realize} shows that $\sigma$ has a unique realization under $\m_2$ and $q+1$ realizations under $\m_3$.

\REM \label{deggen}
The preceding definitions easily generalize to finite unions of periodic orbits in which case the associated combinatorics are products of cyclic permutations with disjoint supports. The definition of degree remains the same and $k \geq \deg(\sigma)$ is a necessary (but no longer sufficient) condition for $\m_k$ having a finite invariant set that realizes $\sigma$.        
\ENDREM    

\subsection{Simulating orbits of $\m_3$}\label{subsec:sp}   
 
Below we describe a process that associates to each periodic {\it point} of $\m_2$ a canonical pair of periodic {\it orbits} of $\m_3$, all with the same combinatorics. First suppose $\{ t_1, \ldots, t_q \}$ is a period $q$ orbit of $\m_2$ with combinatorics $\sigma \in S_q$. By the above discussion, there are precisely $q+1$ realizations $\OO_0, \ldots, \OO_q$ of $\sigma$ under $\m_3$. They are uniquely determined by the condition 
\begin{equation}\label{deploy}
\# (\OO_k \cap [0,1/2[) = k \qquad \text{for} \ 0 \leq k \leq q.
\end{equation}

The following result is a slight restatement of \thmref{A} in the introduction and will be proved at the end of this section:

\THM \label{interlace}
For each $1 \leq k \leq q$, the neighboring orbits $\OO_{k-1}, \OO_k$ interlace. More precisely, if $\OO_k = \{ x_1, \ldots , x_q \}$ and $\OO_{k-1} =\{ y_1, \ldots , y_q \}$, then
\begin{equation}\label{xy}
0\leq x_1<y_1<\cdots<x_k<\frac{1}{2} \leq y_k<\cdots<x_q<y_q<1.
\end{equation}
\ENDTHM

We call $\OO_k, \OO_{k-1}$ the {\bit simulating orbits} and $x_k, y_k$ the {\bit simulating angles} of $t_k$. \vs 

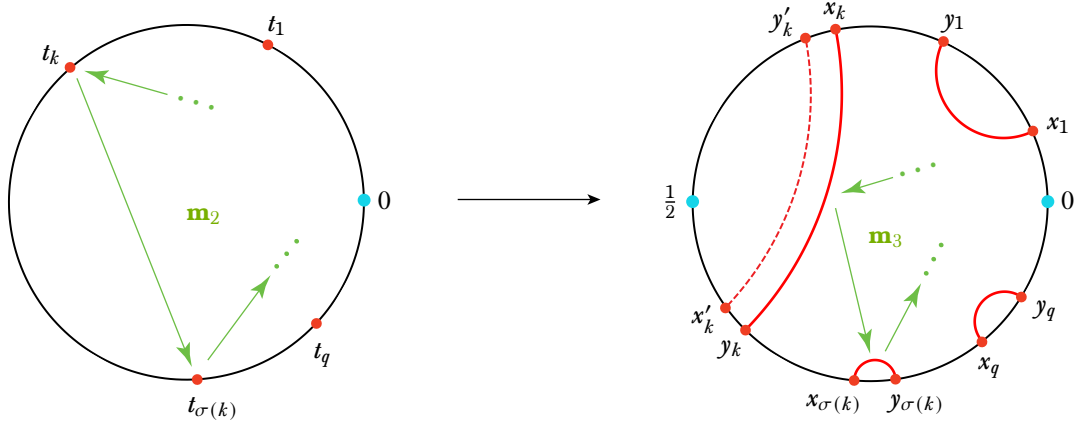
\begin{figure}[t]
	\centering
	\begin{overpic}[width=\textwidth]{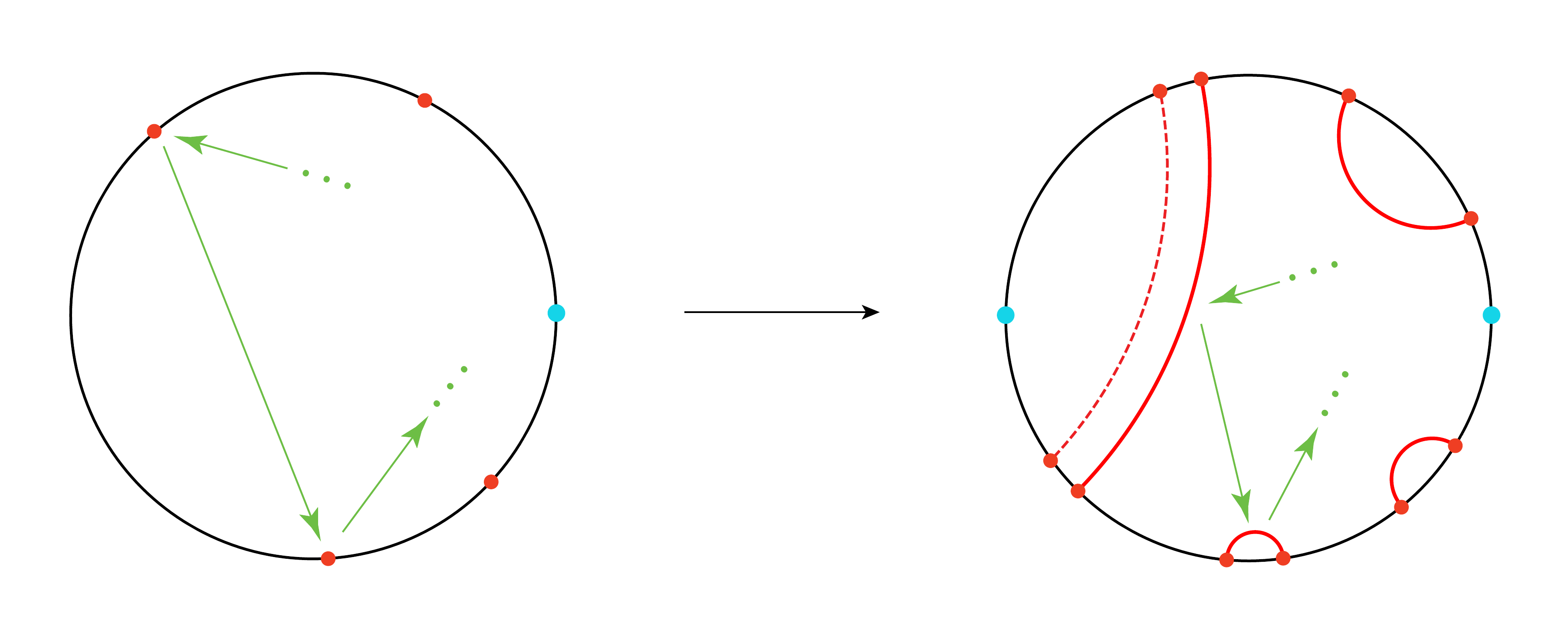}
		
		\put (36.7,19) {\footnotesize $0$}
		\put (96.3,19) {\footnotesize $0$}
		\put (61.5,19) {\footnotesize $\frac{1}{2}$}
		\put (27,34.5) {\footnotesize $t_1$}
		\put (7.3,32) {\footnotesize $t_k$}
		\put (20,1.3) {\footnotesize $t_{\sigma(k)}$}
		\put (31,6) {\footnotesize $t_q$}
		\put (95,26) {\footnotesize $x_1$}
		\put (86,35) {\footnotesize $y_1$}
		\put (75.5,36) {\footnotesize $x_k$}
		\put (66.5,6.5) {\footnotesize $y_k$}
     \put (64,9) {\footnotesize $x'_k$}
		\put (71,35) {\footnotesize $y'_k$}
		\put (74,1.8) {\footnotesize $x_{\sigma(k)}$}
		\put (81.3,1.8) {\footnotesize $y_{\sigma(k)}$}
		\put (89,5) {\footnotesize $x_q$}
		\put (94,10) {\footnotesize $y_q$}
      \put (20,18) {\footnotesize \color{uibgreen} $\m_2$}
      \put (79.5,16) {\footnotesize \color{uibgreen} $\m_3$}		
    \end{overpic}
\caption{\small Basic $\m_3$-dynamics of the intervals defined by the simulating pair $\OO_{k-1},\OO_k$.}  
\label{basic}
\end{figure}

When $t_1=0$ is the fixed point of $\m_2$, the simulating orbits and angles are  
$$
\OO_1= \{ 0 \}, \ \OO_0= \Big\{ \frac{1}{2} \Big\} \qquad \text{and} \qquad x_1=0, \ y_1=\frac{1}{2},
$$
so we may only focus on the case where $\{ t_1, \ldots, t_q \}$ has period $q>1$. Assuming \thmref{interlace} for a moment, we see that the closed intervals $I_i:=[x_i,y_i]$ for $1 \leq i \leq q$ are pairwise disjoint and none of them contains the fixed point $0$. Under $\m_3$ each $I_i$ maps homeomorphically onto $I_{\sigma(i)}$ except for $I_k$ which  contains the fixed point $1/2$ and maps to the whole circle. More precisely, consider the points  
$$
x'_k:=x_k+\frac{1}{3} \qquad \text{and} \qquad y'_k:=y_k-\frac{1}{3}
$$
which appear in the order $x_k<y'_k<1/2<x'_k<y_k$. The two sub-intervals $[x_k,y'_k]$ and $[x'_k,y_k]$ of $I_k$ map homeomorphically to $I_{\sigma(k)}$ while the central sub-interval $[y'_k,x'_k]$ maps homeomorphically onto $\ov{(\R/\Z) \sm I_{\sigma(k)}}$ (see \figref{basic}). \vs

%%%%%%%%%%%%%%%%%%%%%%%%%%%%%%%%%%%%%%%%%%%
%%%%%%%%%%%%%%%%%%%%%%%%%%%%%%%%%%%%%%%%%%% 		
		
It follows from the basic dynamical description above that 
$$
|I_k| = 3^{q-1} |I_{\sigma(k)}| \qquad \text{and} \qquad 3 |I_k| = |I_{\sigma(k)}|+1.  
$$	
This shows $|I_{\sigma(k)}|=1/(3^q-1)$, hence   
\begin{equation}\label{|I|}
|I_{\sigma^i(k)}| = \frac{3^{i-1}}{3^q-1} \qquad \text{for} \ 1 \leq i \leq q.
\end{equation}      
In particular, the shortest interval among $I_1, \ldots, I_q$ is $I_{\sigma(k)}$ and the longest is $I_k$ with 
$$
\frac{1}{3} < |I_k| = \frac{3^{q-1}}{3^q-1} \leq \frac{3}{8}.  
$$
Notice that the longest interval $I_k$ appears in the same cyclic position in the collection $\{ I_1, \ldots, I_q \}$ as the point $t_k$ does in its $\m_2$-orbit $\{ t_1, \ldots, t_q \}$. \vs

When representing orbits of rational angles of the same denominator, we often use the following convenient notation:   
$$
\{ p_1, \ldots, p_n \}/p := \left\{ \frac{p_1}{q}, \ldots, \frac{p_n}{q} \right\}.
$$

\begin{figure}[t]
\centering
\begin{overpic}[width=0.6\textwidth]{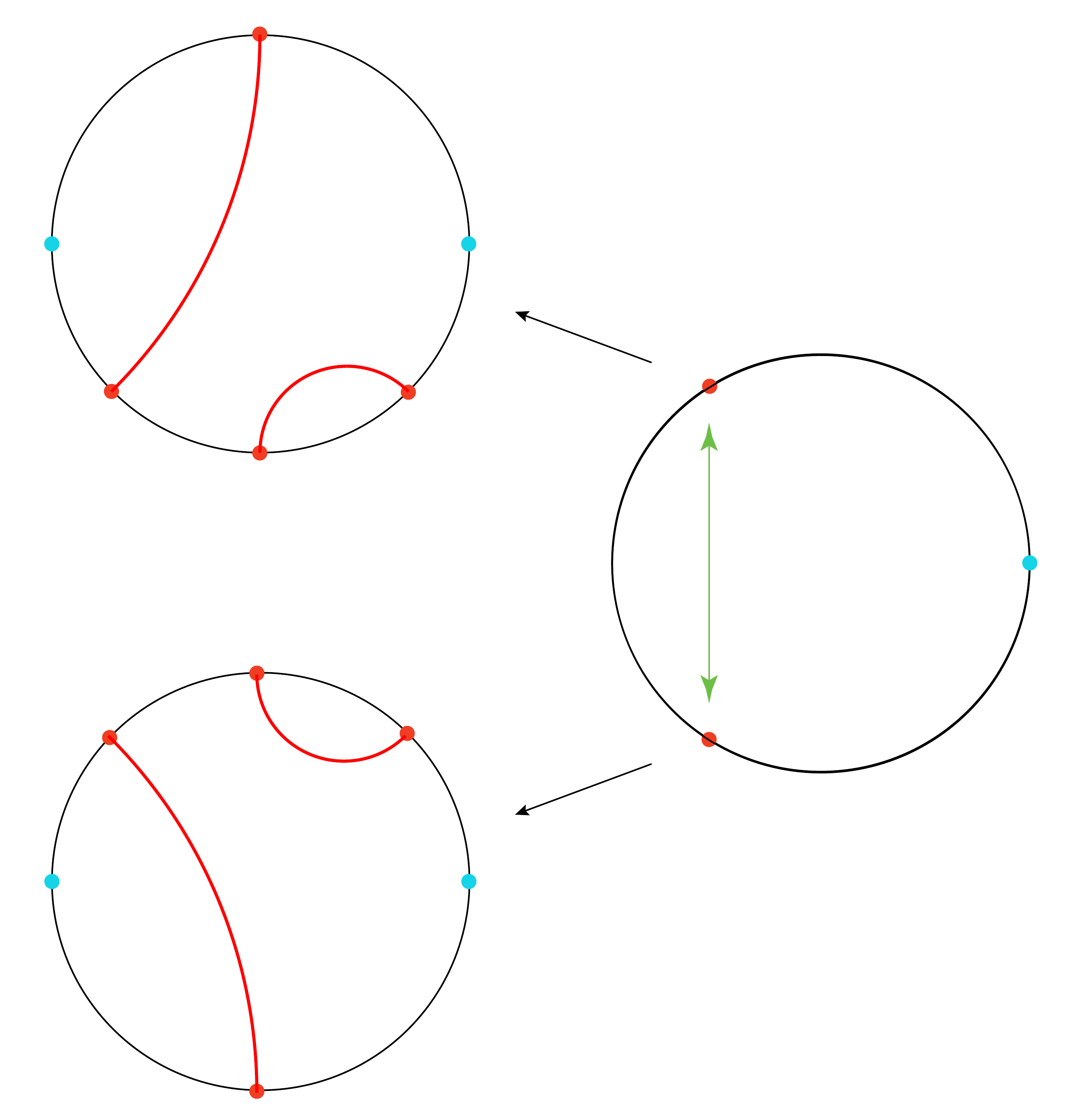}
\put (23,77) {\small $\OO_1 \cup \OO_0$}
\put (23,20) {\small $\OO_2 \cup \OO_1$}

\put (93.5,49) {\tiny $0$}
\put (64.5,63) {\tiny $\tfrac{1}{3}$}
\put (64.5,35) {\tiny $\tfrac{2}{3}$}

\put (18.5,98.5) {\tiny $x_1=2$}
\put (2,62) {\tiny $y_1=5$}
\put (22.5,56) {\tiny $6$}
\put (37.5,62) {\tiny $7$}

\put (37.5,35) {\tiny $1$}
\put (22,41.5) {\tiny $2$}
\put (1,35) {\tiny $x_2=3$}
\put (19,-0.5) {\tiny $y_2=6$}

\end{overpic}
\caption{\small Illustration of \exaref{1/3}. The $\m_2$-orbit $\{ 1, 2 \}/3$ and its two pairs of simulating $\m_3$-orbits, all with the degree $1$ combinatorics $\sigma=(12) \in S_2$ (angles are shown in multiplies of $1/8$).}  
\label{one}
\end{figure}

\EX \label{1/3}
Consider the period $2$ orbit $\{ 1, 2 \}/3$ under $\m_2$ having the combinatorics $\sigma=(12) \in S_2$ with $\deg(\sigma)=1$. This is a rotation cycle of rotation number $1/2$. The three orbits under $\m_3$ that realize $\sigma$ are
$$
\OO_0 = \{ 5,7 \}/8 \qquad   
\OO_1 = \{ 2,6 \}/8 \qquad  
\OO_2 = \{ 1,3 \}/8.  
$$
The resulting simulating orbits and angles are shown in \figref{one}.  
\ENDEX  

%%%%%%%%%%%%%%%%%%%%%%%%%%%%%%%%%%%%%%%%%%%
%%%%%%%%%%%%%%%%%%%%%%%%%%%%%%%%%%%%%%%%%%% 
   
\begin{figure}[t]
\centering
\begin{overpic}[width=0.9\textwidth]{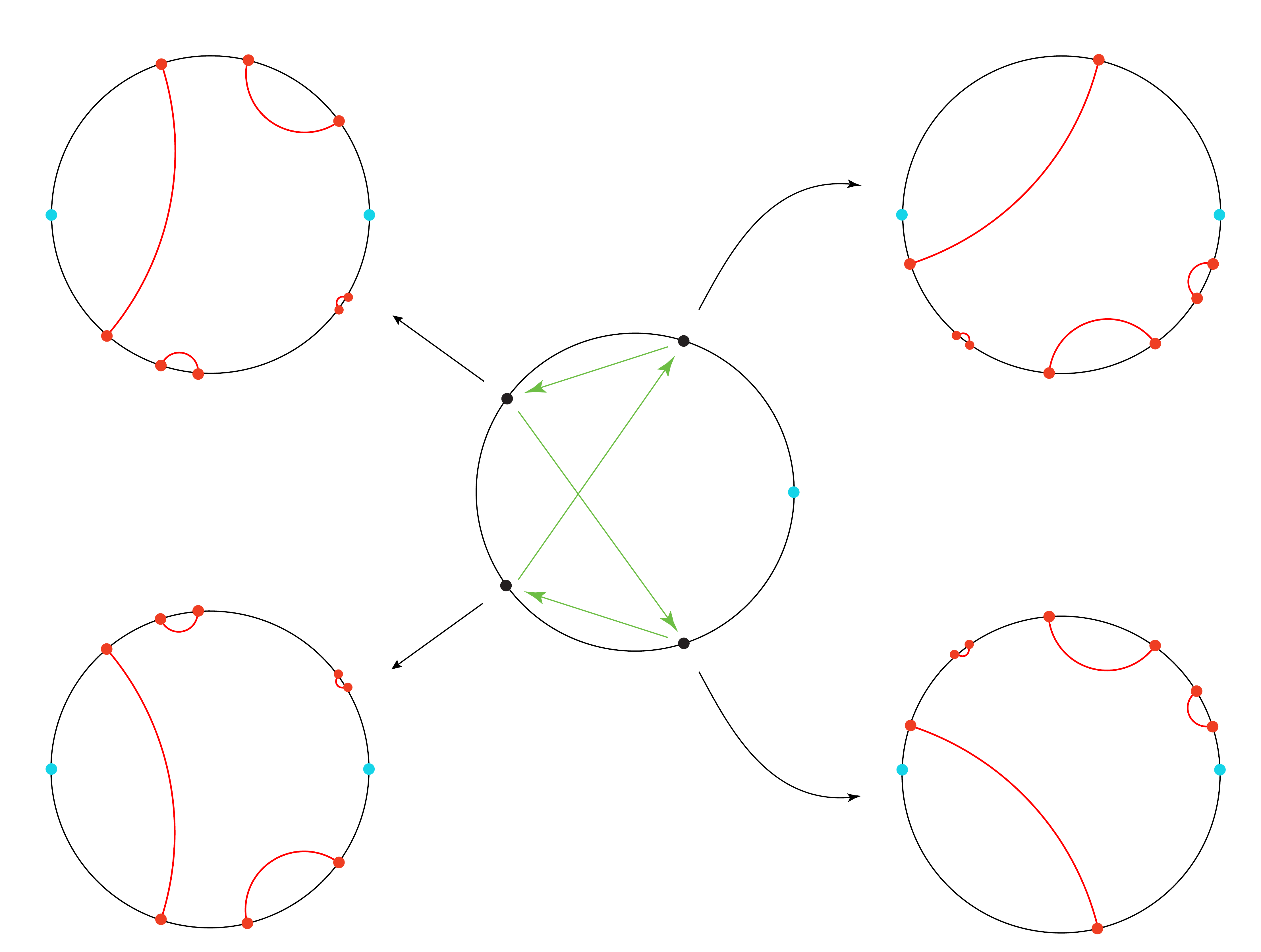}
	\put (83,57) {\footnotesize $\OO_1 \cup \OO_0$}
	\put (16,57) {\footnotesize $\OO_2 \cup \OO_1$}
	\put (16,13) {\footnotesize $\OO_3 \cup \OO_2$}
	\put (83,13) {\footnotesize $\OO_4 \cup \OO_3$}
	
	\put (53.7,45.3) {\tiny $\frac{1}{5}$}
	\put (39.5,40.5) {\tiny $\frac{2}{5}$}
	\put (39.5,31.2) {\tiny $\frac{3}{5}$}
	\put (53.7,26.3) {\tiny $\frac{4}{5}$}
	
	\put (84,72) {\tiny $x_1=17$}
	\put (64.8,53.5) {\tiny $y_1=44$}
	\put (73,47) {\tiny $51$}
	\put (75.5,45.3) {\tiny $52$}
	\put (82,43) {\tiny $59$}
	\put (91,45.5) {\tiny $68$}
	\put (95,49) {\tiny $73$}
	\put (97,54) {\tiny $76$}
	
	\put (97.2,17.5) {\tiny $4$}
	\put (95.1,21.8) {\tiny $7$}
	\put (90.2,25.7) {\tiny $12$}
	\put (82,27.8) {\tiny $21$}
	\put (75.1,26) {\tiny $28$}
	\put (72,23.6) {\tiny $29$}
	\put (65,17.5) {\tiny $x_4=36$}
	\put (83,-0.8) {\tiny $y_4=63$}
	
	\put (27.5,66) {\tiny $8$}
	\put (19,72) {\tiny $17$}
	\put (7,72) {\tiny $x_2=24$}
	\put (1.8,46.4) {\tiny $y_2=51$}
	\put (11,43.3) {\tiny $56$}
	\put (14.7,42.7) {\tiny $59$}
	\put (26,47.9) {\tiny $72$}
	\put (28.5,51.5) {\tiny $73$}

   \put (28.1,20.4) {\tiny $7$}	
	\put (26.5,23) {\tiny $8$}
	\put (14.7,28) {\tiny $21$}
	\put (10.7,27.6) {\tiny $24$}
	\put (2,25) {\tiny $x_3=29$}
	\put (7,-0.3) {\tiny $y_3=56$}
	\put (18.2,-0.6) {\tiny $63$}
	\put (28,6) {\tiny $72$}
\end{overpic}
\caption{\small Illustration of \exaref{2/5}. The $\m_2$-orbit $\{ 1, 2, 3, 4 \}/5$ and its four pairs of simulating $\m_3$-orbits, all with the degree $2$ combinatorics $\sigma=(1243)$ (angles are shown in multiplies of $1/80$).}  
\label{two}
\end{figure}
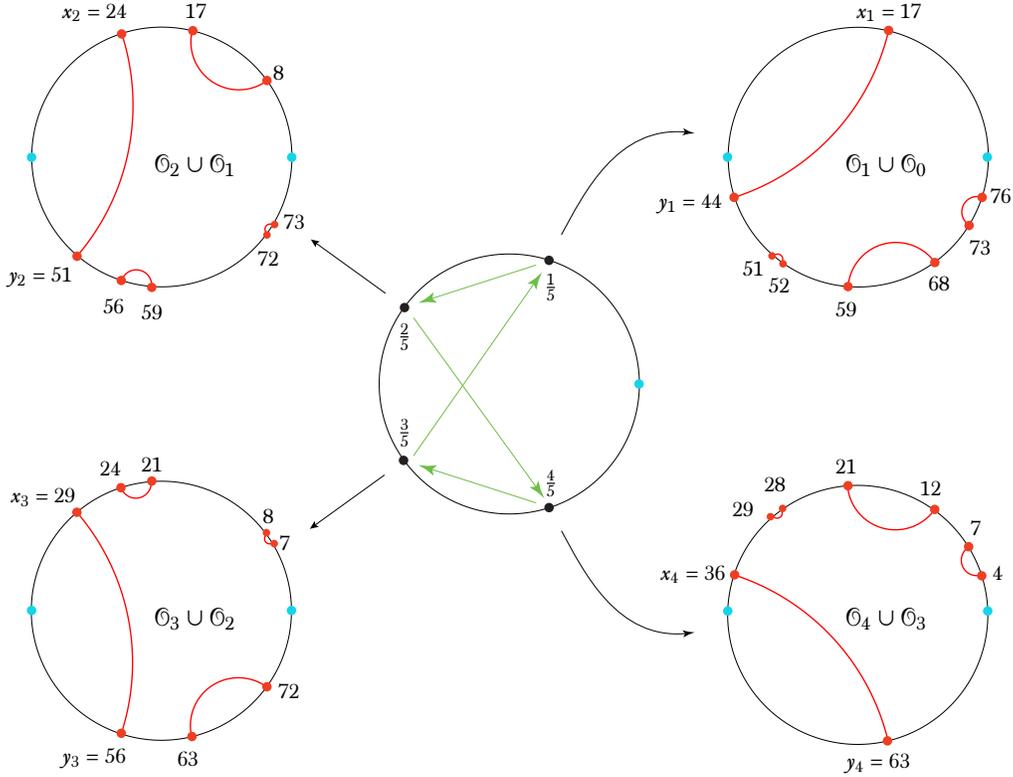

\EX \label{2/5}
Now consider the period $4$ orbit $\{ 1, 2, 3, 4 \}/5$ under $\m_2$ having the combinatorics $\sigma=(1243) \in S_4$ with $\deg(\sigma)=2$. The five orbits under $\m_3$ that realize $\sigma$ are
\begin{align*}
\OO_0 & = \{ 44, 52, 68, 76 \}/80 &   
\OO_1 & = \{ 17, 51, 59, 73 \}/80 &
\OO_2 & = \{ 8, 24, 56, 72 \}/80 \\[2pt]
\OO_3 & = \{ 7, 21, 29, 63 \}/80 & 
\OO_4 & = \{ 4, 12, 28, 36 \}/80. &   
\end{align*}
The resulting simulating orbits and angles are shown in \figref      {two}  
\ENDEX 

%%%%%%%%%%%%%%%%%%%%%%%%%%%%%%%%%%%%%%%%%%% 

\PROOF[Proof of \thmref{interlace}]
Consider the intervals $T_i:=[t_i,t_{i+1}]$ for $1 \leq i \leq q$. Since $\m_2(T_i)=T_{\sigma(i)} \cup \cdots \cup T_{\sigma(i+1)-1}$, we have 
\begin{equation}\label{double-length}
|T_{\sigma(i)}|+\cdots+|T_{\sigma(i+1)-1}|=2|T_i| \qquad \text{for} \ 1 \leq i \leq q.
\end{equation} 
There are two possibilities for how $\m_2$ acts on the interval $T_q$ containing the fixed point $0$: Either $\deg(\sigma)=1$ and $\sigma(q)<\sigma(1)$ in which case $\m_2$ maps $T_q$ over the whole circle, or $\deg(\sigma)=2$ and $\sigma(q)>\sigma(1)$ in which case $\m_2$ maps $T_q$ homeomorphically onto $T_{\sigma(q)} \cup \cdots \cup T_{\sigma(1)-1}$. \vs

Now fix $1 \leq k \leq q$ and take the points $\alpha_1, \ldots, \alpha_q$ and $\beta_1, \ldots, \beta_q$ on the circle such that 
$0<\alpha_1<\beta_1<\cdots<\alpha_q<\beta_q<1$. The intervals $A_i:=[\alpha_i,\beta_i]$ and $B_i:=[\beta_i,\alpha_{i+1}]$ for $1 \leq i \leq q$ form a partition of the circle modulo their endpoints. We may choose these points so that $|A_{\sigma^i(k)}|=3^{i-1}/(3^q-1)$ and $|B_i|=|T_i|/2$ for all $1 \leq i \leq q$. This is possible because 
$$
\sum_{i=1}^q \frac{3^{i-1}}{3^q-1} = \sum_{i=1}^q \frac{|T_i|}{2}= \frac{1}{2}. 
$$ 
Define a covering map $f: \R/\Z \to \R/\Z$ which sends the $A_i$ and $B_i$ affinely to their images as follows: 
$$
f(A_i) = \begin{cases} A_{\sigma(i)} & \quad \text{if} \ i \neq k\\
A_{\sigma(i)} \ \text{plus one full turn}\footnotemark & \quad \text{if} \ i=k  \end{cases}
\footnotetext{Given intervals $[x,y]$ and $[x',y']$ on the circle with $0<y-x<1$ and $0<y'-x'<1$, sending $[x,y]$ affinely to $[x',y']$ {\it plus one full turn} means that at the level of a lift we map $[x,y]$ affinely to $[x',y'+1]$.}
$$
and
$$
f(B_i) = \begin{cases} \ [\beta_{\sigma(i)},\alpha_{\sigma(i+1)}] & \quad \text{if} \ i \neq q \ \text{or if} \ i=q \ \text{and} \ \sigma(q)>\sigma(1) \\
\ [\beta_{\sigma(i)},\alpha_{\sigma(i+1)}] \ \text{plus one full turn} & \quad \text{if} \ i=q \ \text{and} \ \sigma(q)<\sigma(1). \end{cases}
$$
It is easy to see that $f$ is a degree $3$ piecewise affine covering map with $f'|_{A_i}=3$ for all $i$. Moreover, $f(B_i)$ contains the union $B_{\sigma(i)} \cup \cdots \cup B_{\sigma(i+1)-1}$, so by \eqref{double-length}, 
$$
f'|_{B_i} \geq \frac{|B_{\sigma(i)}|+\cdots+|B_{\sigma(i+1)-1}|}{|B_i|} = \frac{|T_{\sigma(i)}|+\cdots+|T_{\sigma(i+1)-1}|}{|T_i|}=2. 
$$
Thus, $f$ is uniformly expanding with slope $\geq 2$. It is well known that this implies the existence of a unique orientation-preserving homeomorphism $\varphi: \R/\Z \to \R/\Z$ which fixes $0$ and satisfies $\varphi \circ f = \m_3 \circ \varphi$ (see e.g. \cite[Theorem 2.4.6]{KH}). Set $\hat{x}_i :=\varphi(\alpha_i), \hat{y}_i:=\varphi(\beta_i)$. Evidently $\OO:= \{ \hat{x}_1, \ldots, \hat{x}_q \}$ and $\OO':= \{ \hat{y}_1, \ldots, \hat{y}_q \}$ are periodic orbits of $\m_3$ with combinatorics $\sigma$. As $[\hat{x}_k,\hat{y}_k]$ maps over itself under $\m_3$, it must contain a fixed point of $\m_3$ which can only be $1/2$ since the other fixed point $0$ is in $(\hat{y}_q,\hat{x}_1)$. This proves the analog of \eqref{xy} for the $\hat{x}_i,\hat{y}_i$. It follows in particular that 
$$
\OO \cap [0,1/2] = \{ \hat{x}_1, \ldots, \hat{x}_k \}  \quad \text{and} \quad \OO' \cap [0,1/2] = \{ \hat{y}_1, \ldots, \hat{y}_{k-1} \}.
$$
We conclude from the characterization \eqref{deploy} that $\OO=\OO_k, \OO'=\OO_{k-1}$, and therefore $\hat{x}_i=x_i$ and $\hat{y}_i=y_i$ for all $i$.    
\ENDPROOF

\COR \label{nothird}
Let $\OO_0, \ldots, \OO_q$ be the realizations of $\sigma \in S_q$ under $\m_3$ as above. Consider the marked interval $[x_k,y_k]$ of the simulating pair $\OO_k, \OO_{k-1}$ which contains the fixed point $1/2$. Then $\OO_j \cap \, ]x_k,y_k[ \, \neq \es$ whenever $j \notin \{ k-1,k \}$.   
\ENDCOR

\PROOF
In the natural cyclic order of orbits, let $u_i$ denote the $i$-th point of $\OO_i$ if $k \leq i \leq q$ and the $(i+1)$-st point of $\OO_i$ if $0 \leq i \leq k-1$. Then \thmref{interlace} and the characterization \eqref{deploy} show that
\[
x_k=u_k<u_{k+1}<\cdots<u_q<\frac{1}{2}<u_0<\cdots<u_{k-2}<u_{k-1}=y_k. \qedhere
\]
\ENDPROOF

\REM \label{projPi}
It can be shown that for each $k$ there is a degree $1$ monotone map $\Pi: \R/\Z \to \R/\Z$ which fixes $0$, collapses every $I_i=[x_i,y_i]$ to $t_i$, and satisfies the semiconjugacy relation 
$$
\Pi \circ \m_3 = \m_2 \circ \Pi  \qquad \text{outside} \ I_k.
$$
The iterated $\m_3$-preimages of $I_k$ form the plateaus of $\Pi$ and their union has full measure in $\R/\Z$.     
\ENDREM 

The $180^\circ$ rotation $\theta \mapsto \theta+1/2 \modd$ commutes with $\m_3$ and exchanges the fixed points $0,1/2$. Considering the $q+1$ realizations $\OO_0, \ldots, \OO_q$ constructed above, it follows that the rotated sets $\OO^\ast_k := \OO_k+1/2$ are also $q$-periodic orbits of $\m_3$, although they generally have a different combinatorics. In fact, if $\rho$ denotes the rotation $i \mapsto i+1 \ (\operatorname{mod} q)$, the combinatorics of $\OO^\ast_k$ is $\rho^{-k} \sigma \rho^k$, where $\sigma$ is the common combinatorics of $\OO_0, \ldots, \OO_q$. In particular, if $\deg(\sigma)=1$, then $\sigma$ commutes with $\rho$ and the $\OO^\ast_k$ have the same combinatorics as the $\OO_k$. It follows from the characterization \eqref{deploy} that 
\begin{equation}\label{invol}
\OO^\ast_k = \OO_{q-k} \qquad \text{for} \ 0 \leq k \leq q \ \text{if} \ \deg(\sigma)=1. 
\end{equation}     
If $\deg(\sigma)=2$, then \eqref{invol} holds only for $k=0$ and $k=q$. This follows from the fact that the rotational symmetry group of such $\sigma$ is trivial (see \cite{PZ1}), so $\rho^{-k} \sigma \rho^k = \rho$ if and only if $k \in \{ 0,q \}$. 

\section{The lemon family of cubic polynomials}

We begin by briefly recalling basic facts about the dynamics of cubic polynomials which we assume the reader is generally familiar with (see \cite{M1} and \cite{BH} for more details).   

\subsection{Green's function and external rays}\label{gfer}

We will be working with cubic polynomials in the normal form 
$$
P=P_{a,b}: z \mapsto z^3 +3a z^2+ b \qquad (a,b \in \C)
$$
with marked critical points at $\omega_1=0,\omega_2=-2a$. More precisely, we have \vs

\bgroup
\def\arraystretch{1.2}
\begin{center}
	\begin{tabular}{|c|c|c|}
		\hline
		critical point & $\om_1=0$ & $\om_2=-2a$ \\
		\hline 
		critical value & $b$ & $4a^3+b$ \\ 
		\hline 
		co-critical point & $-3a$ & $a$ \\
		\hline 
	\end{tabular}
\end{center}
\egroup

\vs

\noindent
(Here by the {\bit co-critical point} associated with $\omega_i$ is meant the unique point $z \neq \omega_i$ with $P(z)=P(\omega_i)$.)  We will think of $\omega_1$ and $\omega_2$ as the ``first'' and ``second'' critical points of $P$, respectively. The space $\PP(3)$ of all such critically marked cubics is isomorphic to $\C^2$ with coordinates $(a,b)$. Note that two distinct cubics $P_{a,b}$ and $P_{a',b'}$ are conjugate by an affine map preserving the critical points markings if and only if $a'=-a, b'=-b$. Thus, the space of all affine conjugacy classes (the so-called ``moduli space'') of cubics is also isomorphic to $\C^2$ with coordinates $(a^2,b^2)$. \vs

Let us fix $P \in \PP(3)$. We denote by $K=K_P$ and $J=J_P$ the filled Julia set and Julia set of $P$, respectively. Let $G=G_P: \C \to [0,+\infty[$ denote the Green's function of $K$ defined by
$$
G(z):= \lim_{n \to \infty} 3^{-n} \log^+|P^{\circ n}(z)|.
$$
It is easy to see that $G$ is subharmonic in $\C$, harmonic outside $K=G^{-1}(0)$, and satisfies 
$$
G(P(z))=3G(z) \qquad \text{for all} \ z \in \C. 
$$
There is a canonical \Bottcher coordinate $\phi=\phi_P$ which is tangent to the identity near $\infty$ and satisfies 
\begin{equation}\label{bottinfty}
\phi(P(z))=(\phi(z))^3.
\end{equation}
It extends to a biholomorphism between $\{ z: G(z) > G^{\ast} \}$ and $\{ z: |z| > \e^{G^{\ast}} \}$, where $G^{\ast}=\max_{i=1,2} G(\omega_i)$, and is related to $G$ by
$$
\log |\phi(z)|=G(z) \qquad \text{if} \ G(z)>G^{\ast}.
$$
In particular, if $K$ is connected, then $G^{\ast}=0$ and $\phi : \C \sm K \to \C \sm \ov{\D}$ is a conformal isomorphism. In this case for each $\theta \in \R/\Z$ the smooth curve $R(\theta)=R_P(\theta):= \phi^{-1}(\{ \e^{s+2\pi \ii \theta}: s>0 \})$ is called the {\bit external ray} (or the {\bit dynamic ray}) of $P$ at angle $\theta$. \vs

The external rays of $P$ can also be defined when $K$ is disconnected. In this case the equipotential curves $G=\text{const.}>0$ form a singular foliation of the basin of infinity $\C \sm K$ and the field lines of $\nabla G$ form an orthogonal singular foliation. We denote by $R(\theta)$ the maximally extended smooth field line of $\nabla G$ given by $s \mapsto \phi^{-1}(\e^{s+2\pi \ii \theta})$ for large $s$. This naturally parametrizes $R(\theta)$ by the potential, so for each $\theta$ there is an $s(\theta) \geq 0$ such that $G(R(\theta,s))=s$ for all $s>s(\theta)$. The set $N=N_P \subset \R/\Z$ of angles $\theta$ for which $s(\theta)>0$ is countable, dense and backward-invariant under $\m_3$. For $\theta \notin N$ the field line $R(\theta)$ extends all the way to the Julia set, just like the ordinary external rays in the connected case. For $\theta \in N$ the field line $R(\theta)$ crashes into a singularity of $\nabla G$ at the potential $s(\theta)>0$ and there is more than one way to extend it to a curve consisting of field lines and singularities on which $G$ defines a homeomorphism onto $]0,+\infty[$. But there are always two special extensions: $R^+(\theta)=\lim_{\alpha \downarrow \theta} R(\alpha)$ which turns immediate right and $R^-(\theta)=\lim_{\alpha \uparrow \theta} R(\alpha)$ which turns immediate left at each singularity they meet (here the angle $\alpha$ tends to $\theta$ avoiding the exceptional set $N$). These are called the {\bit broken external rays} at angle $\theta$. We refer to \cite{PZ2} for a comprehensive treatment of the external rays of disconnected Julia sets. \vs

As a concrete example of the disconnected case, consider the situation where $\om_1 \in K$ but $\om_2 \notin K$. Then there is a unique angle $\theta_0 \in \R/\Z$ such that the rays $R(\theta_0 \pm 1/3)$ crash into $\om_2$ (compare \figref{onehalfdisc} left). In this case, the exceptional set $N$ will consist of all angles $\theta$ for which $\m_3^{\circ n}(\theta) = \theta_0 \pm 1/3$ for some $n \geq 0$. \vs

For $z_0 \in J=\bd K$ we denote by $A_{z_0}$ the set of angles $\theta \in \R/\Z$ for which $R(\theta)$ (if $\theta \notin N)$ or one of $R^\pm(\theta)$ (if $\theta \in N$) lands at $z_0$. Every periodic $\theta$ belongs to $A_{z_0}$ for some repelling or parabolic point $z_0$ whose period under $P$ divides that of $\theta$ under $\m_3$. Conversely, if $z_0$ is a repelling or parabolic point of period $k$, then $A_{z_0}$ is non-empty and $\m_3^{\circ k}: A_{z_0} \to A_{z_0}$ has a well-defined rotation number $\rho$. Moreover, the following dichotomy holds: Either $\rho$ is rational in which case $A_{z_0}$ is a union of finitely many cycles of the same length, or $\rho$ is irrational in which case $A_{z_0}$ is a minimal Cantor set. The first alternative always happens when $K$ is connected or more generally when the connected component of $K$ containing $z_0$ is not a single point (see \cite{DH1} or \cite{M1} for the connected case and \cite{LP} or \cite{PZ2} for the disconnected case).  

\subsection{The lemon family: dynamical plane}\label{lfdyn}
 
Consider the complex line $\per_1(0) \subset \PP(3)$ consisting of all cubics $P_{a,b}$ for which the critical point $\om_1=0$ is fixed, that is, those in the normal form 
$$
P_a:=P_{a,0}: z \mapsto z^3+3az^2 \qquad (a \in \C)
$$
(for simplicity we use the subscript $a$ instead of $P_a$ for the objects associated with the cubics in $\per_1(0)$). The connectedness locus of this family is shown in \figref{lemonfig}. The prominent central region $\HH_0$ in this figure, resembling the shape of a lemon, consists of all $a \in \C$ for which the critical point $\om_2=-2a$ belongs to the immediate basin of attraction of $0$ (equivalently, the Julia set $J_a$ is a Jordan curve). We refer to $\{ P_a \}_{a \in \C}$ as the {\bit lemon family} of cubics. \vs
  
The study of the slice $\per_1(0)$ appears in the works of several authors including Faught \cite{F}, Milnor \cite{M4}, and Roesch \cite{R1, R2}. Our brief description in the dynamical plane is closest to \cite{M4} but differs in some details in order to incorporate the notion of simulating orbits from \S \ref{subsec:sp} into the picture. \vs

Fix some $a \neq 0$. The immediate basin of attraction of $0$ is a Jordan domain $B_a$. There is a canonical \Bottcher coordinate $\beta_a$ which is defined and holomorphic near $0$ and satisfies 
\begin{equation}\label{bottzero}
\beta_a(P_a(z))=(\beta_a(z))^2.
\end{equation}
The symmetry $P_{-a}(z)=-P_a(-z)$ implies
\begin{equation}\label{bottsym}
\beta_{-a}(z)=\beta_a(-z). 
\end{equation}
We distinguish three cases depending on the position of $\om_2=-2a$: \vs 

$\bullet$ {\bf Case 1}. $\om_2 \in B_a$. Then $P_a: B_a \to B_a$ is a degree 3 branched covering and the filled Julia set $K_a= \ov{B_a}$ is a closed topological disk. There is a largest Jordan domain in $B_a$ containing $0$ which maps conformally by $\beta_a$ to a round disk $|z|<r<1$ and contains $\om_2$ on its boundary. However, $\beta_a$ extends analytically to a neighborhood of $\om_2$ using the conjugacy relation \eqref{bottzero} (making $\om_2$ a critical point of this extension).% 
\vs

$\bullet$ {\bf Case 2}. $\om_2 \in K_a \sm B_a$. Then $\beta_a : B_a \to \D$ is a conformal isomorphism. This allows a parametrization $z_a(t):=\beta_a^{-1}(\e^{2\pi \ii t})$ of the Jordan curve $\bd B_a$ by the {\bit internal angle $t \in \R/\Z$} which satisfies
\begin{equation}\label{double}
P_a(z_a(t))=z_a(\m_2(t)).
\end{equation}
In this case $K_a$ is the union of the closed topological disk $\ov{B_a}$ together with a countable collection of disjoint continua called the {\bit (dynamical) limbs} of $K_a$. Each limb meets $\bd B_a$ at a single point called its {\bit root}. We denote by $L_a(t)$ the limb with the root $z_a(t)$. The countable set of roots is backward-invariant under the degree $2$ covering map $P_a: \bd B_a \to \bd B_a$. More precisely, there is a unique internal angle $\tau_a \in \R/\Z$ with the following properties: \vs
\begin{enumerate}[leftmargin=*]
\item[(i)]
$z_a(t)$ is a root if and only if $\m_2^{\circ n}(t) = \tau_a$ for some $n \geq 0$. \vs
\item[(ii)]
$L_a(\tau_a)$ contains $\om_2$. Under $P_a$ the limb $L_a(t)$ maps homeomorphically to $L_a(\m_2(t))$ if $t \neq \tau_a$, and to the entire filled Julia set $K_a$ if $t=\tau_a$. 
\end{enumerate}

\begin{figure}[t]
\centering
\begin{overpic}[width=0.6\textwidth]{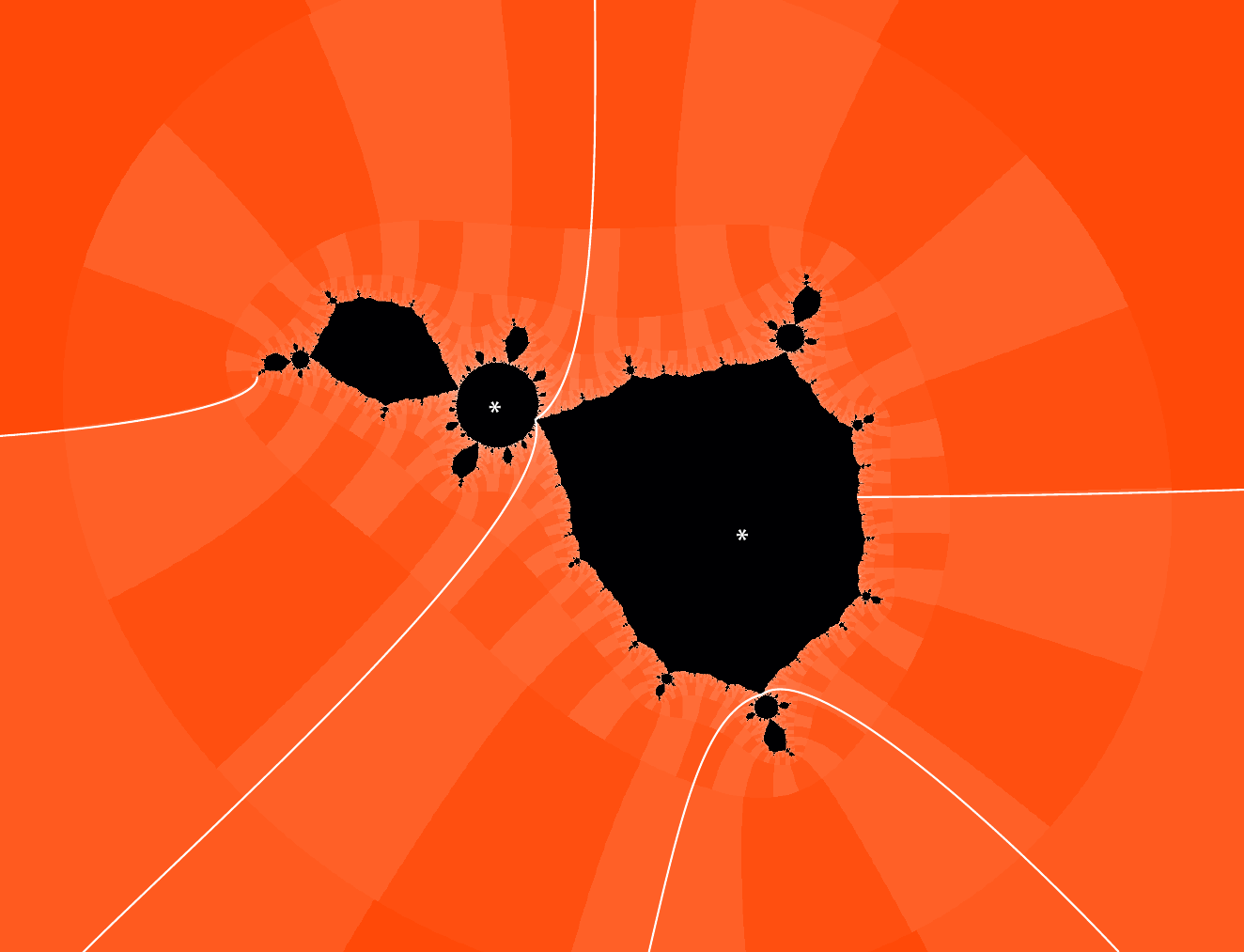}
\put (97.2,38) {\small $0$}
\put (1,43.5) {\small $\frac{1}{2}$}
\put (49,72.3) {\small $x_1=\frac{1}{4}$}
\put (11.4,1.9) {\small $y_1=\frac{5}{8}$}
\put (50,1.9) {\small $\frac{3}{4}$}
\put (90,1.9) {\small $\frac{7}{8}$}
\put (58,31) {\color{white} \tiny $\om_1$}
\put (57,40) {\color{white} \footnotesize $B_a$}
\put (39,45) {\color{white} \tiny $\om_2$}
\put (62.2,36) {\color{white} \tiny $z_a(0)$}
\put (44,41.5) {\color{white} \tiny $z_a(1/3)$}
\put (55,23) {\color{white} \tiny $z_a(2/3)$}
\put (25,65) {\small $W_a(1/3)$}
\put (65,5) {\small $W_a(2/3)$}
\end{overpic}
\caption{\small A cubic $P_a$ in the lemon family. Here $\tau_a=1/3$ has period $q=2$ under $\m_2$ with the simulating angles $x_1=1/4, y_1=5/8$. For this example $\om_2$ is periodic of period $2$. In the parameter space, $P_a$ belongs to the limb $\LL_0(1/3)$.}  
\label{onehalf}
\end{figure}
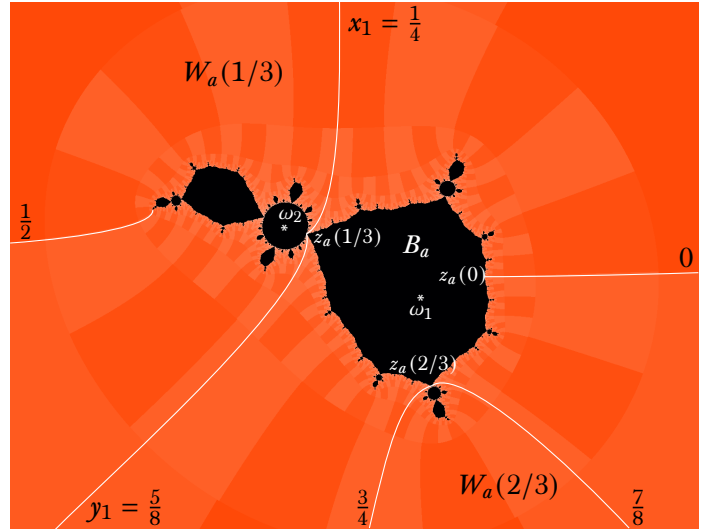

Compare \figref{onehalf}. Because of (ii), we call $L_a(\tau_a)$ the {\bit critical limb} and $\tau_a$ the {\bit critical internal angle}. Note that $\om_2$ can be the root of $L_a(\tau_a)$ (see \corref{rootcp} below). \vs  

By \cite{R1}, for each $\eps>0$ there are at most finitely many limbs $L_a(t)$ with Euclidean diameter $>\eps$. This property can be used to show that every non-root point on $\bd B_a$ is the landing point of a unique external ray while every root point on $\bd B_a$ has at least two external rays landing on it (see \cite[\S 5.2]{Z} for details of this argument in a very similar case). Below we will show that each root has precisely two external rays landing on it. \vs 

The external rays landing at a root $z_a(t)$ divide the plane into finitely many open sectors. By definition, the {\bit (dynamical) wake} $W_a(t)$ is the complement of the closure of the sector containing $B_a$. Notice that $W_a(t)$ contains $L_a(t) \sm \{ z_a(t) \}$ and distinct roots have disjoint wakes. The wake $W_a(t)$ is bounded by two rays $R_a(x),R_a(y)$ landing at $z_a(t)$ which we always label counterclockwise, so $R_a(\theta) \subset W_a(t)$ if and only if $\theta \in \, ]x,y[$. We think of the length of the interval $[x,y]$ as the {\bit angular size} of $W_a(t)$. The critical wake $W_a(\tau_a)$ can be described as the unique wake whose angular size is $\geq 1/3$. The analogue of the property (ii) above holds for wakes: Under $P_a$ the wake $W_a(t)$ maps conformally to $W_a(\m_2(t))$ if $t \neq \tau_a$, and to the entire plane if $t=\tau_a$.

\COR \label{rootcp}
If the critical internal angle $\tau_a$ is not periodic under $\m_2$, then $z_a(\tau_a)=\om_2$ and every root is the landing point of precisely two external rays. 
\ENDCOR

\PROOF
The forward $\m_2$-orbit of $\m_2(\tau_a)$ never hits $\tau_a$, hence  $z_a(\m_2(\tau_a))$ is not a root and has a unique ray landing on it. Since there are at least two rays landing at $z_a(\tau_a)$, it follows that $z_a(\tau_a)$ is a critical point. The fact that the local degree of $P_a$ at $z_a(\tau_a)$ is $2$ now shows that there are precisely two rays landing there, and therefore two rays landing at all other roots since they are all non-critical.
\ENDPROOF

Every periodic point $z_a(t) \in \bd B_a$ of period $q \geq 1$ has combinatorial rotation number zero since the internal ray of angle $t$ landing at $z_a(t)$ is invariant under $P_a^{\circ q}$. It follows in particular that only the external rays at angles $0$ or $1/2$ can land at the fixed point $z_a(0)$. After conjugating $P_a$ by $z \mapsto -z$ if necessary (that is, replacing $P_a$ with $P_{-a}$), we may assume that $P_a$ is {\bit fixed-ray normalized} in the following sense: The external ray $R_a(0)$ lands at $z_a(0)$, and if both $R_a(0),R_a(1/2)$ land at $z_a(0)$, then the wake $W_a(0)$ corresponds to the interval $]0,1/2[$ (as opposed to $]1/2,1[$). \vs

Let us now discuss the case where the critical internal angle $\tau_a$ is periodic of some period $q$ under $\m_2$. If $q=1$ then $\tau_a=0$ and the critical wake $W_a(0)$ is bounded by the rays $R_a(0), R_a(1/2)$ landing at the fixed point $z_a(0)$. Notice that $0,1/2$ are the simulating angles of $0$ in the sense of \S \ref{subsec:sp}. If $q \geq 2$ then $\tau_a \neq 0$ and by our normalization $R_a(0)$ lands at $z_a(0)$ but $R_a(1/2)$ does not. Label the $\m_2$-orbit of $\tau_a$ in positive cyclic order as $0<t_1<\cdots<t_q<1$, so $t_k=\tau_a$ for some $1 \leq k \leq q$. The points $z_i:=z_a(t_i) \in \bd B_a$ form a repelling or parabolic orbit of period $q$. Let $A_i$ denote the finite set of angles of external rays that land on $z_i$. As the $z_i$ have combinatorial rotation number zero under $P_a^{\circ q}$, the iterate $\m_3^{\circ q}$ acts as the identity map on each $A_i$. It follows that all orbits in $A_1 \cup \cdots \cup A_q$ have length $q$ with the same combinatorics under $\m_3$ as that of $\tau_a$ under $\m_2$. For each $1 \leq i \leq q$, there are angles $0<x_i<y_i<1$ in $A_i$ so that the wake $W_a(t_i)$ is bounded by $R_a(x_i), R_a(y_i)$. Evidently $\{ x_1, \ldots, x_q \}$ and $\{ y_1, \ldots, y_q \}$ are distinct orbits in $A_1 \cup \cdots \cup A_q$ with $0<x_1<y_1<\cdots<x_q<y_q<1$. Since the critical wake $W_a(\tau_a)$ maps over itself, it must contain the fixed ray $R_a(1/2)$. This implies $x_k<1/2<y_k$ by the choice of $k$. The discussion of \S \ref{subsec:sp} now shows that $\{ x_1, \ldots, x_q \}=\OO_k$ and $\{ y_1, \ldots, y_q \}=\OO_{k-1}$ are the simulating orbits for $\tau_a$, with the simulating angles $x_k, y_k$. It also follows from uniqueness of simulating orbits that no other orbit can show up in the union $A_1 \cup \cdots \cup A_q$. Thus, $R_a(x_i),R_a(y_i)$ are the only external rays landing at $z_i$, that is, $A_i=\{ x_i,y_i \}$.  \vs

\figref{onehalf} illustrates the filled Julia set and dynamical wakes for a cubic $P_a$ with $\tau_a=1/3$ and the simulating angles $x_1=1/4, y_1=5/8$. \vs

%%%%%%%%%%%%%%%%%%%%%%%%%%%%%%%%%%%%%%%%%%%
%%%%%%%%%%%%%%%%%%%%%%%%%%%%%%%%%%%%%%%%%%%

$\bullet$ {\bf Case 3}. $\om_2 \notin K_a$. Then $P_a$ restricted to a neighborhood of $B_a$ is a quadratic-like map hybrid equivalent to $Q_0:z \mapsto z^2$. The filled Julia set $K_a$ has countably many components homeomorphic to a closed disk (these are the iterated preimages of $B_a$) and uncountably many point components. \vs 

In fact, much of the wake structure described in Case 2 remains intact in the disconnected case. As the \Bottcher map $\beta_a: B_a \to \D$ is a conformal isomorphism, the parametrization $t \mapsto z_a(t)$ of $\bd B_a$ still exists. Let $\Theta_a$ be the set of external angles $\theta$ such that the (smooth or broken) ray $R_a(\theta)$ lands on $\bd B_a$. According to \cite{PZ2}, $\Theta_a$ is an $\m_3$-invariant Cantor set and the landing point map $\zeta_a: \Theta_a \to \bd B_a$ is monotone and surjective. Thus, every point on $\bd B_a$ is the landing point of one or two external rays. The connected components of $(\R/\Z) \sm \Theta_a$, called {\bit gaps}, define the dynamical wakes. More precisely, for each gap $]\theta_1, \theta_2[$ we set $t:=z_a^{-1}(\zeta_a(\theta_1))=z_a^{-1}(\zeta_a(\theta_2))$ and define the wake $W_a(t)$ to be the connected component of $\C \sm (R_a(\theta_1) \cup R_a(\theta_2) \cup \{ z_a(t) \})$ that does not contain $B_a$. The critical internal angle $\tau_a$ is then characterized by the condition that $W_a(\tau_a)$ has angular size $\geq 1/3$, or equivalently, contains the escaping critical point $\om_2$ in its closure. \vs

\begin{figure}[t]
\centering
\begin{overpic}[width=\textwidth]{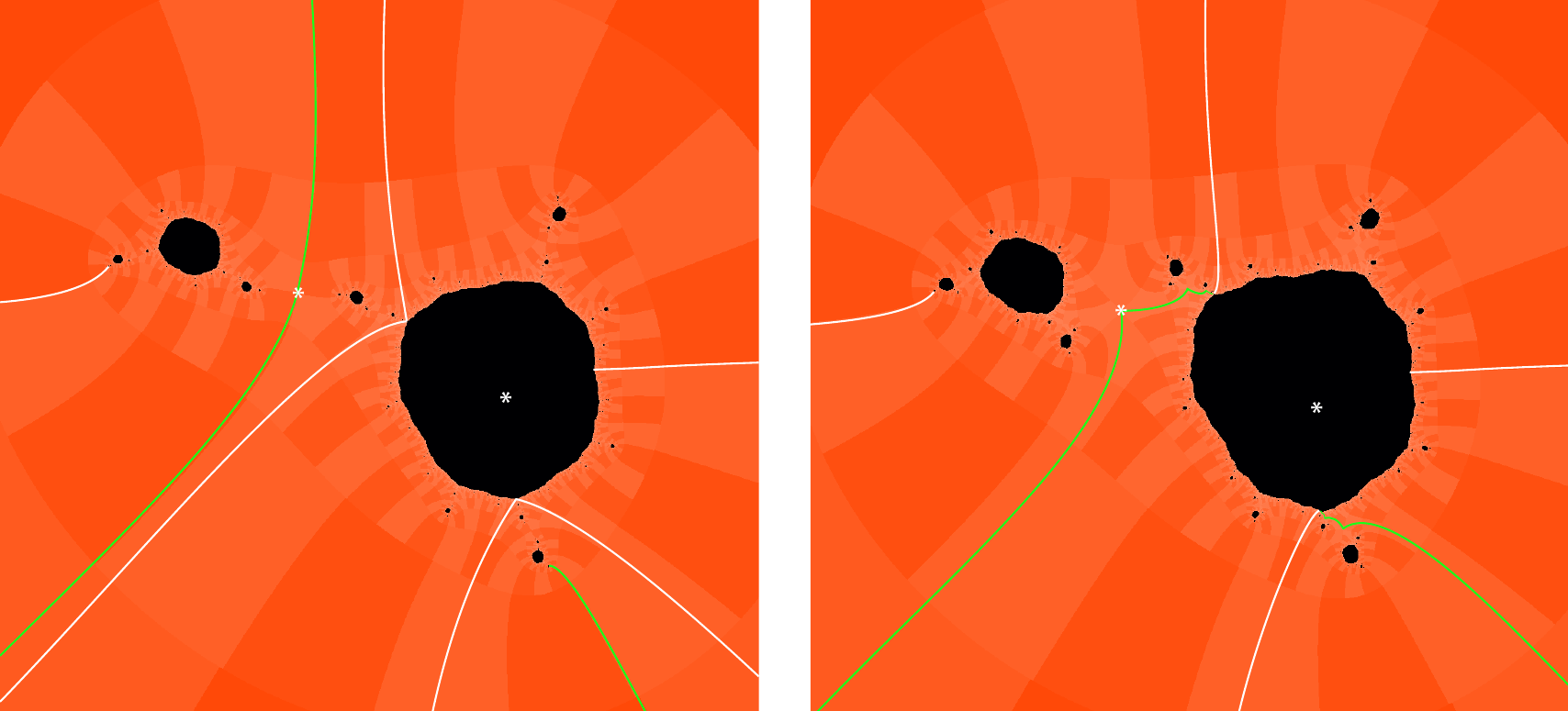}
\put (46.5,23) {\footnotesize $0$}
\put (98.2,22.7) {\footnotesize $0$}
\put (0.5,27.8) {\footnotesize $\frac{1}{2}$}
\put (52,26.5) {\footnotesize $\frac{1}{2}$}
\put (25,43) {\footnotesize $x_1=\frac{1}{4}$}
\put (77.5,43) {\footnotesize $x_1=\frac{1}{4}$}
\put (3,1.5) {\footnotesize $y_1=\frac{5}{8}$}
\put (56,1.5) {\footnotesize $y_1=\frac{5}{8}$}
\put (26,1.5) {\footnotesize $\frac{3}{4}$}
\put (77.5,1.5) {\footnotesize $\frac{3}{4}$}
\put (46.5,4.8) {\footnotesize $\frac{7}{8}$}
\put (98.2,4.5) {\footnotesize $\frac{7}{8}$}
\put (31,18.2) {\color{white} \tiny $\om_1$}
\put (82.8,17.6) {\color{white}\tiny $\om_1$}
\put (30,23) {\color{white} \footnotesize $B_a$}
\put (82,23) {\color{white} \footnotesize $B_a$}
\put (16.9,27.5) {\tiny $\om_2$}
\put (70,26.5) {\tiny $\om_2$}
\put (14,43) {\footnotesize $\xi+\frac{1}{3}$}
\put (0.5,10.5) {\footnotesize $\xi+\frac{2}{3}$}
\put (36.8,1.5) {\footnotesize $3\xi$}
\end{overpic}
\caption{\small Examples of cubics $P_a$ in the lemon family with disconnected Julia sets. Here $\tau_a=1/3$ has the simulating angles $x_1=1/4, y_1=5/8$ as in \figref{onehalf}. Left: $P_a$ is on the parameter ray $\RR(\xi) \subset \WW_0(1/3)$ for some $11/12<\xi<23/24$ and each dynamical wake is bounded by a pair of smooth rays. Right: $P_a \in \RR(23/24) \subset \bd\WW_0(1/3)$ and each dynamical wake is bounded by one smooth and one infinitely broken ray.}  
\label{onehalfdisc}
\end{figure}

\figref{onehalfdisc} shows the dynamical wakes for two cubics $P_a$ with disconnected Julia sets and $\tau_a=1/3$. \vs

%%%%%%%%%%%%%%%%%%%%%%%%%%%%%%%%%%%%%%%%%%%
%%%%%%%%%%%%%%%%%%%%%%%%%%%%%%%%%%%%%%%%%%%

Based on the preceding discussion, it is not hard to verify the following  
 
\THM
Let $t \in \R/\Z$ be periodic under $\m_2$ with simulating angles $x_k,y_k$. Suppose $P_a$ is fixed-ray normalized and $\om_2 \notin B_a$. Then the (smooth or broken) external rays $R_a(x_k), R_a(y_k)$ co-land if and only if $\tau_a=t$.   
\ENDTHM

\subsection{The lemon family: parameter plane}\label{lfpar}

\begin{figure}[t]
\centering
\begin{overpic}[width=0.75\textwidth]{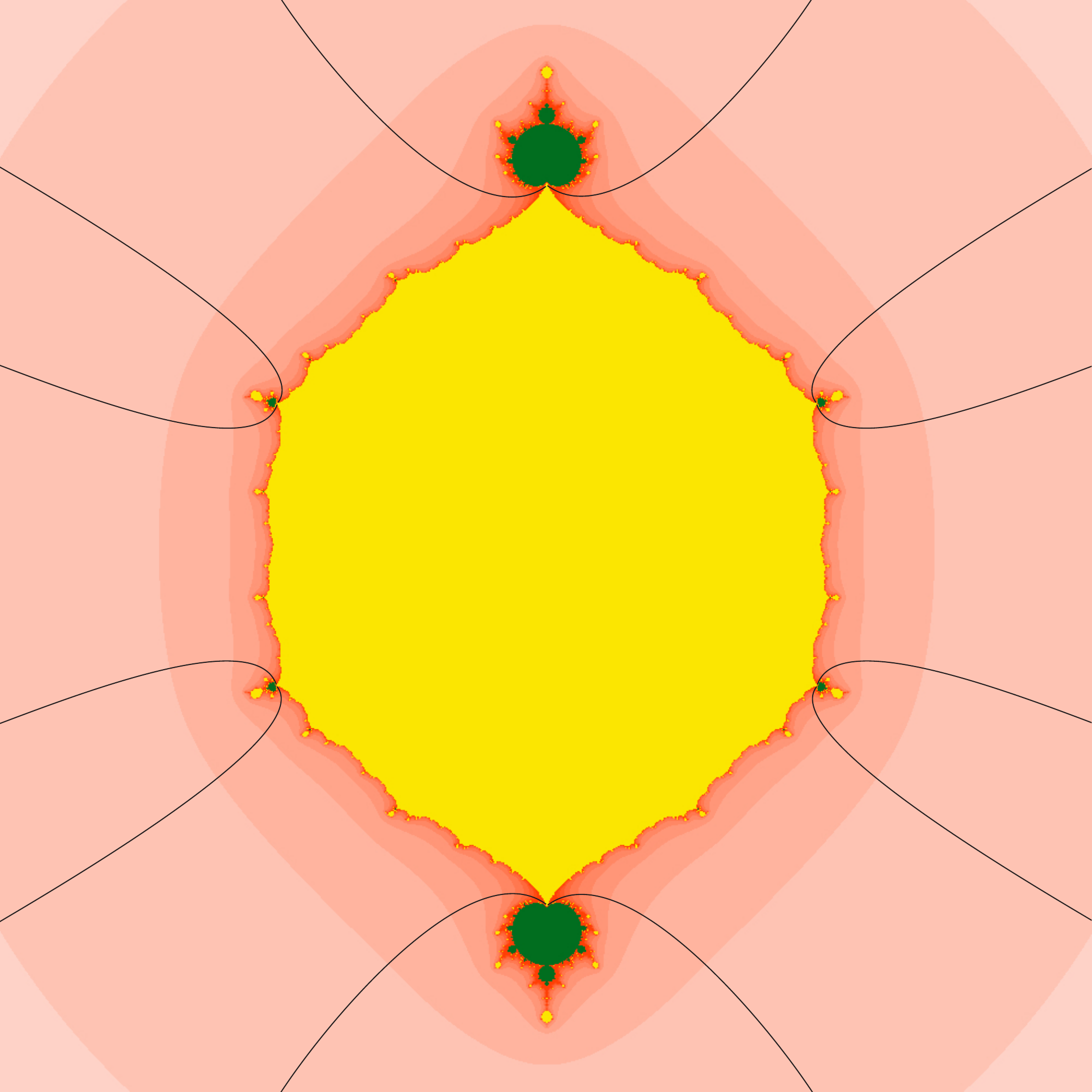}
\put (48,48) {$\HH_0$}
\put (55,5) {\small $\WW_0(0)$}
\put (82,30) {\small $\WW_0(1/3)$}
\put (82,68) {\small $\WW_0(2/3)$}
\put (74,96) {\footnotesize $\frac{1}{6}$}
\put (24,96) {\footnotesize $\frac{1}{3}$}
\put (24,3) {\footnotesize $\frac{2}{3}$}
\put (74,3) {\footnotesize $\frac{5}{6}$}
\put (96,86) {\footnotesize $\frac{1}{12}$}
\put (96,61) {\footnotesize $\frac{1}{24}$}
\put (96,38) {\footnotesize $\frac{23}{24}$}
\put (96,12) {\footnotesize $\frac{11}{12}$}
\put (1,12) {\footnotesize $\frac{7}{12}$}
\put (1,38) {\footnotesize $\frac{13}{24}$}
\put (1,61) {\footnotesize $\frac{11}{24}$}
\put (1,86) {\footnotesize $\frac{5}{12}$}
\end{overpic}
\caption{\small The connectedness locus $\CC_0=\CC(3) \cap \per_1(0)$ of the lemon family and selected wakes $\WW_0(t)$. Each intersection $\ov{\WW_0(t)} \cap \CC_0$ is the limb $\LL_0(t)$.}  
\label{lemonfig}
\end{figure}

\begin{figure}[t]
\centering
\begin{overpic}[width=0.5\textwidth]{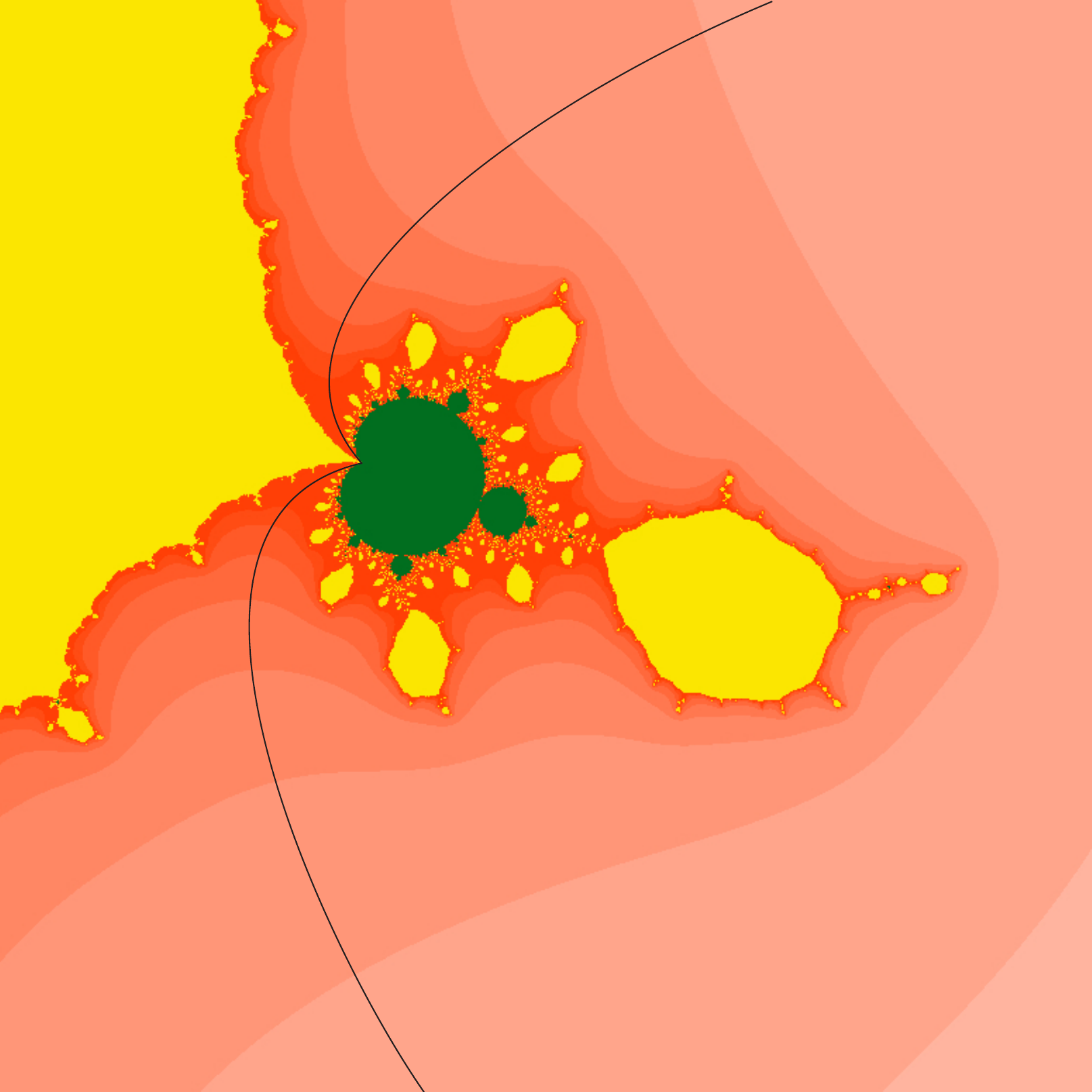}
\put (8,78) {$\HH_0$}
\put (70,10) {$\WW_0(1/3)$}
\put (32.5,55.5) {\scalebox{0.6}{\color{white} $\MM_0(1/3)$}}
\put (70,94) {$\frac{23}{24}$}
\put (39,3) {$\frac{11}{12}$}
\put (19,58.5) {\tiny $a(1/3)$}
\end{overpic}
\caption{\small Details of the limb $\LL_0(1/3)$, with its main Mandelbrot copy $\MM_0(1/3)$ (in green) and various capture components (in yellow) attached to its dyadic tips.}  
\label{onethirdlimb}
\end{figure}
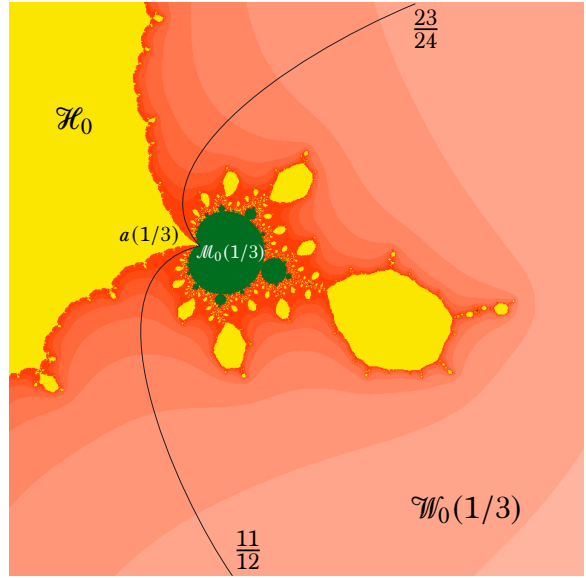

We now look at the complex line $\per_1(0) = \{ P_a \}_{a \in \C}$. When convenient, we identify $P_a \in \per_1(0)$ with the parameter $a \in \C$. The connectedness locus 
$$
\CC_0 := \CC(3) \cap \per_1(0) = \{ a \in \C: K_a \ \text{is connected} \} 
$$
is a full continuum. Every hyperbolic component in $\CC_0$ other than the central component $\HH_0$ is either a {\bit capture component} where some forward iterate of $\om_2$ is in the immediate basin $B_a$ of $\om_1$, or a {\bit Mandelbrot component} where $\om_2$ belongs to the basin of an attracting orbit of $P_a$ other than $\om_1$.\footnote{In the language of \cite{M4}, these components are of {\it disjoint type}.} Looking at Figures \ref{lemonfig} and \ref{onethirdlimb}, one immediately recognizes the capture components as the lemon-shaped bulbs in yellow and the Mandelbrot components as the bulbs in various copies of the Mandelbrot set in green. All hyperbolic components in $\CC_0$ have Jordan curve boundaries \cite{R2}. \vs

Define the map $\kappa: \HH_0 \to \D$ by 
$$
\kappa(a):=\beta_a(\om_2)
$$
(the \Bottcher coordinate of the second critical point in $B_a$). It is easily seen that $\kappa$ is a degree $2$ proper holomorphic map, branched over the center $a=0$, which by \eqref{bottsym} satisfies $\kappa(a)=\kappa(-a)$. Using $\kappa$ we can parametrize each half of $\bd \HH_0$ by the interval $[0,1)$ of internal angles: Let $t \mapsto a(t)$ be the lift of $t \mapsto \e^{2\pi \ii t}$ under $\kappa$ with the initial point $a(0)=-2i/3$ (the lowest point of the lemon in \figref{lemonfig}). One can check that the path $t \mapsto a(t)$ parametrizes the right half of $\bd \HH_0$ from $a(0)=-2i/3$ up to $a(1)=2i/3$, and that the cubic $P_{a(t)}$ is fixed-ray normalized for all $0 \leq t <1$. \vs

%%%%%%%%%%%%%%%%%%%%%%%%%%%%%%%%%%%%%%%%%%%
%%%%%%%%%%%%%%%%%%%%%%%%%%%%%%%%%%%%%%%%%%%

Turning attention to the exterior of $\CC_0$, we can define external rays and a limb structure off of $\HH_0$ similar to the limb structure off of the main cardioid of the Mandelbrot set. The map $\Phi: \per_1(0) \sm \CC_0 \to \C \sm \ov{\D}$ defined by 
$$
\Phi(a):=\phi_a(a)
$$ 
(the \Bottcher coordinate of the escaping co-critical point in $\C \sm K_a$) is a conformal isomorphism with positive derivative at $\infty$. This allows us to define the {\bit parameter rays} 
$$
\RR(\xi):=\Phi^{-1} ( \{ r\e^{2\pi \ii\xi}: r>1 \} ) \qquad (\xi \in \R/\Z).
$$  
The following appears in \cite{M4} (see also \cite{BOT}): \vs

\THM
Let $t \in \R/\Z$ be periodic under $\m_2$ with simulating angles $x_k, y_k$. Then the parameter rays $\RR(x_k-1/3)$ and $\RR(y_k+1/3)$ co-land at $a(t) \in \bd \HH_0$. If $\WW_0(t)$ denotes the connected component of $\C \sm (\RR(x_k-1/3) \cup \RR(y_k+1/3) \cup \{ a(t) \} )$ not containing $\HH_0$, then   
\begin{align*}
\ov{\WW_0(t)} & = \{ a \in \C : \emph{the rays} \ R_a(x_k) \ \emph{and} \ R_a(y_k) \ \emph{co-land} \, \} \\
& = \{ a \in \C : \tau_a=t \}.
\end{align*}
\ENDTHM

We call $\WW_0(t)$ the {\bit (parameter) wake} in $\per_1(0)$ at internal angle $t$. \figref{lemonfig} shows a few examples of wakes. \vs

The co-landing rays $R_a(x_k), R_a(y_k)$ are smooth if $a \in \WW_0(t) \cup \{ a(t) \}$. However, $R_a(x_k)$ is infinitely broken if $a \in \RR(x_k-1/3)$ and $R_a(y_k)$ is infinitely broken if $a \in \RR(y_k+1/3)$ (compare \figref{onehalfdisc} right). The common landing point $z_a(t)$ of $R_a(x_k), R_a(y_k)$ is periodic under $P_a$ of the same period $q$ as that of $t$ under $\m_2$; it is repelling if $a \neq a(t)$ and parabolic with multiplier $1$ if $a=a(t)$. The  orbit of $z_a(t)$ moves holomorphically over $a \in \WW_0(t)$; in particular, the multiplier map $a \mapsto (P_a ^{\circ q})'(z_a(t))$ is analytic in $\WW_0(t)$. \vs 

Define the {\bit (parameter) limb} $\LL_0(t) \subset \per_1(0)$ by
\begin{align*}
\LL_0(t) & := \ov{\WW_0(t)} \cap \CC_0  = \{ a \in \CC_0 : \text{the rays} \ R_a(x_k), R_a(y_k) \ \text{co-land} \, \} \\
& \, =  \{ a \in \CC_0 : \tau_a=t \}.  
\end{align*}
The limb $\LL_0(t)$ is a full continuum and meets $\bd \HH_0$ only at $a(t)$. The parameter $a(t)$ is the root of the {\bit main Mandelbrot copy} $\MM_0(t) \subset \LL_0(t)$. If $q$ is the period of $t$ under $\m_2$, the compact set $\MM_0(t)$ can be described as the set of $a \in \LL_0(t)$ for which $P_a$ is $q$-renormalizable around the critical point $\om_2$. More precisely, there is a homeomorphism $\chi: \MM_0(t) \to \MM$ with $\chi(a(t))=1/4$ having the following property: For every $a \in \MM_0(t)$ there are Jordan domains $U \Subset U'$ containing $\om_2$ such that the restriction $P_a^{\circ q}: U \to U'$ is quadratic-like hybrid equivalent to $Q_{\chi(a)}: z \mapsto z^2+\chi(a)$. \vs

\COR [\cite{R2}] \label{L0}
Let $0 \leq t <1$ have period $q$ under $\m_2$. For every $\nu \in \MM$ there exists a unique parameter $a \in \LL_0(t)$ such that the restriction of $P_a^{\circ q}$ to a neighborhood of $\om_2$ is quadratic-like hybrid equivalent to $Q_\nu$.  
\ENDCOR       

%%%%%%%%%%%%%%%%%%%%%%%%%%%%%%%%%%%%%%%%%%%
%%%%%%%%%%%%%%%%%%%%%%%%%%%%%%%%%%%%%%%%%%%

A parameter $a \in \MM_0(t)$ is called a {\bit dyadic tip} if the forward orbit of $\om_2$ under $P_a^{\circ q}$ eventually hits the ``$\beta$-fixed point'' of the quadratic-like map $P_a^{\circ q}: U \to U'$ (this is simply the point $z_a(t) \in \bd B_a$). Under the aforementioned homeomorphism $\chi$, the dyadic tips of $\MM_0(t)$ correspond to those of $\MM$, i.e., the landing points of the parameter rays $\RR_\MM(p/2^n)$ for $n \geq 1$ and odd $p$. Every dyadic tip of $\MM_0(t)$ is on the boundary of a unique capture component (compare \figref{onethirdlimb}). This phenomenon is a special property of the main Mandelbrot copies in $\CC_0$. Other maximal copies of the Mandelbrot set have no capture components directly attached to them although they are densely accumulated by capture components. \vs 

To complete the picture of the limb structure in $\CC_0$, we need to also consider the left half of $\bd \HH_0$ consisting of $P_a$ for which the conjugate map $P_{-a}$ via $z \mapsto -z$ is fixed-ray normalized. Notice that the involution $z \mapsto -z$ in the dynamical plane corresponds to the automorphism $\theta \mapsto  \theta+1/2$ of $\m_3$ on external angles which exchanges the fixed points $0$ and $1/2$. Thus, we also have wakes $\WW^\ast_0(t)=-\WW_0(t)$ and limbs $\LL^\ast_0(t)=-\LL_0(t)$ characterized by
\begin{align*}
\ov{\WW^\ast_0(t)} & =  \{ a : \text{the rays} \ R_a(x_k+1/2) \ \text{and} \ R_a(y_k+1/2) \ \text{co-land} \} \\
\LL^\ast_0(t) & = \ov{\WW^\ast_0(t)} \cap \CC_0. 
\end{align*}
We can now state the {\it limb structure theorem} in $\per_1(0)$, proved by Roesch in \cite{R2}, as follows: 
\begin{equation}\label{ls}
\CC_0 = \ov{\HH_0} \cup \bigcup_{t \ \text{periodic under} \ \m_2}  \big( \LL_0(t) \cup \LL^\ast_0(t) \big).
\end{equation}

\section{Lemon limbs and renormalization} 

Throughout this section $t \in \R/\Z$ is periodic under $\m_2$ with period $q$ and simulating angles $x_k,y_k$. 

\subsection{Lemon limbs $\LL(t)$}\label{ssec:lemlim}

The combinatorial characterization of parameter limbs in the one-dimensional slice $\per_1(0)$ motivates the following definition in the full parameter space $\PP(3)$: 

\DEF
The {\bit lemon limb} $\LL(t)$ is defined by
$$
\LL(t) := \{ P \in \CC(3): \text{the rays} \ R_P(x_k) \ \text{and} \ R_P(y_k) \ \text{co-land} \}. 
$$
\ENDEF

\noindent 
Thus, $\LL_0(t)=\LL(t) \cap \per_1(0)$. \vs

For $P \in \LL(t)$ the co-landing point $z_k$ of the ray pair $R_P(x_k), R_P(y_k)$ is periodic with period dividing $q$. Moreover, $z_k$ is either repelling or parabolic with $(P^{\circ q})'(z_k) = 1$. \vs

Similar to what we saw in the slice $\per_1(0)$ in \S \ref{lfpar}, we can define the {\bit rotated lemon limbs} $\LL^\ast(t)$ as the image of $\LL(t)$ under the involution $(a,b) \mapsto (-a,-b)$ of $\PP(3)$. In other words, 
$$
\LL^\ast(t) := \{ P \in \CC(3) : \text{the rays} \ R_P(x_k+1/2) \ \text{and} \ R_P(y_k+1/2) \ \text{co-land} \}. 
$$ 
It can be shown that distinct lemon limbs are disjoint: $\LL(t) \cap \LL(s) = \es$ whenever $t \neq s$. However, a limb $\LL(t)$ and a rotated limb $\LL^\ast(s)$ may well intersect. These statements and some related ideas and results are discussed in the Appendix. \vs

One subtlety in the definition of lemon limbs (in contrast to the slices $\LL_0(t)$ or the primary limbs of the Mandelbrot set) is that $\LL(t)$ is {\it not} compact.   

\begin{figure}[t]
\centering
\begin{overpic}[width=\textwidth]{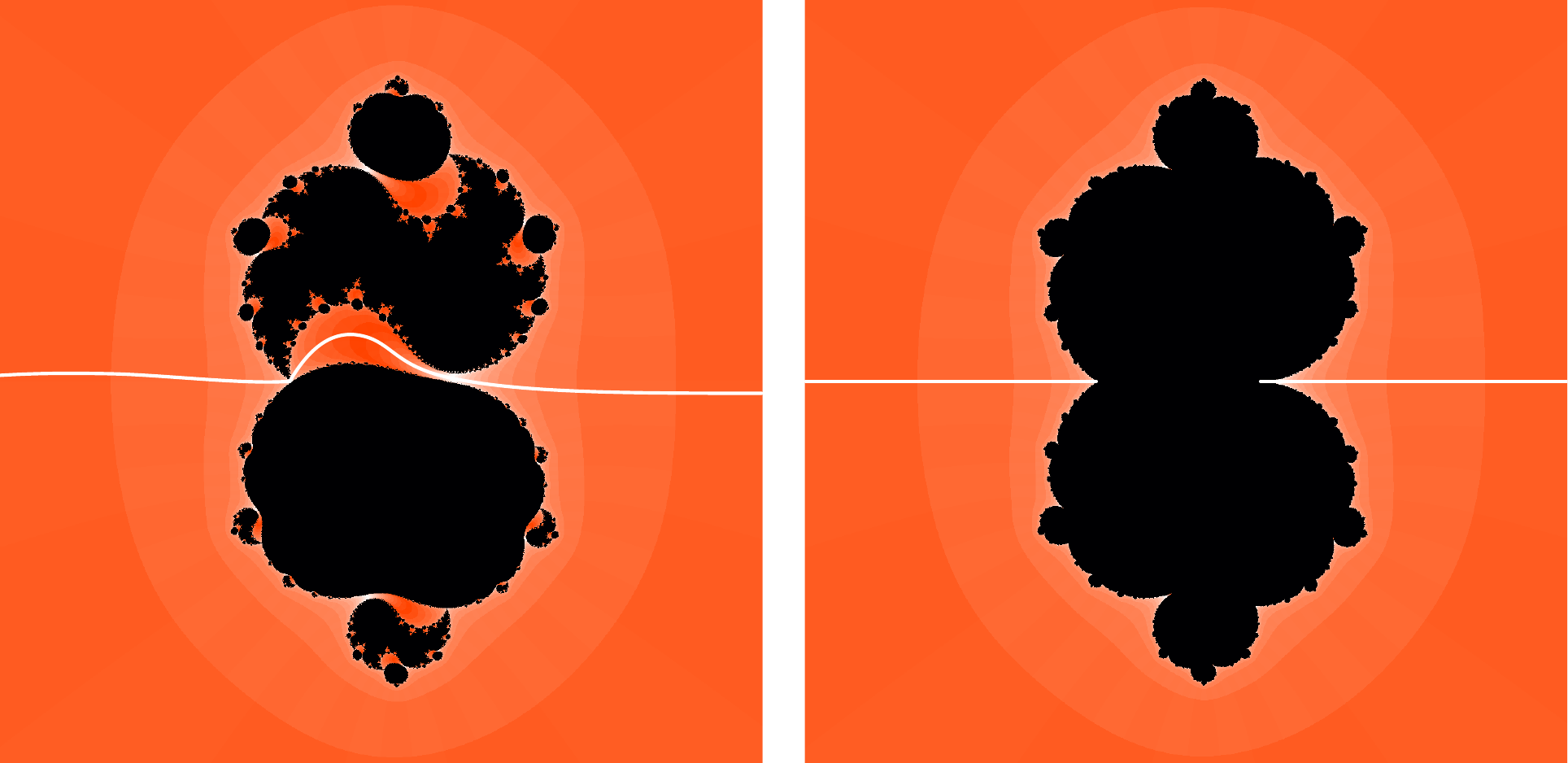}
\put (46.5,24.5) {\footnotesize $0$}
\put (1,27) {\footnotesize $\frac{1}{2}$}
\put (98,25.3) {\footnotesize $0$}
\put (52.5,26.5) {\footnotesize $\frac{1}{2}$}
\put (10,40) {$K_{f_n}$}
\put (62,40) {$K_f$}
\end{overpic}
\caption{\small Illustration of non-compactness of lemon limbs: The co-landing condition defining a limb can fail to persist under taking limits (see \exaref{noncompact}).}  
\label{hetero}
\end{figure}
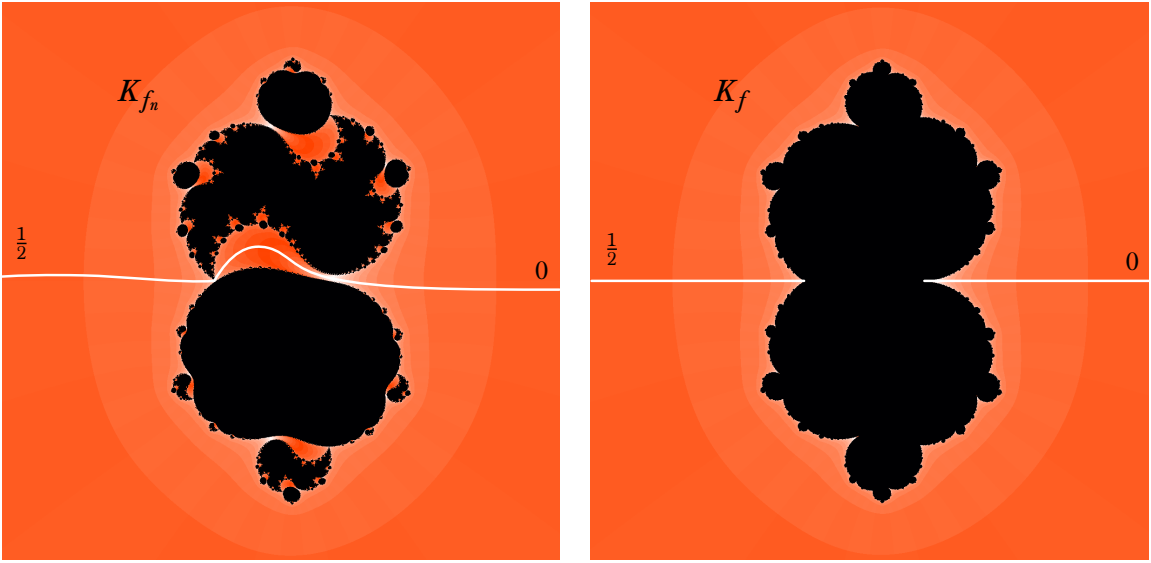

\EX \label{noncompact}
There is a sequence of cubic maps
$$
f_n(z):=z^3-(\sqrt{2}-\eps_n)z^2+\frac{3}{2}z 
$$
with $\eps_n \to 0$ for which the filled Julia set $K_{f_n}$ is connected and the rays $R_{f_n}(0),R_{f_n}(1/2)$ co-land at the repelling fixed point $z=0$. This is certainly true if $\eps_n>0$ (compare \cite[Theorem D]{PZ3}) and it likely holds whenever $\re(\eps_n)>0$ and $\im(\eps_n)/\re(\eps_n) \to 0$ sufficiently fast. But for the cubic $f=\lim_{n \to \infty} f_n$ the ray $R_{f}(0)$ lands at $z=0$ while $R_{f}(0)$ lands at the parabolic fixed point $z=1/\sqrt{2}$ (see \figref{hetero}). 
\ENDEX 

We refer the reader to \remref{hetconj} for a discussion of the characterization of cubics in $\ov{\LL(t)} \sm \LL(t)$. 

%%%%%%%%%%%%%%%%%%%%%%%%%%%%%%%%%%%%%%%%%%%
%%%%%%%%%%%%%%%%%%%%%%%%%%%%%%%%%%%%%%%%%%%

\subsection{Central and peripheral renormalizaions}\label{cpren}

For the rest of the subsequent sections we will be working under our usual assumptions: $0 \leq t_1 < \cdots < t_q <1$ is an $\m_2$-orbit of period $q$ with combinatorics $\sigma \in S_q$. We fix a marked angle $t=t_k$ and let $\OO_k=\{ x_1, \ldots, x_q \}, \OO_{k-1}=\{ y_1, \ldots, y_q \}$ be the simulating orbits for $t$. We also fix a cubic $P \in \LL(t)$ and for simplicity of notation drop the subscript $P$ from all the associated objects. Thus, for $1 \leq i \leq q$ the external rays $R(x_i),R(y_i)$ co-land at a point $z_i \in J$. The {\bit co-landing orbit} $\{ z_1, \ldots, z_q \}$ has period dividing $q$, with each $z_i$ either repelling or parabolic with $(P^{\circ q})'(z_i)=1$. \vs

Define the dynamical wake $W_i$ as the connected component of $\C \sm (R(x_i) \cup R(y_i) \cup \{ z_i \})$ not containing $\om_1$ (this is the analog of $W_a(t_i)$ in the lemon family; see \S\ref{lfdyn}). Under $P$ the wake $W_i$ maps conformally onto $W_{\sigma(i)}$ if $i \neq k$ and onto the entire plane if $i=k$. It follows that $\om_2 \in W_k$ and $\om_1 \in \C \sm \bigcup_{i=1}^q \ov{W_i}$. Note that unlike the lemon family the wakes $W_i$ can now have intersecting boundaries. This happens precisely when the period of the co-landing orbit $\{ z_1, \ldots, z_q \}$ is a proper divisor of $q$ (see \S \ref{clor} and \figref{cheb}). \vs

\begin{figure}[t]
	\centering
	\begin{overpic}[width=\textwidth]{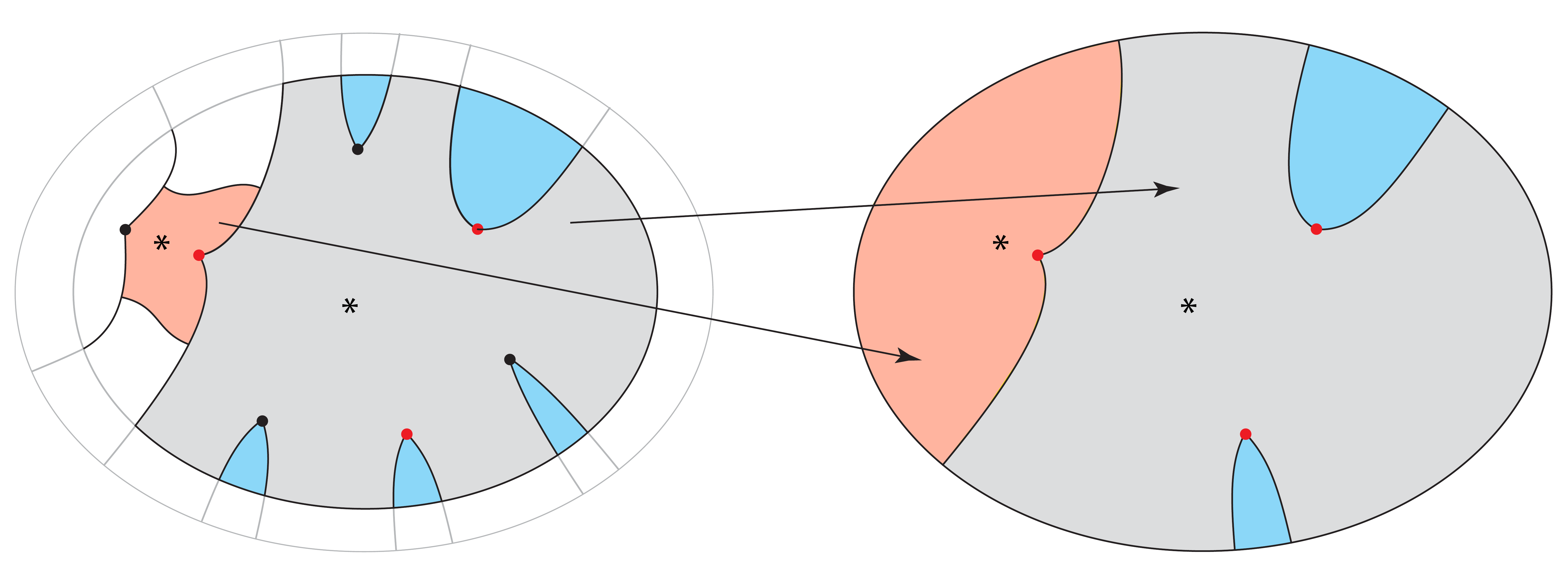}
\put (74.8,15.5) {\tiny $\om_1$}
\put (62.8,19.7) {\tiny $\om_2$}
\put (83,20.5) {\tiny $z_1$}
\put (67,19.5) {\tiny $z_k$}
\put (79,10.3) {\tiny $z_{\sigma(k)}$}
\put (92,31) {\tiny $x_1$}
\put (83.3,35) {\tiny $y_1$}
\put (70,35.3) {\tiny $x_k$}
\put (58.4,6) {\tiny $y_k$}
\put (75.8,0.1) {\tiny $x_{\sigma(k)}$}
\put (82,0.5) {\tiny $y_{\sigma(k)}$}
\put (21,15.5) {\tiny $\om_1$}
\put (9.2,19.7) {\tiny $\om_2$}
\put (29,20.5) {\tiny $z_1$}
\put (13.6,19.5) {\tiny $z_k$}
\put (25,10.3) {\tiny $z_{\sigma(k)}$}
\put (5.8,20.4) {\tiny $z'_k$}
\put (38.2,31) {\tiny $x_1$}
\put (29.5,35) {\tiny $y_1$}
\put (16.2,35.3) {\tiny $x_k$}
\put (4.9,6) {\tiny $y_k$}
\put (0,12.6) {\tiny $x'_k$}
\put (8,32.3) {\tiny $y'_k$}
\put (22.6,0.1) {\tiny $x_{\sigma(k)}$}
\put (28.8,0.5) {\tiny $y_{\sigma(k)}$}
\put (49,24) {\footnotesize $P$}
\put (49,16.5) {\footnotesize $P^{\circ q}$}
\put (24,25) {\footnotesize $U_0$}
\put (75,29) {\footnotesize $U'_0$}
\put (12.2,25.8) {\footnotesize $V_0$}
\put (65,29) {\footnotesize $V'_0$}
	\end{overpic}
	\caption{\small Dynamical wakes and quadratic-like restrictions of a typical cubic $P \in \LL(t)$.}  
	\label{wakes}
\end{figure}
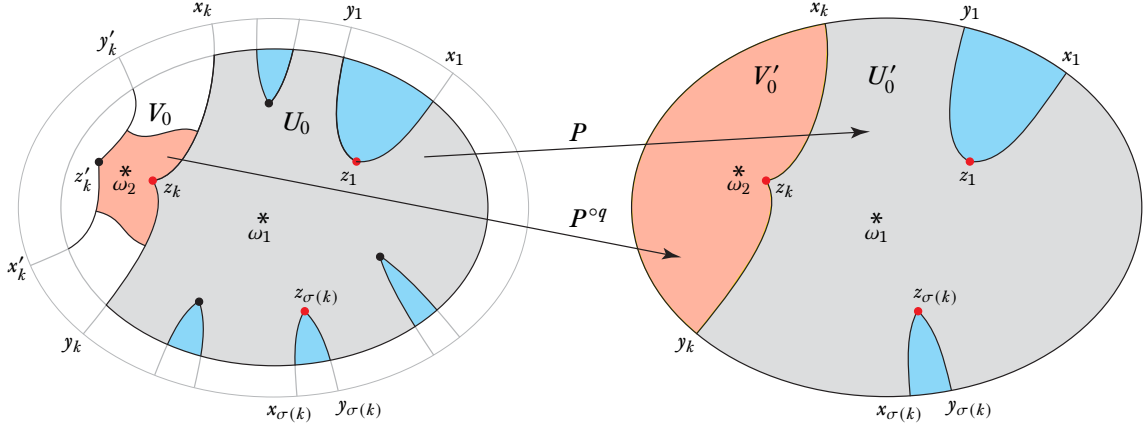

Fix a potential level $s>0$. For every $1 \leq i \leq q$ there is a univalent pull-back $V_i$ of $W_{\sigma(i)}$ outside $\bigcup_{i=1}^q W_i$ (this is the analog of $W_a(t_i+1/2)$ in the lemon family). Let $U'_0$ be the open subset of $\C \sm \bigcup_{i=1}^q \ov{W_i}$ consisting of all points at Green's potential $<s$, and let $U_0$ be the open subset of $\C \sm \bigcup_{i=1}^q (\ov{W_i} \cup \ov{V_i})$ consisting of all points at Green's potential $<s/3$. Then the restriction $P: U_0 \to U'_0$ is proper of degree $2$ (see \figref{wakes}). If $\{ z_1, \ldots, z_q \}$ is ``peripherally repelling,'' i.e., repelling or parabolic with no immediate basin in $\bigcup_{i=1}^q W_i$, we can apply a standard thickening procedure to obtain a quadratic-like restriction $P: U \to U'$ between slightly larger Jordan domains $U, U'$ (see \cite[p. 57]{M2}). Notice that when the period of $\{ z_1, \ldots, z_q \}$ is $<q$ the sets $U_0,U'_0$ are disconnected, but the thickening procedure will always produce connected domains $U, U'$. The domains $U, U'$ can be chosen arbitrarily close to  $U_0,U'_0$ in the Hausdorff metric. According to \cite[Theorem 5.11]{Mc1} the filled Julia set of $P: U \to U'$ does not depend on the choice of the thickened domains, so it can be described as the {\bit central filled Julia set} 
$$
\kcen := \{ z \in K : P^{\circ n}(z) \in \C \sm \bigcup_{i=1}^q W_i \ \text{for all} \ n \geq 0 \}.
$$
Observe that by our normalizations the landing point of the fixed ray $R_P(0)$ always belongs to $\kcen$ while the landing point of the fixed ray $R_P(1/2)$ does not belong to $\kcen$ unless $t=0$ in which case the two landing points coincide. \vs  

%%%%%%%%%%%%%%%%%%%%%%%%%%%%%%%%%%%%%%%%%%%
%%%%%%%%%%%%%%%%%%%%%%%%%%%%%%%%%%%%%%%%%%%

Recall from \S \ref{subsec:sp} that the angles $x'_k:=x_k+1/3, y'_k:=y_k-1/3$ appear in the order $x_k<y'_k<1/2<x'_k<y_k$. Let $W'_k$ be the sub-wake of $W_k$ bounded by the rays $R(y'_k),R(x'_k)$ which co-land at the unique preimage $z'_k$ of $P(z_k)=z_{\sigma(k)}$ in $W_k$. Let $V'_0$ be the sub-domain of $W_k$ consisting of all points at Green's potential $<s$ and $V_0$ be the sub-domain of $W_k \sm \ov{W'_k}$ consisting of all points at Green's potential $<s/3^q$ (see \figref{wakes}). It is easy to see that $\om_2 \in V_0$ and $P^{\circ q}: V_0 \to V'_0$ is proper of degree $2$. If $\{ z_1, \ldots, z_q \}$ is ``centrally repelling,'' i.e., repelling or parabolic with no immediate basin outside $\bigcup_{i=1}^q W_i$, the thickening procedure gives a quadratic-like restriction $P^{\circ q}: V \to V'$ between slightly larger Jordan domains $V, V'$. Again, the filled Julia set of $P^{\circ q}: V \to V'$ does not depend on the choice of the thickened domains, so it can be described as the {\bit peripheral filled Julia set} 
$$
\kper := \{ z \in K : P^{\circ nq}(z) \in \ov{W_k} \sm W'_k \ \text{for all} \ n \geq 0 \}.
$$
Observe that $\kcen$ contains the entire co-landing orbit $\{ z_1, \ldots, z_q \}$ and $\kcen \cap \kper= \{ z_k \}$.
 
\DEF
We call $P \in \LL(t)$ {\bit centrally renormalizable} if $\kcen$ is connected and {\bit peripherally renormalizable} if $\kper$ is connected. We denote the set of all such cubics by $\lcen(t)$ and $\lper(t)$, respectively.    
\ENDEF

By the straightening theorem of Douady and Hubbard \cite{DH2}, when $P$ is centrally renormalizable the quadratic-like restriction $P: U \to U'$ is hybrid equivalent to the quadratic $Q_c$ for a unique $c \in \MM$. Thus, there is a central straightening map 
$$
\chicen: \lcen(t) \to \MM \qquad \text{given by} \qquad  \chicen(P)=c.
$$
Similarly, when $P$ is peripherally renormalizable the quadratic-like restriction $P^{\circ q}: V \to V'$ is hybrid equivalent to $Q_\nu$ for a unique $\nu \in \MM$. This defines a peripheral straightening map 
$$
\chiper: \lper(t) \to \MM \qquad \text{given by} \qquad  \chiper(P)=\nu.
$$ 
(Notice that to simplify the notation we have suppressed the dependence of these maps on $t$.) \vs 
 
In the subsequent discussions we will invoke the following result on the correspondence between the external rays of a polynomial and those of its renormalizations: 

\THM[Ray correspondence]\label{raycor}
Let $P: \C \to \C$ be a polynomial of degree $D \geq 3$ and $\varphi$ be a hybrid equivalence between a polynomial-like restriction $f=P|_U: U \to U'$ and a polynomial $Q$ of degree $d$ with $K_Q$ connected. Let $I$ be the set of angles of external rays of $P$ that accumulate on $K_f$. Then \vs
\begin{enumerate}
\item[(i)]
$I$ is compact and totally disconnected and there is a continuous surjection $\Pi: I \to \R/\Z$ which satisfies the semiconjugacy relation $\Pi \circ \m_D = \m_d \circ \Pi$ and is unique up to postcomposition with some rotation $t \mapsto t+j/(d-1) \ \modd$. \vs
\item[(ii)]
$R_P(\theta)$ lands at $z \in K_f$ if and only if $R_Q(\Pi(\theta))$ lands at $\varphi(z) \in K_Q$. In this case, $\Pi(\theta)=\Pi(\theta')$ for some $\theta' \neq \theta$ if and only if $R_P(\theta')$ also lands at $z$ and one of the two components of $\C \sm (R_P(\theta) \cup R_P(\theta')\cup \{ z\})$ does not meet $K_f$. 
\end{enumerate}    
\ENDTHM

This is proved in \cite{PZ2} in the case $K_P$ is disconnected and in \cite{Le} in the general case. The surjection $\Pi$ in the above result can be extended to a degree $1$ monotone map $\Pi: \R/\Z \to \R/\Z$ by defining it to be constant on each complementary component of $I$. The resulting extension still satisfies the conjugacy relation $\Pi \circ \m_D = \m_d \circ \Pi$ outside the complementary components of $I$ which have length $>1/D$. 

\LEM \label{renorm}
Consider $P \in \LL(t)$ with the co-landing orbit $\{ z_1, \ldots, z_q \}$. \vs
\begin{enumerate}
\item[(i)]
If $P \in \lcen(t)$ and $\varphi$ is any hybrid equivalence between the restriction $P: U \to U'$ and $Q_c$, then for $1 \leq i \leq q$ the external ray $R_c(t_i)$ lands at the point $w_i:=\varphi(z_i) \in K_c$.%
\vs

\item[(ii)]
If $P \in \lper(t)$ and $\psi$ is any hybrid equivalence between the restriction $P^{\circ q}: V \to V'$ and $Q_\nu$, then $\psi(z_k)$ and $\psi(z'_k)$ are the $\beta$-fixed point of $Q_\nu$ and its preimage, respectively.% 
\vs
\end{enumerate}
\ENDLEM

\PROOF
(i) Let $\Pi: I \to \R/\Z$ be the unique semiconjugacy of \thmref{raycor} associated with the central renormalization which satisfies $\Pi(0)=0$. Thus, $\theta \in I$ if and only if $R(\theta)$ accumulates on $\kcen$. Since $\kcen \cap W_i = \es$ for every $1 \leq i \leq q$, we have $\Pi(x_i)=\Pi(y_i)$. It suffices to check that $\Pi(x_i)=\Pi(y_i)=t_i$, for then another application of \thmref{raycor} will show that $R_c(t_i)$ lands at $w_i$. Note that the result would be trivial assuming $\Pi(x_1), \ldots, \Pi(x_q)$ are distinct, for then by monotonicity of $\Pi$ these $q$ points would form an $\m_2$-orbit with combinatorics $\sigma$ and there is only one such orbit namely $\{ t_1, \ldots, t_q \}$. The subtlety here is to rule out any identification among $\Pi(x_1), \ldots, \Pi(x_q)$ that may change the combinatorics of the projection. \vs

Suppose $\Pi(x_i)=\Pi(y_i)$ is not $t_i$ for some (hence all) $i$. Then $\Pi(\OO_k)=\Pi(\OO_{k-1})$ is an $\m_2$-orbit disjoint from $\{ t_1, \ldots, t_q \}$. Choose any $u_1 \in \Pi^{-1}(t_1)$ and consider its $\m_3$-orbit $\{ u_1, \ldots, u_q \} \subset I$, which has the combinatorics $\sigma$ by monotonicity of $\Pi$, and is disjoint from $\OO_k \cup \OO_{k-1}$. By \corref{nothird}, $u_j \in \, ]x_k,y_k[$ for some $j$. Again by monotonicity of $\Pi$, this implies $t_j=\Pi(u_j)=\Pi(x_k)=\Pi(y_k)$, contradicting the assumption that $\Pi(\OO_k)=\Pi(\OO_{k-1})$ is disjoint from $\{ t_1, \ldots, t_q \}$. \vs  

(ii) The argument is similar to but easier than (i). This time let $\Pi: I \to \R/\Z$ be a semiconjugacy (between $\m_3^{\circ q}$ and $\m_2$) given by \thmref{raycor} associated with the peripheral renormalization. Then $\Pi(x_k)=\Pi(y_k)$ since $\kper \subset \ov{W_k} \sm W'_k$. Now $x_k, y_k$ are fixed under $\m_3^{\circ q}$, so $\Pi(x_k)=\Pi(y_k)=0$. Thus, $\psi(z_k)$ is the landing point of the fixed ray $R_\nu(0)$ of $Q_\nu$ and $\psi(z'_k)$ is its other preimage.  
\ENDPROOF 

\subsection{Renormalizable cubics and straightening}

\DEF
The {\bit main renormalization locus} $\NN(t)$ is the set of all cubics in $\LL(t)$ that are both centrally and peripherally renormalizable:
$$
\NN(t)=\lcen(t) \cap \lper(t).
$$ 
\ENDEF 

Note that for every $P \in \NN(t)$ the co-landing orbit $\{ z_1, \ldots, z_q \}$ is repelling. \vs

The following characterization of $\NN(t)$ is immediate:

\COR \label{nchar}
Let $P \in \LL(t)$ have a repelling co-landing orbit. Then $P \in \NN(t)$ if and only if $P^{\circ n}(\om_1) \in \C \sm \bigcup_{i=1}^q \ov{W_i}$ and $P^{\circ nq}(\om_2) \in W_k \sm \ov{W'_k}$ for all $n \geq 0$.  
\ENDCOR

Notice that since $\{ z_1, \ldots, z_q \}$ is repelling, so is the corresponding orbit $\{ w_i \}$ of $Q_c$ and the $\beta$-fixed point of $Q_\nu$, where $c=\chicen(P), \nu=\chiper(P)$. It follows from \lemref{renorm} and the general theory of the Mandelbrot set that $\nu \neq 1/4$ and $c \neq c(t_i)$ for all $1 \leq i \leq q$, where $c(t_i)$ is the landing point of the external ray $\RR_\MM(t_i)$ of the Mandelbrot set (see the discussion in \S \ref{qpms}). Set  
$$
\check{\MM}_t := \MM \sm \{ c(t_1), \ldots, c(t_q) \}.
$$   
In particular, $\check{\MM}_0=\MM \sm \{ 1/4 \}$. \vs  

\thmref{B} in the introduction claims that the straightening map
$$
\Chi: \NN(t) \to \check{\MM}_t \times \check{\MM}_0
$$ 
defined by $\Chi(P)=(\chicen(P), \chiper(P))$ is a homeomorphism.  
The inverse $\Chi^{-1}$ provides a homeomorphic embedding of the product $\check{\MM}_t \times \check{\MM}_0$ into the lemon limb $\LL(t)$. In the terminology of \S \ref{lfpar}, the main Mandelbrot copy $\MM_0(t) \subset \LL_0(t)$ is the closure of $\Chi^{-1}(\{ 0 \} \times \check{\MM}_0)$. An analogous result holds for the rotated limbs $\LL^\ast(t)$ since the corresponding straightening map $\Chi^\ast: \NN(t) \to \check{\MM}_t \times \check{\MM}_0$ satisfies $\Chi^\ast(P_{a,b})=\Chi(P_{-a,-b})$. 

\begin{figure}[t]
	\centering
	\begin{overpic}[width=0.6\textwidth]{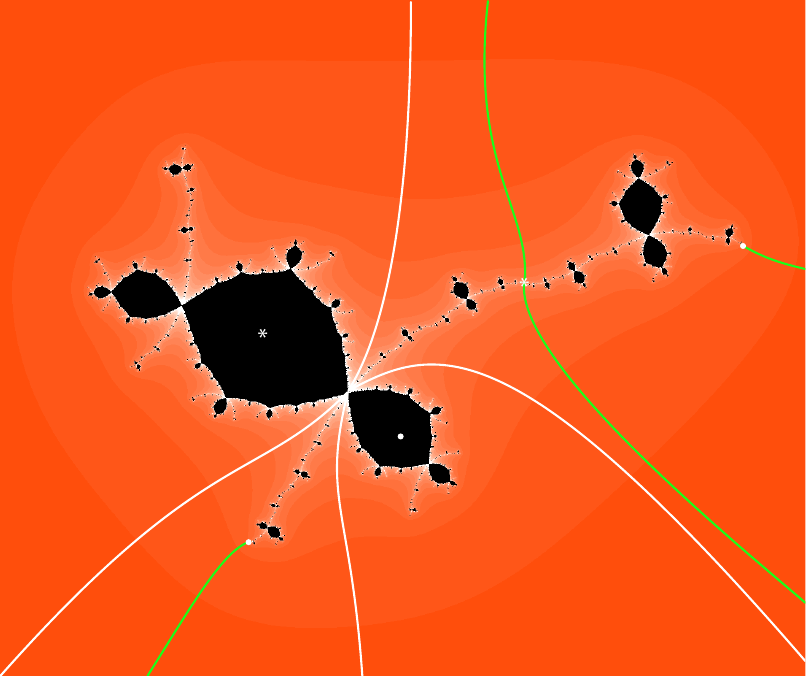}
		\put (48,79.5) {\footnotesize $\frac{1}{4}$}
		\put (5.3,2) {\footnotesize $\frac{5}{8}$}
		\put (46,2) {\footnotesize $\frac{3}{4}$}
		\put (93.7,2) {\footnotesize $\frac{7}{8}$}
		\put (40,36) {\color{white} \tiny $z_1$}
		\put (44,33.5) {\color{white} \tiny $z_2$}
		\put (31,40) {\color{white} \tiny $\om_2$}
		\put (61.6,46.1) {\tiny $\om_1$}
        \put (97.5,47.5) {\footnotesize $0$}
        \put (22,2) {\footnotesize $\frac{2}{3}$}
        \put (61,79.5) {\footnotesize $\frac{2}{9}$}
        \put (97,13.5) {\footnotesize $\frac{8}{9}$}
	\end{overpic}
	\caption{\small Filled Julia set of the cubic $\Chi^{-1}(-2,0)$ in the limb $\LL(1/3)$. Here the critical point $\om_1=0$ maps to a fixed point after two iterates, while $\om_2=-2a$ is periodic with period $2$. The points of the co-landing orbit have merged into a repelling fixed point $z_1=z_2$.}  
	\label{cheb}
\end{figure}
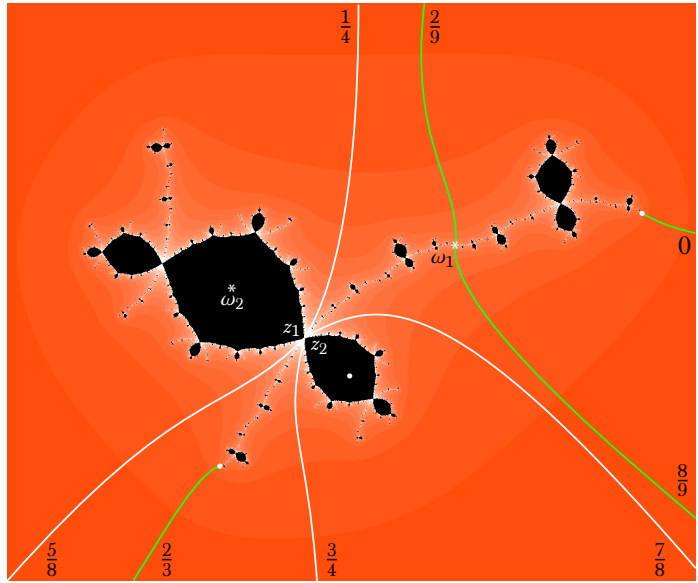

\EX \label{cheb+bas}
\figref{cheb} shows the filled Julia set of the cubic $P=\Chi^{-1}(-2,0)$ in the limb $\LL(1/3)$ which can also be described as the intertwining of the Chebyshev map $z \mapsto z^2-2$ and $z \mapsto z^2-1$. Here the points of the co-landing orbit of $P$ merge into a repelling fixed point $z_1=z_2$ (compare \thmref{coal=}).   
\ENDEX

%%%%%%%%%%%%%%%%%%%%%%%%%%%%%%%%%%%%%%%%%%%
%%%%%%%%%%%%%%%%%%%%%%%%%%%%%%%%%%%%%%%%%%%

Continuity of $\Chi$ follows from the well-known continuity of the straightening map in degree $2$ and the fact that the domains of the quadratic-like restrictions $P: U \to U'$ and $P^{\circ q}: V \to V'$ can be chosen to depend continuously on $P$ (see \cite{DH2} and \cite[Proposition 4.7]{Mc2}). As pointed out in the introduction, the proof of \thmref{B} will be completed in three stages in \S \ref{injchi}, \S \ref{propchi}, and \S \ref{surg}.

\section{A closer look at the co-landing orbit} 

Suppose $P \in \lcen(t)$ and $\chicen(P)=c$. When $c$ is in the main cardioid of the Mandelbrot set, the co-landing orbit $\{ z_1, \ldots, z_q \}$ of $P$ has period $q$ and each $z_i$ is the landing point of precisely two external rays at angles $x_i,y_i$. When $c$ is outside the main cardioid, this picture may be violated in two mutually exclusive ways: The orbit $\{ z_1, \ldots, z_q \}$ may have a third cycle of rays landing on it, or it may coalesce into an orbit of lower period. These cases will be studied in \S \ref{eaor} and \S \ref{clor}, respectively. As both situations can be directly translated into the dynamical plane of $Q_c$, we begin by recalling a few relevant facts about periodic orbits of quadratic polynomials and the external rays landing on them. These descriptions have been known for the most part since the seminal work of Douady and Hubbard \cite{DH1}, but our brief presentation is closer to Milnor's in \cite{M3} to which we refer for further details. 

\subsection{Quadratic portraits and the Mandelbrot set}\label{qpms} 

Let $c \in \MM$ and $Z$ be a repelling or parabolic orbit of $Q_c$ of period $p$. For each $w \in Z$ consider the finite set $A_w \subset \Q/\Z$ of angles $t$ such that the ray $R_c(t)$ lands at $w$. In the language of \cite{M3} the collection $\{ A_w: w \in Z \}$ is the {\bit portrait} associated with $Z$. Each restriction $\m_2: A_w \to A_{Q_c(w)}$ is an order-preserving bijection, so the first return map $\m_2^{\circ p}: A_w \to A_w$ has a well-defined (combinatorial) rotation number $\rho$ independent of $w$. It follows that $\bigcup_{w \in Z} A_w$ is a disjoint union of $N \geq 1$ cycles of common period $q=rp$ for some $r \geq 1$. Thus, the cardinality of each $A_w$ is $N r$ and the rotation number $\rho$ is a reduced fraction of the form $s/r$ with $0 \leq s <r$. The portrait of $Z$ is necessarily of one of the following types: \vs

1. {\bit Trivial}: $r=N=1 \Rightarrow p=q, \rho=0$. In this case the orbit period equals the ray period and each $w \in Z$ is the landing point of a unique external ray. \vs

2. {\bit Primitive}: $r=1, N=2 \Rightarrow p=q, \rho=0$. Again the orbit period equals the ray period and each $w \in Z$ is the landing point of precisely two external rays belonging to different cycles. \vs           

3. {\bit Satellite}: $r\geq 2, N=1 \Rightarrow p<q, \rho \neq 0$. In this case the orbit period is a proper divisor of the ray period and each $w \in Z$ is the landing point of $r$ external rays all of which belong to the same cycle. \vs

Every non-zero periodic angle $t \in \Q/\Z$ under $\m_2$ has a unique $\MM$-{\bit partner}, i.e., an angle $\hat{t} \in \Q/\Z$ of the same period such that the parameter rays $\RR_{\MM}(t),\RR_{\MM}(\hat{t})$ co-land. Their common landing point $c(t)=c(\hat{t})$ is then the root of a unique hyperbolic component of $\MM$. The connected component of $\MM \sm \{ c(t) \}$ not containing $c=0$ union the root point $c(t)$ defines the Mandelbrot limb $\LL_{\MM}(t,\hat{t})$ which is characterized by the following condition: For $c \in \MM$ the dynamic rays $R_c(t),R_c(\hat{t})$ co-land at a repelling point of $Q_c$ if and only if $c \in \LL_{\MM}(t,\hat{t}) \sm \{ c(t) \}$ and at a parabolic point of $Q_c$ if and only if $c=c(t)$. The portrait of the $Q_c$-orbit $Z$ of this common landing point is non-trivial and its {\bit characteristic arc}, i.e., the shortest arc among all components of $(\R/\Z) \sm A_w$ as $w$ runs through $Z$, is bounded by $t,\hat{t}$. Conversely, if the portrait of a periodic orbit $Z$ of some $Q_c$ is non-trivial, it has a well-defined characteristic arc bounded by a pair $t,\hat{t}$ of $\MM$-partners, and $c \in \LL_{\MM}(t,\hat{t})$. \vs

Every periodic orbit $\{ t_1, \ldots, t_q \}$ under $\m_2$ with period $q>1$ contains at most one pair of $\MM$-partners. In other words, if $\hat{t}_i$ denotes the $\MM$-partner of $t_i$, then the relation $\hat{t}_i=t_j$ can hold for at most one unordered pair $i,j$. Let $c \in \MM$ and consider the set $Z$ of the landing points of $R_c(t_1), \ldots, R_c(t_q)$. If $c \in \LL_{\MM}(t_i,\hat{t}_i)$ for some $i$, the dynamic rays $R_c(t_i), R_c(\hat{t}_i)$ co-land and the portrait of $Z$ is primitive if $t_i,\hat{t}_i$ belong to different $\m_2$-orbits and is satellite if $t_i,\hat{t}_i$ belong to the same $\m_2$-orbit. Conversely, if the portrait of $Z$ is non-trivial, then $c \in \LL_{\MM}(t_i,\hat{t}_i)$ as we pointed out above. \vs 

We summarize the above observations in the following 

\COR \label{portsum}
Let $\{ t_1, \ldots, t_q \}$ be a periodic orbit of $\m_2$ with period $q>1$ and for each $1 \leq i \leq q$ let $\hat{t}_i$ be the unique $\MM$-partner of $t_i$. Take $c \in \MM$ and consider the set $Z$ of the landing points of the dynamic rays $R_c(t_1), \ldots, R_c(t_q)$. Then the portrait of $Z$ is \vs 

$\bullet$ primitive if $c \in \LL_{\MM}(t_i,\hat{t}_i)$ for some $i$ and $t_i,\hat{t}_i$ belong to different $\m_2$-orbits; \vs

$\bullet$ satellite if $c \in \LL_{\MM}(t_i,\hat{t}_i)$ for some $i$ and $t_i,\hat{t}_i$ belong to the same $\m_2$-orbit; \vs

$\bullet$ trivial otherwise. 
\ENDCOR

\EX
The period $q=4$ orbits of $\m_2$ are listed below: \vs

\bgroup
\def\arraystretch{1.2}
\begin{center}
	\begin{tabular}{|c|c|c|c|}
		\hline
		orbit & combinatorics in $S_4$ & degree & rotation number  \\
		\hline 
		$\{ 1,2,4,8 \}/15$ & $(1234)$ & $1$ & $1/4$\\ 
		\hline
		$\{ 3,6,9,12 \}/15$ & $(1243)$ & $2$ & $-$ \\ 
		\hline
		$\{ 7,11,13,14 \}/15$ & $(1432)$ & $1$ & $3/4$\\ 
		\hline 
	\end{tabular}
\end{center}
\egroup

\vs 

\noindent
The $\MM$-partner pairs (in multiples of $1/15$) are 
$$
(1,2), \ (3,4), \ (6,9), \ (7,8), \ (11,12), \ (13,14).
$$
The degree $2$ orbit $\{ 3,6,9,12 \}/15=\{ 1,2,3,4 \}/5$ is part of a non-trivial portrait in the three limbs $\LL_\MM(3/15,4/15)$, $\LL_\MM(11/15,12/15)$ and $\LL_\MM(6/15,9/15)$, and generates a trivial portrait elsewhere. In the first two limbs the portrait is primitive while in the third limb it is satellite and the corresponding landing points form an orbit of lower period $p=2$.      
\ENDEX

\EX
The period $q=5$ orbits of $\m_2$ are listed below: \vs

\bgroup
\def\arraystretch{1.2}
\begin{center}
	\begin{tabular}{|c|c|c|c|}
		\hline
		orbit & combinatorics in $S_5$ & degree & rotation number  \\
		\hline 
		$\{ 1,2,4,8,16 \}/31$ & $(12345)$ & $1$ & $1/5$\\
		\hline
		$\{ 5,9,10,18,20 \}/31$ & $(13524)$ & $1$ & $2/5$\\ 
		\hline
		$\{ 3,6,12,17,24 \}/31$ & $(12354)$ & $2$ & $-$\\  
		\hline
		$\{ 7,14,19,25,28 \}/31$ & $(12543)$ & $2$ & $-$\\ 
		\hline
		$\{ 11,13,21,22,26 \}/31$ & $(14253)$ & $1$ & $3/5$\\ 
		\hline
		$\{ 15,23,27,29,30 \}/31$ & $(15432)$ & $1$ & $4/5$\\ 
		\hline
	\end{tabular}
\end{center}
\egroup

\vs 

\noindent
The $\MM$-partner pairs (in multiples of $1/31$) are 
\begin{align*}
(1,2),(3,4),(5,6),(7,8),(9,10),(11,12),(13,18),(14,17),\\
(15,16),(19,20),(21,22),(23,24),(25,26),(27,28),(29,30)
\end{align*}
Observe that the two orbits of degree $2$ do not contain any pair of $\MM$-partners, so they can never be part of a satellite portrait. 
\ENDEX

\subsection{External angles of the co-landing orbit}\label{eaor}   

For each $P \in \LL(t)$ the two cycles of external rays with angles in $\OO_k=\{ x_1, \ldots, x_q \}$ and $\OO_{k-1}=\{ y_1, \ldots, y_q \}$ co-land on the orbit $\{ z_1, \ldots, z_q \}$. It is known that the number of distinct cycles of external rays that land on a periodic orbit of a polynomial of degree $D$ is at most $D$ (in fact, the sharp bound is one  plus the number of critical values of the polynomial; see \cite{M3} for the quadratic case already encountered in \S \ref{qpms}, and \cite{K} for the higher degree case). It follows that in our cubic case there can be at most one additional cycle of rays landing on $\{ z_1, \ldots, z_q \}$. The following lemma describes when such third cycle exists for centrally renormalizable cubics: 

\THM [Third cycle comes from central primitive portraits]\label{3rd}
Let $P \in \lcen(t)$ and $\chicen(P)=c$. Take any hybrid equivalence $\varphi$ between $P: U \to U'$ and the quadratic $Q_c$. The following conditions are equivalent: \vs
\begin{enumerate}
\item[(i)]
There is a third cycle of rays landing on $\{ z_1, \ldots, z_q \}$. \vs
\item[(ii)]
The orbit $\{ \varphi(z_1), \ldots, \varphi(z_q) \}$ of $Q_c$ has a primitive portrait. \vs
\item[(iii)]
$c \in \LL_\MM(t_i,\hat{t}_i)$ for some $1 \leq i \leq q$, where $t_i,\hat{t}_i$ do not belong to the same $\m_2$-orbit. 
\end{enumerate}
\ENDTHM

The condition (iii) puts $c$ in the union of either $q-1$ or $q$ Mandelbrot limbs according as the $\m_2$-orbit $\{ t_1, \ldots, t_q \}$ contains a pair of $\MM$-partners or not (see \S \ref{qpms}). Note also that these equivalent conditions imply that $\{ z_1, \ldots, z_q \}$ has exact period $q$, so no merging in the co-landing orbit can occur.    

\PROOF
As the equivalence (ii) $\Leftrightarrow$ (iii) is part of \corref{portsum}, it suffices to verify (i) $\Leftrightarrow$ (ii). Consider the semiconjugacy $\Pi:I \to \R/\Z$ of \thmref{raycor} associated with $P:U \to U'$ and normalized by $\Pi(0)=0$. By (the proof of) \lemref{renorm}, $\Pi(x_i)=\Pi(y_i)=t_i$ and $R_c(t_i)$ lands at $w_i:=\varphi(z_i)$ for every $1 \leq i \leq q$. \vs

Let us first prove (ii) $\Rightarrow$ (i). Suppose the portrait of $\{ w_1, \ldots, w_q \}$ is primitive. Then every $w_i$ is the landing point of exactly two external rays $R_c(t_i)$ and $R_c(s_i)$, where the $\m_2$-orbits $\{ t_1, \ldots, t_q \}$ and $\{ s_1, \ldots, s_q \}$ are disjoint (notice that this choice of labeling no longer guarantees that the $s_i$ are in positive cyclic order). Choose $u_1 \in \Pi^{-1}(s_1)$ and define $u_{\sigma^i(1)}:=\m_3^{\circ i}(u_1)$, so $\Pi(u_i)=s_i$ for every $1 \leq i \leq q$. Then $\{ u_1, \ldots, u_q \}$ is an $\m_3$-orbit of period $q$ disjoint from $\OO_k \cup \OO_{k-1}$. By \thmref{raycor}, for every $i$ the external ray $R(u_i)$ lands at $z_i$. \vs

Let us now prove (i) $\Rightarrow$ (ii). Suppose there is a third cycle of rays $R(u_1), \ldots, R(u_q)$ such that $R(u_i)$ lands at $z_i$ for every $i$. Setting $s_i:=\Pi(u_i)$, it follows from \thmref{raycor} that $R_c(s_i)$ lands at $w_i$ (again, keep in mind that the $u_i$ or $s_i$ are no longer in positive cyclic order). We need to check that $\{ s_1, \ldots, s_q \}$ is an $\m_2$-orbit of period $q$ disjoint from $\{ t_1, \ldots, t_q \}$. \vs

Set $\OO:= \{ u_1, \ldots, u_q \}$. First we observe that 
\begin{equation}\label{noui}
\OO \cap \bigcup_{i=1}^q \, ]x_i,y_i[ = \es.
\end{equation}
In fact, if this fails, by taking a suitable iterate under $\m_3$ we can find some $u \in \OO \cap \, ]x_k,y_k[$. As $\ov{W_k} \cap \{ z_1, \ldots, z_q \}= \{ z_k \}$, the ray $R(u)$ must land at $z_k$, so it cannot be contained in the sub-wake $W'_k$. Thus $u \in \, ]x_k,y'_k[ \, \cup \, ]x'_k,y_k[$, which implies $\m_3(u) \in \, ]x_{\sigma(k)},y_{\sigma(k)}[$, which in turn shows that every $]x_i,y_i[$ contains precisely one element of $\OO$. It follows that $\OO$ has the combinatorics $\sigma$ and $\#(\OO \cap [0,1/2[)=k$ or $k-1$. By the characterization \eqref{deploy}, $\OO=\OO_k$ or $\OO=\OO_{k-1}$, which is a contradiction. \vs    

Now suppose the orbit $\{ s_1, \ldots, s_q \}$ coincides with $\{ t_1, \ldots, t_q \}$. Then we can re-label the angles in $\OO$ as $v_1, \ldots, v_q$, where $\Pi(v_i)=t_i$ for every $i$. From monotonicity of $\Pi$ we see that the $v_i$ are in positive cyclic order and $\OO$ has the combinatorics $\sigma$. It then follows from \eqref{noui} and $\Pi(v_k)=t_k$ that $v_k \in \, ]y_{k-1},x_k[$ or $v_k \in \, ]y_k,x_{k+1}[$. This implies $\#(\OO \cap [0,1/2[)=k$ or $k-1$, which as before leads to a contradiction. 
\ENDPROOF

\EX \label{3cycles}
For $t=2/5$ the simulating orbits are $\OO_2=\{ 8,24,56,72 \}/80$ and $\OO_1= \{ 17,51,59,73 \}/80$ with combinatorics $\sigma=(1243)$ (see \exaref{2/5}). If $P \in \lcen(2/5)$ and $\chicen(P)=c \in \LL_\MM(3/15,4/15)$, \thmref{3rd} shows that there is a third cycle of external rays that lands on the co-landing orbit of $P$. It is easy to see that the angles of this third cycle form the $\m_3$-orbit $\OO= \{ 2,6,18,54 \}/80$ with combinatorics $\tau=(1234)$. Under the monotone surjection $\Pi$ the orbits $\OO_2 \cup \OO_1$ and $\OO$ map to the $\m_2$-orbits $\{ 3,6,9,12 \}/15$ and $\{ 1,2,4,8 \}/15$, respectively.
\ENDEX

\subsection{Merging of the co-landing orbit}\label{clor}

The co-landing orbit $\{ z_1, \ldots, z_q \}$ of a cubic $P \in \LL(t)$ can merge into an orbit of period $<q$. When $P$ is centrally renormalizable, this situation has a direct interpretation in terms of the quadratic polynomial which represents the central hybrid class of $P$. In fact, the following characterization follows immediately from the discussion of \S \ref{qpms}, especially \corref{portsum}: 

\THM [Merging comes from central satellite portraits]\label{coal=}
Let $P \in \lcen(t)$ and $\chicen(P)=c$. Take any hybrid equivalence $\varphi$ between $P: U \to U'$ and the quadratic $Q_c$. The following conditions are equivalent: \vs 
\begin{enumerate}
\item[(i)]
The co-landing orbit $\{ z_1, \ldots, z_q \}$ of $P$ merges into an orbit of period $<q$. \vs
\item[(ii)]
The orbit $\{ \varphi(z_1), \ldots, \varphi(z_q) \}$ of $Q_c$ has a satellite portrait. \vs
\item[(iii)]
The $\m_2$-orbit $\{ t_1, \ldots, t_q \}$ contains a (necessarily unique) pair of $\MM$-partners $t_i,t_j=\hat{t}_i$ and $c \in \LL_\MM(t_i,t_j)$. 
\end{enumerate} 
\ENDTHM

\begin{figure}[t]
\centering
\begin{overpic}[width=\textwidth]{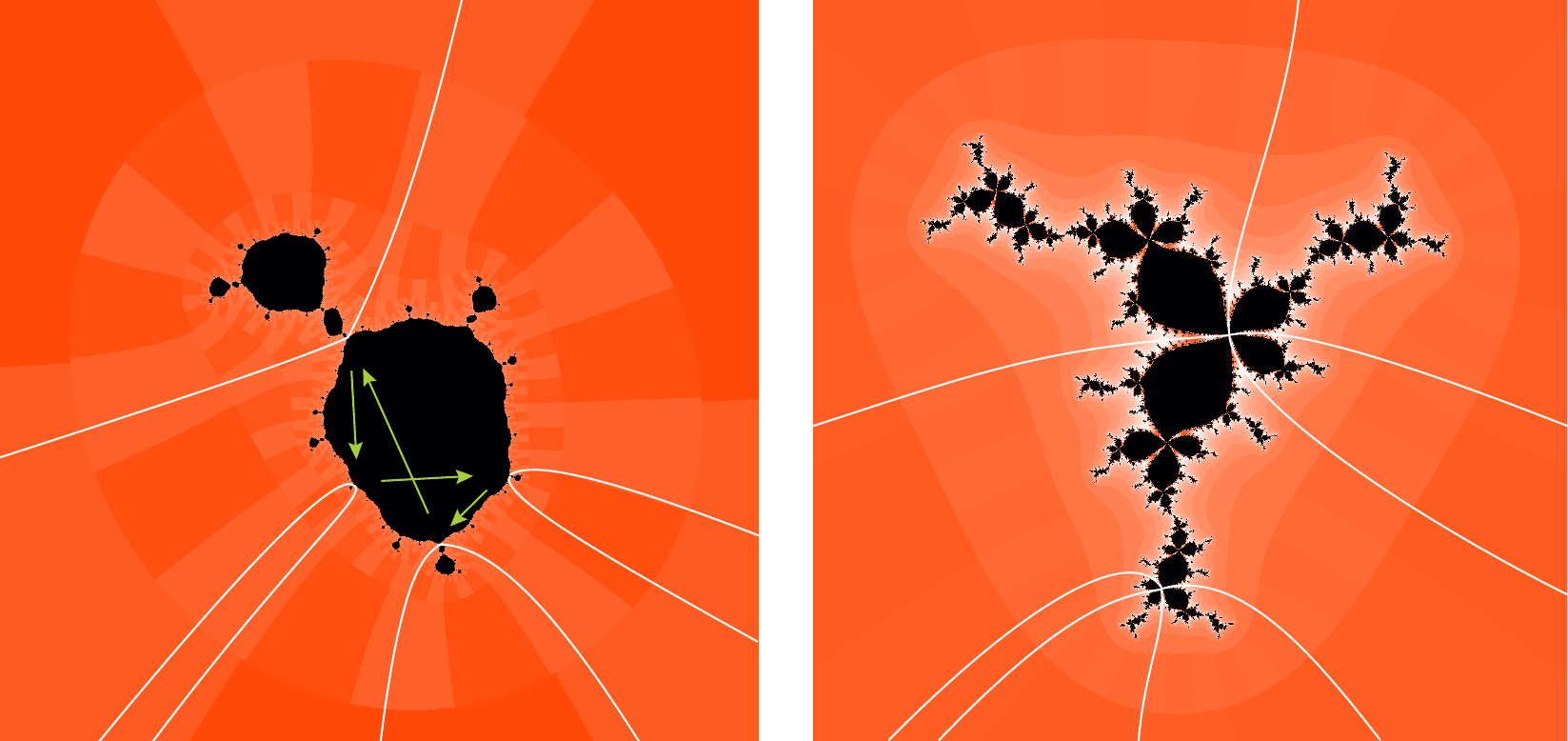}
\put (29.5,44.5) {\footnotesize $\frac{17}{80}$}
\put (0.8,20) {\footnotesize $\frac{44}{80}$}
\put (5,1.4) {\footnotesize $\frac{51}{80}$}
\put (12.5,1.4) {\footnotesize $\frac{52}{80}$}
\put (25,1.4) {\footnotesize $\frac{59}{80}$}
\put (41,1.4) {\footnotesize $\frac{68}{80}$}
\put (46,4) {\footnotesize $\frac{73}{80}$}
\put (46,14.8) {\footnotesize $\frac{76}{80}$}
\put (22.3,24.7) {\color{white} \tiny $z_1$}
\put (22.5,16.5) {\color{white} \tiny $z_2$}
\put (26.8,13.3) {\color{white} \tiny $z_3$}
\put (30.5,17) {\color{white} \tiny $z_4$}
\put (83,44.5) {\footnotesize $\frac{17}{80}$}
\put (52.5,22) {\footnotesize $\frac{44}{80}$}
\put (55,1.4) {\footnotesize $\frac{51}{80}$}
\put (63,1.4) {\footnotesize $\frac{52}{80}$}
\put (73,1.4) {\footnotesize $\frac{59}{80}$}
\put (88,1.4) {\footnotesize $\frac{68}{80}$}
\put (97.5,7) {\footnotesize $\frac{73}{80}$}
\put (97.5,22) {\footnotesize $\frac{76}{80}$}
\put (76.7,26.3) {\color{white} \tiny $z_1$}
\put (74,10.5) {\color{white} \scalebox{0.5}{$z_2$}}
\put (74.4,8.8) {\color{white} \scalebox{0.5}{$z_3$}}
\put (79,24.8) {\color{white} \tiny $z_4$}
\end{overpic}
\caption{\small Illustration of \thmref{coal=} for two cubics $P \in \lcen(1/5)$ with $\chicen(P)=c$. Left: $c=0$ and the co-landing orbit $\{ z_1, z_2, z_3, z_4 \}$ has period $4$. Right: $c \in \LL_{\MM}(2/5,3/5)$ and the co-landing orbit merges into a period $2$ orbit $\{ z_1=z_4, z_2=z_3 \}$.}  
	\label{wee}
\end{figure}
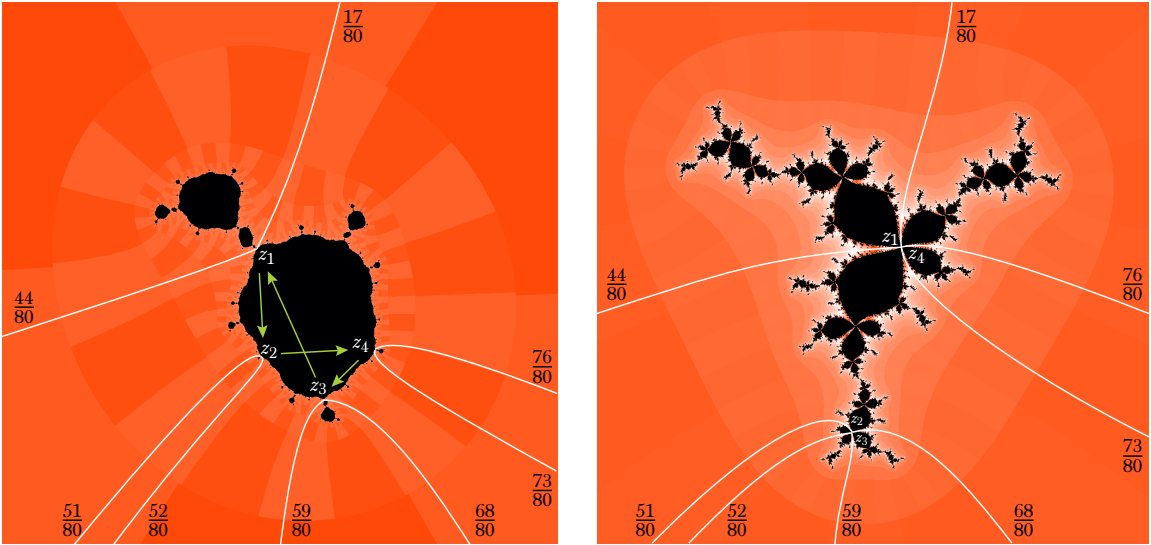

\EX \label{pmerg}
We have already encountered merging in \exaref{cheb+bas} and \figref{cheb}, with $P \in \lcen(1/3)$ and $\chicen(P)=c=-2 \in \LL_\MM(1/3,2/3)$. Here is a more interesting example: For $t=1/5$ the simulating orbits are $\OO_1= \{ 17,51,59,73 \}/80$ and $\OO_0=\{ 44, 52, 68, 76 \}/80$, with combinatorics $\sigma=(1243)$ (see \exaref{2/5}). The $\m_2$-orbit $\{ 1, 2, 3, 4 \}/5$ contains the $\MM$-partner pair $(2/5,3/5)$. By \thmref{coal=}, if $P \in \lcen(1/5)$ and $\chicen(P)=c \in \LL_\MM(2/5,3/5)$, the co-landing orbit of $P$ merges into an orbit whose period is a proper divisor of $4$. This period cannot be $1$ since $\OO_0, \OO_1$ are not rotation cycles, so the period must be $2$. \figref{wee} shows the filled Julia sets and relevant rays for $c=0$ and $c \in \LL_\MM(2/5,3/5)$. 
\ENDEX

%%%%%%%%%%%%%%%%%%%%%%%%%%%%%%%%%%%%%%%%%%%
%%%%%%%%%%%%%%%%%%%%%%%%%%%%%%%%%%%%%%%%%%%

An obstruction to merging in $\LL(t)$ can be formulated in purely combinatorial terms. Suppose the co-landing orbit $\{ z_1, \ldots, z_q \}$ of $P \in \LL(t)$ has period $p=q/r$, a proper divisor of $q$. Then each $z_i$ is the landing point of $r$ external rays with angles in $\OO_k=\{ x_1, \ldots, x_q \}$ and $r$ external rays with angles in $\OO_{k-1}=\{ y_1, \ldots, y_q \}$. Moreover, as the fixed points of $P^{\circ p}$, the $z_i$ have a common combinatorial rotation number of the form $s/r$. It follows that the combinatorics $\sigma \in S_q$ of the $\m_2$-orbit $\{ t_1, \ldots, t_q \}$ is {\bit dynamically reducible}, i.e., it has the property that the power $\sigma^p$ splits into a product of $p$ ``unlinked'' rotation $r$-cycles, all of the same rotation number $s/r$. The condition of being unlinked is stronger than being disjoint and can be described geometrically as follows: Identifying $i \in \Z/(q\Z)$ with $i/q \in \R/\Z$, we can think of elements of $S_q$ as acting on $q$ points on the circle labeled in positive cyclic order. With this convention, two $r$-cycles $(a_1 \cdots a_r), (b_1 \cdots b_r) \in S_q$ are unlinked if the entire set $\{ a_1, \ldots, a_r \}$ is contained in a single connected component of $\R/\Z \sm \{ b_1, \ldots, b_r \}$.

\THM [Combinatorial characterization of merging]
There is a cubic $P \in \lcen(t)$ whose co-landing orbit $\{ z_1, \ldots, z_q \}$ merges into an orbit of period $<q$ if and only if the combinatorics $\sigma \in S_q$ of the $\m_2$-orbit of $t$ is dynamically reducible.  
\ENDTHM

\PROOF[Sketch of the proof]
As noted above, the existence of such $P$ implies that $\sigma$ is dynamically reducible. Conversely, suppose that for some proper divisor $p=q/r$ of $q$ the power $\sigma^p$ is the product $\pi_1 \cdots \pi_p$ of unlinked rotation $r$-cycles with a common rotation number. We can then form the associated ``formal portrait''  in the sense of \cite{M3} on the set $\{ t_1, \ldots, t_q \}$ in which the angles corresponding to each $\pi_j$ are grouped in one class. The characteristic arc of this portrait determines a Mandelbrot limb in which the portrait is realized by the rays $R_c(t_1), \ldots, R_c(t_q)$. Any $P \in \lcen(t)$ with $\chicen(P)$ within this Mandelbrot limb (whose existence is guaranteed by \thmref{B}) will have a merged co-landing orbit. 
\ENDPROOF

\REM
The decomposition of a dynamically reducible cycle into unlinked cycles of shorter lengths is unique. Hence the corresponding portrait (and therefore the pattern of merging) is uniquely determined by $\sigma$.
\ENDREM

\EX
We conclude with a few examples that illustrate the above theorem. \vs
\begin{enumerate}[leftmargin=*]
\item[$\bullet$]
Every $q$-cycle $\sigma \in S_q$ with $\deg(\sigma)=1$ is dynamically reducible (here $r=q, p=1$). The corresponding parameter $c$ belongs to a primary limb attached to the main cardioid of $\MM$, with $R_c(t_1), \ldots, R_c(t_q)$ all co-landing at the $\alpha$-fixed point of $Q_c$. Any cubic $P \in \lcen(t)$ with $\chicen(P)=c$ has its co-landing orbit merged into a single fixed point. \vs 
	
\item[$\bullet$]
The $4$-cycle $\sigma=(1243)$ is dynamically reducible with the decomposition $\sigma^2=(14)(23)$. The corresponding formal portrait $\{ \{ 1/5,4/5 \}, \{ 2/5, 3/5 \} \}$ is realized by any $c \in \LL_\MM(2/5,3/5)$. The co-landing orbit of corresponding cubics merge into an orbit of period $2$, as in \exaref{pmerg} and \figref{wee}. \vs   
	
\item[$\bullet$]
If $q$ is a prime number and $\deg(\sigma)=2$, then $\sigma$ is dynamically irreducible, so no merging is possible. The $5$-cycle $\sigma=(12354)$, realized by the $\m_2$-orbit $\{ 3,6,12,17,24 \}/31$, provides an example of this type. \vs
	
\item[$\bullet$]
A more interesting example is the $6$-cycle $\sigma=(123465)$ with $\deg(\sigma)=2$, realized by the $\m_2$-orbit $\{ 1,2,4,8,11,16 \}/21$, which is dynamically irreducible since $\sigma^2=(136)(245)$ and $\sigma^3=(14)(26)(35)$ and neither of these decompositions is unlinked. Again, no merging is ever possible with this combinatorics.  
\end{enumerate}
\ENDEX

\section{Injectivity of straightening} \label{injchi}

The goal of this section is to prove the following 

\THM
The straightening map $\Chi: \NN(t) \to \check{\MM}_t \times \check{\MM}_0$ is injective. 
\ENDTHM

Given $P, \tilde{P} \in \NN(t)$ with $\Chi(P)=\Chi(\tilde{P})$, there are natural hybrid equivalences between the central and peripheral renormalizations of $P$ and $\tilde{P}$ that can be glued together to obtain a quasiconformal homeomorphism $\Phi_0: \C \to \C$ on the entire plane which partially conjugates $P$ to $\tilde{P}$. A standard pull-back argument then allows to promote $\Phi_0$ to a sequence $\Phi_n: \C \to \C$ of uniformly quasiconformal homeomorphisms which conjugate $P$ to $\tilde{P}$ on larger and larger domains. The limit $\Phi=\lim_{n \to \infty} \Phi_n$ is a global quasiconformal conjugacy and satisfies $\bar{\bd} \Phi=0$ a.e. outside the {\bit residual Julia set} of $P$ defined by   
$$
K_\infty := K \sm \bigcup_{n \geq 0} P^{-n}(\kcen \cup \kper).
$$ 
The main technical point then is to show that $K_\infty$ has measure zero, for that would imply $\Phi$ is an affine map. \vs 

The construction of $\Phi$ is standard. Take $P,\tilde{P} \in \NN(t)$ and let $\phi, \tilde{\phi}$ and $G,\tilde{G}$ denote their \Bottcher coordinates and Green's functions. If $\Chi(P)=\Chi(\tilde{P})$, the quadratic-like restriction pairs $P:U \to U', \tilde{P}: \tilde{U} \to \tilde{U}'$ and $P^{\circ q}: V \to V', \tilde{P}^{\circ q}: \tilde{V} \to \tilde{V}'$ are hybrid equivalent. Recall from \S \ref{cpren} that $U, U'$ and $V,V'$ are obtained by slightly thickening the open sets $U_0, U'_0$ and $V_0, V'_0$. We have corresponding pairs of open sets $\tilde{U}_0, \tilde{U}'_0$ and $\tilde{V}_0, \tilde{V}'_0$ for $\tilde{P}$. There is a hybrid equivalence $\varphi: U_0' \to \tilde{U}'_0$ conjugating $P$ to $\tilde{P}$ which we may assume satisfies $\phi= \tilde{\phi} \circ \varphi$ on $\bd U_0'$. Similarly, there is a hybrid equivalence $\psi: V_0' \to \tilde{V}'_0$ conjugating $P^{\circ q}$ to $\tilde{P}^{\circ q}$ which satisfies $\phi= \tilde{\phi} \circ \psi$ on $\bd V'_0$. Define a homeomorphism $\Phi_0: \C \to \C$ by
$$
\Phi_0 =
\begin{cases}
\tilde{\phi}^{-1} \circ \phi & \quad \text{on} \ G^{-1}([s,\infty[), \\
\varphi & \quad \text{on} \ U'_0, \\
\psi  & \quad \text{on} \ V'_0, \\
\tilde{P}^{-i} \circ \psi \circ P^{\circ i} & \quad \text{on} \  W_{\si^{-i}(k)}\cap G^{-1}([0,s[), \ 1 \leq i \leq q-1.
\end{cases}
$$
Here $\tilde{P}^{-i}$ in the last line is the inverse branch which maps the wake $\tilde{W}_k$ of $\tilde{P}$ to $\tilde{W}_{\si^{-i}(k)}$. Evidently $\Phi_0$ is an $M$-quasiconformal homeomorphism for some $M>1$ and satisfies $\bar{\bd} \Phi_0=0$ a.e. on 
$$
G^{-1}([s,\infty[) \cup \kcen \cup \bigcup_{i=0}^{q-1} P^{\circ i}(\kper).
$$
Since $\Phi_0$ conjugates $P$ and $\tilde{P}$ on their respective postcritical sets, we can define $\Phi_n: \C \to \C$ inductively by lifting $\Phi_{n-1}\circ P$ under $\tilde{P}$, so $\Phi_n$ is $M$-quasiconformal, coincides with $\Phi_{n-1}$ on the set
$$
G^{-1}([3^{-n}s,\infty[) \cup \bigcup_{\ell=0}^{n-1} P^{-\ell}(\kcen \cup \bigcup_{i=0}^{q-1} P^{\circ i}(\kper)),
$$
and satisfies $\bar{\bd} \Phi_n=0$ a.e. there. Thus, after passing to a subsequence, $\{ \Phi_n \}$ converges to an $M$-quasiconformal homeomorphism $\Phi: \C \to \C$ which conjugates $P$ to $\tilde{P}$ everywhere and satisfies $\bar{\bd} \Phi=0$ a.e. on the set 
$$
(\C \sm K) \cup \bigcup_{\ell=0}^{\infty} P^{-\ell}(\kcen \cup \bigcup_{i=0}^{q-1} P^{\circ i}(\kper))=\C \sm K_\infty.
$$
Assuming $K_\infty$ has measure zero, it follows that $\Phi$ is conformal a.e. in $\C$, hence is an affine map. Since $\Phi$ is tangent to the identity at $\infty$ and $\Phi(0)=0$, we must have $\Phi=\operatorname{id}$ and $P=\tilde{P}$. \vs

We now address the measure zero statement that was essential in the above argument:  

\THM \label{K=null}
For every $P\in \NN(t)$ the residual Julia set $K_\infty$ has Lebesgue measure $0$.
\ENDTHM

Compare \cite[Lemma 6.1]{IK} where a similar result is proved in a more general and abstract setting. Below we give a concrete independent proof based on hyperbolic geometry estimates in \lemref{lipcor}. See also \cite[Theorem 6.1]{EY} for the special case of intertwining cubics with a restriction on their quadratic renormalizations. \vs

Let $\Om_0$ denote the complement of the compact invariant set $\kcen \cup \bigcup_{i=0}^{q-1} P^{\circ i}(\kper)$. Evidently $\Om_0$ is a backward-invariant domain which avoids the postcritical set of $P$. For $n \geq 1$ set $\Om_n:=P^{-n}(\Om_0)$. Then $\Om_n \subsetneq \Om_{n-1}$ and each $\Om_n$ is conformally isomorphic to the punctured disk, the puncture being at $\infty$. For simplicity we will write $\lambda_n$ for the hyperbolic metric $\la_{\Om_n}$, $B_n(w,r)$ for the hyperbolic ball $B_{\Om_n}(w,r)$, and $\area_n(\cdot)$ for the hyperbolic area in $\Om_n$. By the Schwarz lemma, $\lambda_n>\lambda_{n-1}$ in $\Om_n$ and in particular $\la_n>\la_0$ in $\Om_n$. Since the restriction $P^{\circ n}: (\Om_n,\lambda_n) \to (\Om_0, \lambda_0)$ is a covering map and therefore a local isometry, it follows that $P^{\circ n}: (\Om_n,\lambda_0) \to (\Om_0, \lambda_0)$ is expanding (although not uniformly). \vs

The residual Julia set $K_\infty$ is the set of points in $J=\bd K$ which stay in $\Om_0$ under the iterations of $P$:
$$
K_\infty = J \cap \bigcap_{n \geq 0} \Om_n.
$$
Recall that for each $1 \leq i \leq q$, $W_i$ is the dynamical wake bounded by the pair $R(x_i), R(y_i)$ co-landing at $z_i$ and $W'_k$ is the sub-wake of $W_k$ bounded by the pair $R(y'_k), R(x'_k)$ co-landing at $z'_k$. Observe that $W'_k \subset \Om_0$.  

\LEM \label{hole}
There are constants $0<\eps<r$ such that for every $w \in W'_k$ sufficiently close to $z'_k$ the ball $B_0(w,r)$ contains a ball $B_0(u,\eps)$ disjoint from $K$. 
\ENDLEM

\PROOF
Set $j:=\sigma(k)$ and $\rho:=(P^{\circ q})'(z_j)$. Let $\psi: \C \to \C$ be the extended linearizer of $P^{\circ q}$ at $z_j$ that satisfies $\psi(0)=z_j$, $\psi'(0)=1$ and $\psi(\rho z)=P^{\circ q}(\psi(z))$ for all $z$. Let $E, \ga_x, \ga_y$ denote the connected components of $\psi^{-1}(P(\kper)), \psi^{-1}(R(x_j)), \psi^{-1}(R(y_j))$ containing $0$. Set $\Om := \C \sm E$ and define $W$ as the sub-sector of $\Om$ bounded by $\ga_x,\ga_y$. The domains $\Om, W$ as well as the rays $\ga_x, \ga_y$ are invariant under the linear map $z \mapsto \rho z$.% 
\vs

A small round disk $D$ centered at $0$ maps isomorphically under $\psi$ to a small neighborhood of $z_j$, which in turn maps isomorphically under a local branch $g$ of $P^{-1}$ to a small neighborhood $D'$ of $z'_k$. The isomorphism $g \circ \psi : D \to D'$ carries $\Om \cap D$ to $\Om_0 \cap D'$ and $W \cap D$ to $W'_k \cap D'$. Since the hyperbolic metrics of $\Om_0$ and $\Om_0 \cap D'$ are asymptotically equal near $z'_k$ and those of $\Om$ and $\Om \cap D$ are asymptotically equal near $0$, it suffices to prove the existence of constants $0<\eps<r$ such that for every $w \in W$ the ball $B_\Om(w,r)$ contains a ball $B_\Om(u,\eps)$ disjoint from $\psi^{-1}(K)$. This is easy to see once we pass to the quotient space $\tilde{\Om}:=\Om/\langle z \mapsto \rho z \rangle$ which is a conformal cylinder (of modulus less than $2\pi/\log |\rho|$) embedded in the complex torus $(\C \sm \{ 0 \})/\langle z \mapsto \rho z \rangle$. The rays $\ga_x, \ga_y$ project to two homotopically non-trivial loops $\tilde{\ga}_x, \tilde{\ga}_y$ in $\tilde{\Om}$ and the sub-cylinder bounded by them is precisely the quotient $\tilde{W}$ of $W$. As $\tilde{\ga}_x, \tilde{\ga}_y$ have annular neighborhoods of modulus $\pi/(q \log 3)$ that are disjoint from the quotient $\tilde{K}$ of $\psi^{-1}(K)$, we can find an $\eps>0$ such that $B_{\tilde{\Om}}(u,\eps) \cap \tilde{K}=\es$ for all $u \in \tilde{\ga}_x \cup \tilde{\ga}_y$. Let $\delta>0$ be the hyperbolic diameter of $\tilde{W}$ in $\tilde{\Om}$ and $r:=\delta+\eps$. Then for every $w \in \tilde{W}$ the ball $B_{\tilde{\Om}}(w,r)$ contains the $\eps$-neighborhood of $\tilde{W}$, hence it contains the ball $B_{\tilde{\Om}}(u,\eps)$ for every $u \in \tilde{\ga}_x \cup \tilde{\ga}_y$. 
\ENDPROOF

For the proof of \thmref{K=null} we also need a priori estimates on the geometry of the inclusion of one hyperbolic Riemann surface into another, which we formulate in Lemmas \ref{lip} and \ref{lipcor} below.  

\LEM \label{lip}
Let $X \subsetneq Y$ be hyperbolic Riemann surfaces such that the inclusion $i: X \hookrightarrow Y$ induces an injection $\pi_1(X) \to \pi_1(Y)$. Then $\log \| i' \| = \log (\la_Y/\la_X) : X \to \, ]-\infty, 0[$ is a Lipschitz function (with constant $2$). 
\ENDLEM  

\PROOF
Fix $x_0,x_1 \in X$, take holomorphic universal covering maps $p_X: (\D,0) \to (X,x_0)$ and $p_Y: (\D,0) \to (Y,x_0)$, and lift $i$ to a holomorphic map $f: (\D,0) \to (\D,0)$. By the assumption on the fundamental groups, $f$ is injective. Take a minimal geodesic in $X$ joining $x_0$ to $x_1$ and lift it under $p_X$ to a geodesic in $\D$ joining $0$ to some $w$. Then $\dist_X(x_0,x_1)=\dist_{\D}(0,w)= \log((1+|w|)/(1-|w|))$. Now use Koebe's distortion theorem to estimate 
\begin{align*}
\frac{\| i'(x_0) \|}{\| i'(x_1) \|} & = \frac{\| f'(0) \|}{\| f'(w) \|} = \frac{\la_{\D}(w)}{\la_{\D}(f(w))} \ \frac{|f'(0)|}{|f'(w)|} = \frac{1-|f(w)|^2}{1-|w|^2} \ \frac{|f'(0)|}{|f'(w)|} \\
& \leq \frac{1}{1-|w|^2} \ \frac{(1+|w|)^3}{1-|w|} = \left(\frac{1+|w|}{1-|w|} \right)^2 = \exp(2 \dist_{\D}(0,w)).
\end{align*}
Taking logarithms gives  
$$
\log \| i'(x_0) \| - \log \| i'(x_1) \| \leq 2 \dist_X(x_0,x_1),
$$
and the result follows by symmetry. 
\ENDPROOF

\LEM \label{lipcor}
In the situation of \lemref{lip}, fix $x_0 \in X$ and let $r>0$ be smaller than the injectivity radius of $X$ at $x_0$, so the ball $B_X(x_0,r)$ is embedded. \vs
\begin{enumerate}
\item[(i)]
For every $0<\eps<r$ there is a constant $C(r,\eps)>0$ such that if $B_X(x,\eps) \subset B_X(x_0,r)$, then 
$$
\frac{\area_Y(B_X(x,\eps))}{\area_Y(B_X(x_0,r))} \geq C(r,\eps). \vs
$$
\item[(ii)]
The ball $B_X(x_0,r)$ has bounded geometry in $Y$ in the following sense: 
\begin{align*}
\rin & := \max \{ s: B_Y(x_0,s) \subset B_X(x_0,r) \} \geq \tfrac{1}{2} \, (1-\e^{-2r}) \, \| i'(x_0) \|, \\[2pt]
\rout & := \min \, \{ s: B_Y(x_0,s) \supset B_X(x_0,r) \} \leq \tfrac{1}{2} \, (\e^{2r}-1) \ \ \| i'(x_0) \|.  
\end{align*}
In particular, $\rin/\rout \geq \e^{-2r}$.
\end{enumerate}  
\ENDLEM

\PROOF 
We will use the following inequalities given by \lemref{lip}: 
$$
\e^{-2 \dist_X(z,x_0)} \leq \frac{\| i'(z) \|}{\| i'(x_0) \|} \leq \e^{2 \dist_X(z,x_0)} \qquad \text{for all} \ z \in X.
$$ 

(i) We have 
\begin{align*}
\area_Y(B_X(x,\eps)) & = \iint_{B_X(x,\eps)} \la_Y^2(z) \, |dz|^2 = \iint_{B_X(x,\eps)} \| i'(z) \|^2 \la_X^2(z) \, |dz|^2 \\
& \geq \e^{-4r} \| i'(x_0) \|^2 \area_X(B_X(x,\eps)), 
\end{align*}
and
\begin{align*}
\area_Y(B_X(x_0,r)) & = \iint_{B_X(x_0,r)} \la_Y^2(z)\,  |dz|^2 = \iint_{B_X(x_0,r)} \| i'(z) \|^2 \la_X^2(z) \, |dz|^2 \\
& \leq \e^{4r} \| i'(x_0) \|^2 \area_X(B_X(x_0,r)). 
\end{align*}
Thus, 
\begin{align*}
\frac{\area_Y(B_X(x,\eps))}{\area_Y(B_X(x_0,r))} & \geq \e^{-8r} \ \frac{\area_X(B_X(x,\eps))}{\area_X(B_X(x_0,r))} = \e^{-8r} \ \frac{\area_{\D}(B_{\D}(0,\eps))}{\area_{\D}(B_{\D}(0,r))} \\
& = \e^{-8r} \, \frac{\cosh(\eps)-1}{\cosh(r)-1} =: C(r,\eps). \vs
\end{align*} 

(ii) Find $p,q \in \bd B_X(x_0,r)$ such that $\dist_Y(x_0,p)=\rin$ and $\dist_Y(x_0,q)=\rout$. Take the $Y$-geodesic $\ga: [0,r'] \to \ov{B_Y(x_0,\rin)}$ with $\ga(0)=x_0, \ga(r')=p$, parametrized so that $\| \ga'(t) \|_X = \la_X(\ga(t)) \, |\ga'(t)|=1$ for all $t$. Then $\dist_X(x_0,\ga(t)) \leq t$, and in particular $r \leq r'$. Hence,   
\begin{align*}
\rin & = \int_0^{r'} \la_Y(\ga(t)) \, |\ga'(t)| \ dt = \int_0^{r'} \| i'(\ga(t)) \| \ dt \geq \| i'(x_0) \| \int_0^{r'} \e^{-2t} \, dt \\
& \geq \| i'(x_0) \| \int_0^r \e^{-2t} \, dt = \frac{1}{2}(1-\e^{-2r}) \, \| i'(x_0) \|. 
\end{align*}  
Similarly, take the $X$-geodesic $\eta: [0,r] \to \ov{B_X(x_0,r)}$ with $\eta(0)=x_0, \eta(r)=q$, parametrized so that $\| \eta'(t) \|_X = \la_X(\eta(t)) \, |\eta'(t)|=1$ for all $t$. Then $\dist_X(x_0,\eta(t))=t$ and  
\begin{align*}
\rout & \leq \int_0^r \la_Y(\eta(t)) \, |\eta'(t)| \ dt = \int_0^r \| i'(\eta(t)) \| \ dt \leq \| i'(x_0) \| \int_0^r \e^{2t} \, dt \\
& = \frac{1}{2}(\e^{2r}-1) \, \| i'(x_0) \|. \hfill \qedhere
\end{align*}
\ENDPROOF

\PROOF[Proof of \thmref{K=null}]
Assume by way of contradiction that $K_\infty$ has positive measure. Since the forward orbit of a.e. point in $J$ converges to the postcritical set of $P$ (see for example \cite[Theorem 3.9]{Mc1}), we can find a Lebesgue density point $w_0 \in K_\infty$ whose orbit $\{ w_n :=P^{\circ n}(w_0) \}$ converges to the postcritical set. Evidently the orbit $\{ w_n \}$ must visit the sub-wake $W'_k$ infinitely often. Since the postcritical set can meet $\ov{W'_k}$ only at $z'_k$, it follows that some subsequence of $\{ w_n \}$ converges to $z'_k$, i.e., there is an increasing sequence $S \subset \N$ such that $w_n \to z'_k$ as $n \in S$ tends to $\infty$. By \lemref{hole} there are constants $0<\eps<r$ such that for all large $n \in S$ the ball $B_0(w_n,r)$ contains a ball of radius $\eps$ disjoint from $K$. Note that $B_0(w_n,r)$ is embedded since $\Om_0$ is a disk punctured at $\infty$. Lifting under the covering map $P^{\circ n}: (\Om_n, \la_n) \to (\Om_0, \la_0)$ then shows that the ball $B_n(w_0,r)$ is embedded and contains a ball $B_n(u_n,\eps)$ disjoint from $K$. By \lemref{lipcor}(i) (applied to $X=\Om_n,Y=\Om_0$),  
$$
\frac{\area_0(B_n(u_n,\eps))}{\area_0(B_n(w_0,r))} \geq C(r,\eps).
$$
On the other hand, the inclusion $i_n: \Om_n \to \Om_0$ satisfies $\| i'_n(w_0) \| = \la_0(w_0)/\la_n(w_0) \to  0$ since the Euclidean distance between $w_0$ and $\bd \Om_n$ tends to $0$ as $n \to \infty$. By \lemref{lipcor}(ii), $B_n(w_0,r)$ is squeezed between two balls $B_0(w_0,s_n)$ and $B_0(w_0,r_n)$, where $s_n/r_n \geq \e^{-2r}$ and $0<s_n \leq r_n \leq \tfrac{1}{2}(\e^{2r}-1) \, \| i'_n(w_0) \| \to 0$ as $n \to \infty$. It follows that
\begin{align*}
\frac{\area_0(B_n(u_n,\eps))}{\area_0(B_0(w_0,r_n))} 
& \geq  C(r,\eps) \,  \frac{\area_0(B_0(w_0,s_n))}{\area_0(B_0(w_0,r_n))} = C(r,\eps) \,  \frac{\cosh(s_n)-1}{\cosh(r_n)-1}\\
& = C(r,\eps) \,  \frac{s_n^2(1+O(s_n^2))}{r_n^2(1+O(r_n^2))} \geq \frac{1}{2} \e^{-4r} C(r,\eps)
\end{align*}
for all large $n \in S$. Since $\la_0$ is asymptotically a multiple of the Euclidean metric near $w_0$, it follows that the ball $B_0(w_0,r_n)$ is uniformly round, i.e., has Euclidean area comparable to (Euclidean diameter)$^2$, and that 
$$
\frac{\area(B_0(w_0,r_n) \sm K)}{\area(B_0(w_0,r_n))} \geq \frac{\area(B_n(u_n,\eps))}{\area(B_0(w_0,r_n))} \geq  \frac{1}{4} \e^{-4r} C(r,\eps)
$$
for all large $n \in S$. This contradicts the hypothesis that $w_0$ is a density point of $K_\infty$. 
\ENDPROOF

\section{Continuity of the inverse straightening} \label{propchi}

So far we have seen that $\Chi: \NN(t) \to \Chi(\NN(t))$ is a continuous bijection. The goal of this section is to prove continuity of the inverse of this map: 

\THM \label{proper}
If $P_n, P_\infty \in \NN(t)$ and $\Chi(P_n)\to \Chi(P_\infty)$, then $P_n \to P_\infty$.  
\ENDTHM

The proof of \thmref{proper} has two main ingredients. The first is the fact that the combinatorial rotation number of a repelling orbit controls its multiplier. Suppose a degree $d$ polynomial map with connected Julia set has a repelling fixed point $\zeta$ of multiplier $\lambda$ at which $N$ cycles of external rays of common combinatorial rotation number $0 \leq p/q <1$ land. Then, according to {\bit Yoccoz's inequality} (see \cite{Hu} and \cite{P}), there is a suitable branch of $\log \lambda$ that satisfies 
$$
\frac{|\log \lambda - 2\pi \ii p/q|^2}{\re(\log \lambda)} \leq \frac{2 \log d}{Nq}.  
$$
This means that $\log \lambda$ belongs to the closed disk of radius $\log d/(Nq)$ tangent to the imaginary axis at $2\pi \ii p/q$. It follows from this inequality that if $\lambda \to 1$, then $p/q \to 0$. Conversely, if $p/q$ remains non-zero and $p/q \to 0$, then $\lambda \to 1$. We invoke this observation in the following  

\LEM \label{same}
Let $f_n$ be a sequence of degree $d$ polynomials with quadratic-like restrictions $f_n: U_n \to U'_n$ and with both $f_n, f_n|_{U_n}$ having connected Julia sets. Take a hybrid equivalence $\varphi_n$ between $f_n|_{U_n}$ and a quadratic polynomial $Q_{a_n}$. If $\zeta_n$ is the $\alpha$-fixed point of $f_n|_{U_n}$ and the multiplier $\lambda_n:=f'_n(\zeta_n)$ tends to $1$ as $n \to \infty$, then the same is true of the multiplier $\mu_n:=Q'_{a_n}(\varphi_n(\zeta_n))$, hence $a_n \to 1/4$.   
\ENDLEM

\PROOF
As the multiplier of a non-repelling fixed point is invariant under a hybrid equivalence, we may assume that $\zeta_n$ and therefore $\varphi_n(\zeta_n)$ is repelling for every $n$. It follows from $\lambda_n \to 1$ and Yoccoz's inequality that the combinatorial rotation number of $\zeta_n$ tends to $0$, so the same must be true of the $\alpha$-fixed point $\varphi_n(\zeta_n)$ of $Q_{a_n}$. This combinatorial rotation number cannot be equal to zero since $a_n \in \MM$. Another application of Yoccoz's inequality now shows $\mu_n \to 1$.  
\ENDPROOF

The second ingredient in the proof of \thmref{proper} is the behavior of periodic external rays of polynomial maps under Hausdorff limits, which is the subject of our recent work in \cite{PZ3}. Since presenting the main results of that paper in full generality would divert our focus from the task at hand, we instead distill a special case that applies to our cubic setup. Let $\theta \in \R/\Z$ have period $q$ under $\m_3$. Suppose $P_n \to P$ in $\CC(3)$ and $\zeta_n$ and $\zeta$ are the landing points of the external rays $R_n:=R_{P_n}(\theta)$ and $R:=R_{P}(\theta)$, respectively. After passing to a subsequence, we may assume $\zeta_n \to \zeta_\infty$ and the closed rays $\ov{R}_n:=R_n \cup \{ \zeta_n \}$ converge in the Hausdorff metric to a continuum $L$ which necessarily contains $\ov{R}=R \cup \{ \zeta \}$ and $\zeta_\infty$. If $\zeta$ is repelling, then $\zeta=\zeta_\infty$ and $L=\ov{R}$ (this follows from the well-known {\bit stability} of rays landing on repelling points; see for example \cite{DH1} or \cite{GM}). However, if $\zeta$ is parabolic, $L$ can be strictly larger than $\ov{R}$. In this case, every point of $L \sm \ov{R}$ in the Fatou set belongs to a $P^{\circ q}$-invariant analytic arc in $L$ connecting $\zeta$ to $\zeta_\infty$ on which the action of $P^{\circ q}$ is conjugate to a translation. We call such an arc a ``homoclinic'' if $\zeta=\zeta_\infty$ and a  ``heteroclinic'' otherwise.    

\THM \label{het}
If $L \neq \ov{R}$, then $\zeta$ is a parabolic fixed point of $P^{\circ q}$ with an immediate attracting basin $B$ of exact period $q$ such that one of the following holds: \vs 

\begin{enumerate}
\item[(i)]
$\zeta=\zeta_\infty $. Then $L=R \cup \{ \zeta \} \cup \gamma_1 \cup \gamma_2 \cup \cdots$ where the $\gamma_i$ are homoclinic arcs (possibly one but at most countably many) in $B$ connecting $\zeta$ to itself. They can be labeled so that $\gamma_{i+1}$ is inside $\gamma_i$, so there is always an outermost homoclinic $\gamma_1$ in the nested chain. \vs    

\item[(ii)]
$\zeta \neq \zeta_\infty$. Then $L=R \cup \{ \zeta \} \cup \gamma \cup \{ \zeta_\infty \}$, where $\gamma$ is a unique heteroclinic arc in $B$ connecting $\zeta$ to $\zeta_\infty$. In this case \vs
  
\begin{enumerate}
\item[$\bullet$]
$\zeta_\infty$ is a repelling fixed point of $P^{\circ q}$.  \vs

\item[$\bullet$]
If $p=q/r$ is the period of $\zeta$, either $r=1$ or $r>1$ in which case the period of $\zeta_\infty$ must be $q$. The multiplier $(P^{\circ p})'(\zeta)$ is a primitive $r$th root of unity and the multiplicity of $\zeta$ as a fixed point of $P^{\circ q}$ is $r+1$. \vs  

\item[$\bullet$]
Each of the two connected components of $B \sm \gamma$ contains at least one critical point of $P^{\circ q}$. In particular, both critical points of $P$ must belong to the disjoint union $B \cup P(B) \cup \cdots \cup P^{\circ q-1}(B)$. 
\end{enumerate}
\end{enumerate} 
\ENDTHM  

See Theorems A, B, C as well as the Basic Structure Lemma and Lemma 6.2 in \cite{PZ3}. \vs           

The proof of \thmref{proper} begins as follows. After passing to a subsequence, we may assume $P_n$ converges to some $P \in \CC(3)$ and the co-landing orbit $\{ z_{1,n}, \ldots, z_{q,n} \}$ of $P_n$ converges to a periodic orbit $\{ z_1, \ldots, z_q \}$ of $P$, where each $z_i$ is fixed under $P^{\circ q}$. The critical points of $P$ are $\om_1=\lim_{n \to \infty} \om_{1,n}$ and $\om_2=\lim_{n \to \infty} \om_{2,n}$. We will show that the rays $R(x_k), R(y_k)$ of $P$ co-land at $z_k$ (so $P \in \LL(t)$) and that $z_k$ is repelling. In particular, the closed rays $\ov{R}_n(x_k), \ov{R}_n(y_k)$ of $P_n$ converge in the Hausdorff metric to the closed rays $\ov{R}(x_k), \ov{R}(y_k)$. Thus, the domains of the quadratic-like restrictions $P_n: U_n \to U'_n$ and $P_n^{\circ q}: V_n \to V'_n$ can be chosen to converge to those of $P: U \to U'$ and $P^{\circ q}: V \to V'$. It follows from \corref{nchar} that $P \in \NN(t)$. Continuity of $\Chi$ then shows $\Chi(P)=\Chi(P_\infty)$, which implies $P=P_\infty$ by injectivity of $\Chi$. 

\LEM \label{onelands}
At least one of the rays $R(x_k), R(y_k)$ lands at $z_k$.
\ENDLEM

\PROOF
Suppose $R(x_k), R(y_k)$ land at $\zeta^-, \zeta^+$ distinct from $z_k$. By \thmref{het}, $\zeta^\pm$ are parabolic and the Hausdorff limits of $\ov{R}_n(x_k)$ and $\ov{R}_n(y_k)$ are of the form $R(x_k) \cup \{ \zeta^- \} \cup \gamma^- \cup \{ z_k \}$ and $R(y_k) \cup \{ \zeta^+ \} \cup \gamma^+ \cup \{ z_k \}$. Here $\gamma^\pm$ are disjoint heteroclinic arcs connecting $\zeta^\pm$ to $z_k$ and residing in period $q$ parabolic basins $B^\pm$, where $B^-=B^+$ if and only if $\zeta^-=\zeta^+$. Consider the dynamical wakes $W_{i,n} \ (1 \leq i \leq q)$ and the sub-wake $W'_{k,n} \subset W_{k,n}$ of $P_n$. Define the ``wakes'' $W_i \ (1 \leq i \leq q)$ and $W'_k$ of $P$ as the interior of the Hausdorff limit of the closures $\ov{W_{i,n}}$ and $\ov{W'_{k,n}}$, respectively. It follows from the basic action of $P_n$ on its wakes that the $W_i$ are pairwise disjoint and we have the following diagram for the action of $P$: 
$$
W_k \sm \ov{W'_k} \stackrel{2:1}{\longrightarrow} W_{\sigma(k)} \stackrel{1:1}{\longrightarrow} \cdots \stackrel{1:1}{\longrightarrow} W_{\sigma^{q-1}(k)} \stackrel{1:1}{\longrightarrow} W_k.
$$
In particular, $P^{\circ q}: W_k \sm \ov{W'_k} \to W_k$ is a proper map of degree $2$ with a single critical point at $\om_2$. Moreover, the forward orbit of $\om_1$ never enters the union $\bigcup_{i=1}^{q} W_i$. We distinguish three cases: \vs

%%%%%%%%%%%%%%%%%%%%%%%%%%%%%%%%%%%%%%%%%%%
\begin{figure}[t]
	\centering
	\begin{overpic}[width=0.85\textwidth]{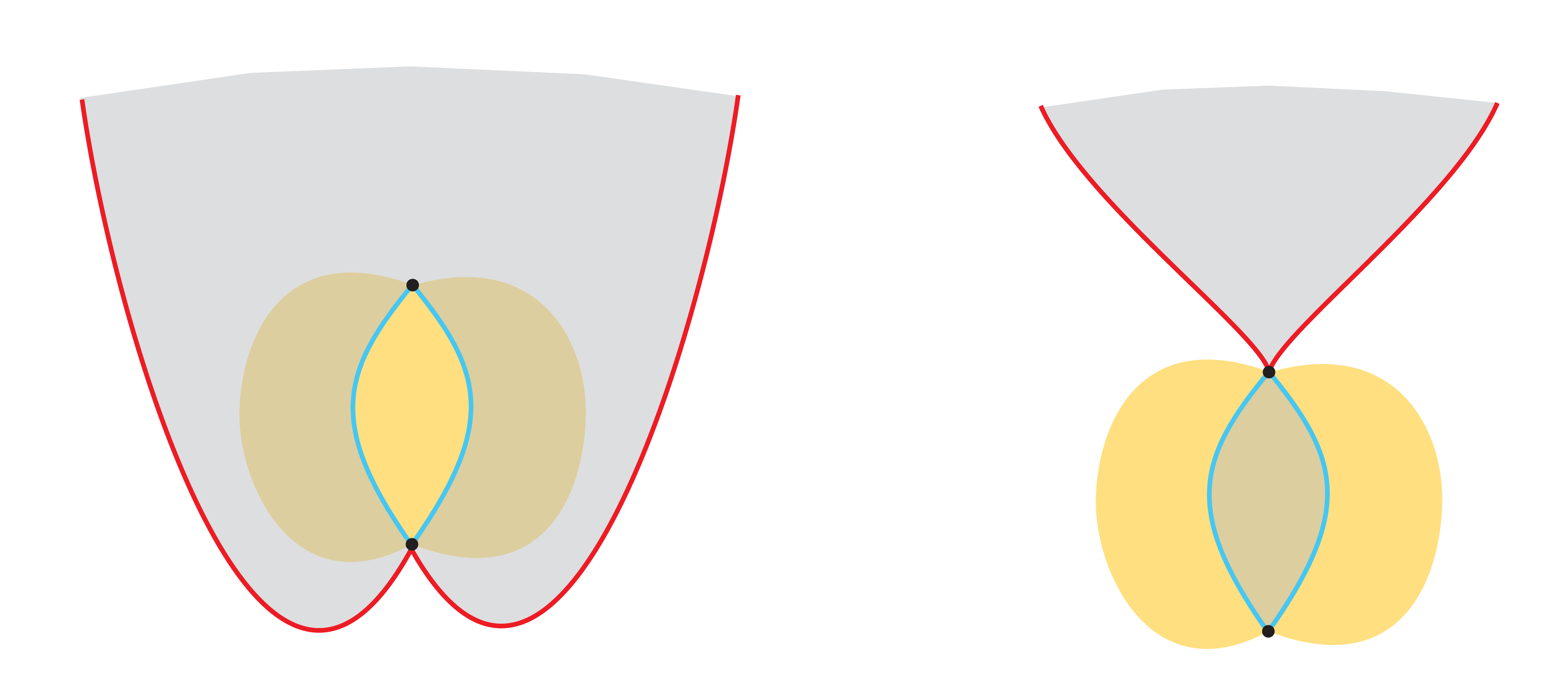}
		\put (25,33) {\footnotesize $W_k$}
		\put (34,26) {\footnotesize $B$}
		\put (25.4,18) {\footnotesize $\Om$}
		\put (28.9,12) {\tiny $\gamma^-$}
		\put (21.3,12) {\tiny $\gamma^+$}
		\put (32.5,18) {\footnotesize $D^-$}
		\put (17,18) {\footnotesize $D^+$}	
		\put (25.5,6.4) {\tiny $\zeta$}
		\put (25.3,27.5) {\tiny $z_k$}   
		\put (45,40) {\footnotesize $R(x_k)$}
		\put (1,40) {\footnotesize $R(y_k)$}
		\put (89.5,19.5) {\footnotesize $B$}
		\put (80,12) {\footnotesize $\Om$}		
		\put (80,33) {\footnotesize $W_k$}
		\put (93,39) {\footnotesize $R(x_k)$}
		\put (63,39) {\footnotesize $R(y_k)$}
		\put (83.7,7) {\tiny $\gamma^-$}
		\put (75.6,7) {\tiny $\gamma^+$}
		\put (87,12) {\footnotesize $D^-$}
		\put (72,12) {\footnotesize $D^+$}	
		\put (80,2) {\tiny $z_k$}
		\put (80.5,23) {\tiny $\zeta$}
	\end{overpic}
	\caption{\small The cases 2 (left) and 3 (right) in the proof of \lemref{onelands}. In either case $W_k$ is the region shaded in gray and $B$ is the parabolic basin of $\zeta$ shaded in yellow.}  
	\label{cases}
\end{figure}
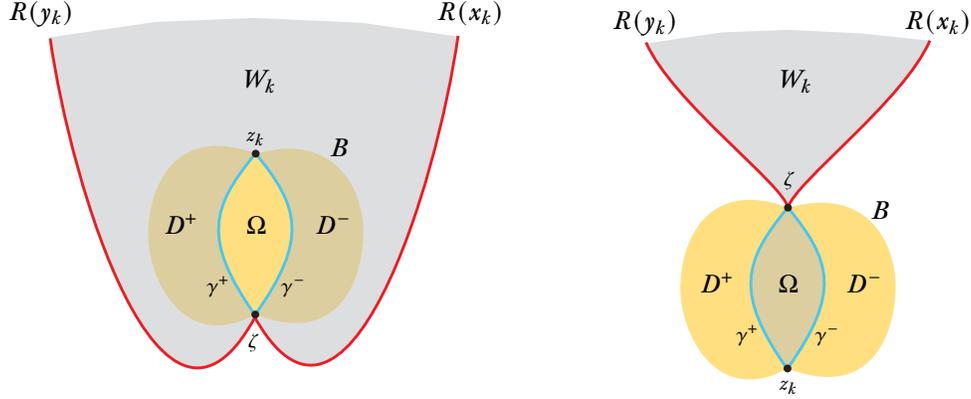
%%%%%%%%%%%%%%%%%%%%%%%%%%%%%%%%%%%%%%%%%%%

{\it Case 1.} $\zeta^- \neq \zeta^+ \Rightarrow B^- \cap B^+ = \es$. This is easily seen to be impossible since by \thmref{het} both critical points $\om_1$ and $\om_2$ of $P$ would have to be simultaneously contained in the (disjoint) orbits of $B^\pm$. \vs

{\it Case 2.} $\zeta^-=\zeta^+=\zeta \Rightarrow B^-=B^+=B$ and the domain $\Om$ bounded by $\ga^\pm$ is disjoint from $W_k$ (\figref{cases} left). The complement $B \sm (\ga^- \cup \ga^+)$ has three connected components, one being $\Om$ and the other two contained in $W_k$ and labeled $D^\pm$ having $\ga^\pm$ on their boundary. For integer $i$ modulo $q$, set $B_i:=P^{\circ i}(B)$ and let $\Om_i:=P^{\circ i}(\Om)$ be the domain in $B_i$ bounded by $\ga^\pm_i:=P^{\circ i}(\ga^\pm)$. Let $D_i^\pm$ denote the remaining components of $B_i \sm (\ga^-_i \cup \ga^+_i)$ having $\ga^\pm_i$ on their boundary. Evidently $D^-_i \cup D^+_i \subset W_{\sigma^i(k)}$ but $\Om_i \cap \bigcup_{j=1}^q W_j= \es$. Since $P^{\circ q}: \bd \Om \to \bd \Om$ is a homeomorphism, $P^{\circ q}: \Om \to \Om$ must be a conformal isomorphism, so $P$ maps each $\Om_i$ conformally onto $\Om_{i+1}$. On the other hand, $D^\pm_i$ maps conformally onto $D^{\pm}_{i+1}$ if $D^\pm_i$ does not contain a critical point, and it maps as a degree $2$ branched covering onto $B_{i+1}$ otherwise. By \thmref{het}, both $\om_1, \om_2$ belong to some iterated images of $D^\pm$. Since $\om_2$ is in $W_k$, it must belong to one of  $D^\pm$, say $\om_2 \in D^-$. In that case $\om_1 \in D^+_i$ for some $i$. This implies $\om_1 \in W_{\sigma^i(k)}$, which is a contradiction. \vs   

{\it Case 3.} $\zeta^-=\zeta^+=\zeta \Rightarrow B^-=B^+=B$ and the domain $\Om$ bounded by $\ga^\pm$ is contained in $W_k$ (\figref{cases} right). We can define $B_i, \Om_i, D_i^\pm$ as in {\it Case 2} except that now $\Om_i \subset W_{\sigma^i(k)}$ and $(D^-_i \cup D^+_i) \cap \bigcup_{j=1}^q W_j = \es$. Again by \thmref{het}, both $\om_1, \om_2$ belong to the iterated images of $D^\pm$. Assuming without loss of generality that $\om_1 \in D^-_i$ for some $i$, it follows that $\om_2 \in D^+_j$ for some $j$. This implies $D^+_j \cap W_k \neq \es$, which is a contradiction. 
\ENDPROOF

\LEM \label{bothland}
Both rays $R(x_k),R(y_k)$ land at $z_k$, hence $P \in \LL(t)$.
\ENDLEM

\PROOF
By \lemref{onelands} we may assume without loss of generality that $R(y_k)$ lands at $z_k$. Suppose $R(x_k)$ lands at a point $\zeta \neq z_k$. Then, by \thmref{het}, $\zeta$ is parabolic and the Hausdorff limit of $\ov{R}_n(x_k)$ is of the form $R(x_k) \cup \{ \zeta \} \cup \gamma \cup \{ z_k \}$, where $\gamma$ is a heteroclinic arc connecting $\zeta$ to $z_k$ within a period $q$ parabolic basin $B$. Notice that in this case $z_k$ is repelling and the Hausdorff limit of $\ov{R}_n(y_k)$ is $\ov{R}(y_k)$. Let $p=q/r$ be the period of $\zeta$ under $P$. Then the multiplier $(P^{\circ p})'(\zeta)$ is a primitive $r$th root of unity and the multiplicity of $\zeta$ as a fixed point of $P^{\circ q}$ is $r+1$. Thus, the parabolic point $\zeta$ bifurcates into a nearby fixed point $\zeta_n$ and an $r$-cycle $C_n$ of the map $P_n^{\circ p}$. According to \cite[Theorem 5.1]{PZ4}, for all sufficiently large $n$ the $P_n^{\circ q}$-invariant ray $R_n(x_k)$ which enters and eventually exists a fixed small disk centered at $\zeta$ must separate the points of $C_n \cup \{ \zeta_n \}$ within that disk. This implies that the ray pair $R_n(x_k), R_n(y_k)$ that co-land on $z_{k,n}$ separate these $r+1$ points. In particular, there is a point $\zeta'_n \in C_n \cup \{ \zeta_n \}$ that belongs to $W_{k,n} \sm \ov{W'}_{k,n}$, hence to the peripheral filled Julia set of $P_n$. Under straightening, $\zeta'_n$ must correspond to the $\alpha$-fixed point of the peripheral renormalization since $z_{k,n}$ corresponds to its $\beta$-fixed point. It follows from \lemref{same} that $\chiper(P_\infty)=\lim_{n \to \infty} \chiper(P_n)= 1/4$. This contradicts the hypothesis $\Chi(P_\infty) \in \check{\MM}_t \times \check{\MM}_0$. 
\ENDPROOF

\PROOF[Proof of \thmref{proper}] 
In view of \lemref{bothland}, it remains to prove that $z_k$ is repelling. For simplicity let $(c_n,\nu_n):=\Chi(P_n) \to (c,\nu):=\Chi(P_\infty)$. Assuming $z_k$ is parabolic, we reach a contradiction by showing that either $c=c(t_i)$ for some $1 \leq i \leq q$, or $\nu=1/4$. \vs

First suppose $z_k$ is parabolic but the dynamical wake $W_k$ of $P$ defined by the co-landing pair $R(x_k),R(y_k)$ does not contain an immediate attracting basin of $z_k$. Then the thickening procedure gives a central quadratic-like restriction $P: U \to U'$ (see \S \ref{cpren}) with connected Julia set which is hybrid equivalent to some quadratic polynomial $Q$. Under this hybrid conjugacy the co-landing orbit $\{ z_1, \ldots, z_q \}$ corresponds to an  orbit $\{ w_1, \ldots, w_q \}$ of $Q$ which is still parabolic. For large $n$ there is a component $U_n$ of $P_n^{-1}(U')$ containing $\om_{1,n}$ such that the restriction $P_n: U_n \to U'$ is quadratic-like with connected Julia set, necessarily hybrid equivalent to $Q_{c_n}$. Continuity of the straightening map $\chicen$ then shows that $Q=Q_c$. By \lemref{renorm}(i), the rays $R_c(t_i)$ land at the parabolic point $w_i$ for all $1 \leq i \leq q$. This implies $c=c(t_i)$ for some $i$, a contradiction. \vs

Similarly, if $z_k$ is parabolic with no immediate attracting basin outside $W_k$, the thickening procedure gives a peripheral quadratic-like restriction $P^{\circ q}: V \to V'$ with connected Julia set which is hybrid equivalent to some quadratic polynomial $Q$. Evidently this quadratic has a parabolic fixed point of multiplier $1$, so $Q=Q_{1/4}$. For large $n$ there is a component $V_n$ of $P_n^{-q}(V')$ containing $\om_{2,n}$ such that the restriction $P^{\circ q}_n: V_n \to V'$ is quadratic-like with connected Julia set, necessarily hybrid equivalent to $Q_{\nu_n}$. By continuity of the straightening map $\chiper$ we have $\nu_\infty=\lim_{n \to \infty} \nu_n =1/4$, which is again a contradiction. \vs

Finally, suppose $z_k$ is parabolic with a pair of immediate attracting basins on both sides of the ray pair $R(x_k),R(y_k)$: a central basin $B^-$ outside $W_k$ and a peripheral basin $B^+$ inside $W_k$. (Notice that at this point we do not know if the Hausdorff limits of $\ov{R}_n(x_k)$ and $\ov{R}_n(y_k)$ are $\ov{R}(x_k)$ and $\ov{R}(y_k)$, so a priori $B^\pm$ may contain homoclinic arcs in these Hausdorff limits.) Now there are no quadratic-like restrictions near $z_k$ in the classical sense, but there is a restriction $P^{\circ q}: \Om \to \Om'$ between suitable Jordan domains containing $z_k$ that is a degree $2$ parabolic-like map in the sense of Lomonaco, with the required invariant arcs being the intersections of $R(x_k),R(y_k)$ with $\Om'$ (see \cite{Lo} for details). In fact, we can choose $\Om'$ to be a domain bounded by an equipotential of sufficiently small potential in $W_k$ starting on $R(x_k)$ and ending on $R(y_k)$, completed with an arc in $\C \sm W_k$ that stays sufficiently close to $z_k$. We can then take $\Om$ to be the connected component of $P^{-q}(\Om')$ containing $z_k$. The filled Julia set of the parabolic-like map $P^{\circ q}: \Om \to \Om'$ is the set of points in $K_P$ whose forward $P^{\circ q}$-orbit remain in $\ov{W_k} \sm W'_k$ (in particular, it contains the basin $B^+$). By the straightening theorem for parabolic-like maps, $P^{\circ q}|_\Om$ is hybrid equivalent to a restriction of a unique quadratic rational map of the form $g_A: z \mapsto z+z^{-1}+A$, with $z_k$ corresponding to the parabolic fixed point $\infty$ of $g_A$. Since $z_k$ is a degenerate parabolic point of $P^{\circ q}$, the same must be true of $\infty$ for $g_A$, which implies $A=0$. It follows that $P^{\circ q}|_\Om$ is hybrid conjugate to the restriction of $g_0: z \mapsto z+z^{-1}$ to a neighborhood of the closed left (or right) half-plane and in particular, $(\ov{W_k} \sm W'_k) \cap K_P \subset \Om$ does not contain any fixed point of $P^{\circ q}$ other than $z_k$. Now consider the unique fixed point $\zeta_n$ of $P_n^{\circ q}$ in $W_{k,n} \sm \ov{W'_{k,n}}$ corresponding to the $\alpha$-fixed point of the peripheral renormalization. Let $\zeta$ be an accumulation point of $\{ \zeta_n \}$. Then $\zeta$ belongs to any Hausdorff limit of $\ov{W_{k,n}} \sm W'_{k,n}$ as $n \to \infty$. Such a limit is necessarily contained in $\ov{W_k} \sm W'_k$ union the central basin $B^-$ at $z_k$ and its partner at the co-landing point $z'_k$ of the rays $R(y'_k), R(x'_k)$ (to account for the possibility of homoclinic arcs in $B^-$). Since $\zeta$ is a fixed point of $P^{\circ q}$, we must have $\zeta \in \ov{W_k} \sm W'_k$. It follows from what we established above that $\zeta=z_k$. This proves $\zeta_n \to z_k$, so $(P_n^{\circ q})'(\zeta_n) \to 1$ and therefore $\nu_\infty = \lim_{n \to \infty} \nu_n = 1/4$ by \lemref{same}, which is a contradiction.
\ENDPROOF 

\REM \label{hetconj}
The proof of \lemref{onelands} shows that for every cubic $P \in \ov{\LL(t)} \sm \LL(t)$ exactly one of the rays $R(x_k),R(y_k)$ lands at a parabolic fixed point $\zeta$ of $P^{\circ q}$ of multiplier $1$. The other ray necessarily lands at a repelling fixed point of $P^{\circ q}$ on the boundary of an immediate basin $B$ of $\zeta$ having exact period $q$. Moreover, the first-return map $P^{\circ q}: B \to B$ is proper of degree $3$ or $4$ depending on whether both $\om_1,\om_2$ belong to $B$ or they belong to different basins in the orbit of $B$. In Milnor's terminology \cite{M4}, these cases correspond respectively to $P$ being on the boundary of an ``adjacent'' or ``bi-transitive'' hyperbolic component in $\CC(3)$. The failure of co-landing of $R(x_k),R(y_k)$ is caused by the creation of a heteroclinic arc in the Hausdorff limit of one of these rays. This mechanism can be forced, for example, by taking the limit of suitable sequences $P_n \in \NN(t)$ with $\chicen(P_n) \to c(t_i)$ for some $1 \leq i \leq q$ and $\chiper(P_n) \to 1/4$. It seems reasonable to speculate that all cubics of the form described above are in fact in $\ov{\LL(t)} \sm \LL(t)$. 
\ENDREM

\section{Surjectivity of straightening}\label{surg}

The final piece in the proof of \thmref{B} is the following 

\THM [Surjectivity of $\Chi$] \label{surject}
For every pair $(c,\nu) \in \check{\MM}_t \times \check{\MM}_0$ there is a cubic polynomial $P \in \NN(t)$ such that $\chicen(P)=c$ and $\chiper(P)=\nu$. 
\ENDTHM 

The proof of this theorem is achieved in three steps that can be loosely described as follows: \vs

{\it Step 1.} Start from the quadratic polynomial $Q_c$ and the unique cubic polynomial $P_a \in \LL_0(t)$ which is peripherally renormalizable around the second critical point $\om_2$ with a quadratic-like restriction hybrid equivalent to $Q_\nu$ (\corref{L0}). Construct a new holomorphic map $f:\C \sm \ov{\La} \to \C$ by changing the external class of $Q_c$ (i.e., the map $Q_0:z \mapsto z^2$) to the external class of the quadratic-like map obtained from the central renormalization of $P_a$ around its superattracting basin $B_a$. Here $\La$ is the ``hole'' corresponding to the component $B'_a$ of $P_a^{-1}(B_a)$ other than $B_a$. \vs

{\it Step 2.} Change the conformal structure on $\La$ and all its iterated preimages under $f$ if necessary to create a large enough modulus necessary for quasiconformal interpolation. \vs
	
{\it Step 3.} Restrict to a topological disk bounded by a suitable equipotential curve and perform a quasiconformal interpolation on the sub-wake $W_k' \subset W_k$ of $f$ containing $\La$ to obtain a quasiregular polynomial-like map of degree $3$ with a central renormalization hybrid equivalent to $Q_c$ and a peripheral renormalization hybrid equivalent to $Q_\nu$. This quasiregular map has an invariant conformal structure of bounded dilatation, so we can apply the straightening theorem to obtain the desired cubic polynomial $P$. \vs

The full details of the above steps are given in \S \ref{mec} through \S \ref{pfsur}. The short section \S \ref{qcint} recalls some definitions and relevant known facts on quasiconformal interpolation that are needed in the proof of \thmref{surject}.

\subsection{Dynamical cylinders and interpolation on sectors} \label{qcint}

We begin with some generalities. Let $f$ be a degree $d$ polynomial with connected filled Julia set $K$. Suppose $z_0$ is a repelling point of $f$ of period $p$ and combinatorial rotation number $s/r$. The external rays landing at $z_0$ have a common period $q=pr$ under $f$ and fall into $N$ cycles for some $1 \leq N \leq d-1$, forming $Nr$ rays altogether. There are $Nr$ connected components of $K \sm \{ z_0 \}$ whose intersections with a small neighborhood of $z_0$ are permuted around by $f^{\circ p}$ with combinatorial rotation number $s/r$. We may label these rays and components as $R_1 \ldots, R_{Nr}$ and $K_1, \ldots, K_{Nr}$ in positive cyclic order, so $K_j$ lies between $R_j$ and $R_{j+1}$ (taking indices modulo $Nr$). \vs

Now take a small Jordan domain $D$ containing $z_0$ in which $f$ is linearizable. Then $D \Subset D':=f^{\circ q}(D)$ and the biholomorphism $f^{\circ q}:D \to D'$ is conjugate to $z \mapsto \la z$, where $\la:=(f^{\circ q})'(z_0)$. The quotient $T:=(D \sm \{ z_0 \})/f^{\circ q}$ is a complex torus conformally isomorphic to $(\D \sm  \{ 0 \})/\langle z \mapsto \la z \rangle$. In particular, the conformal isomorphism class of $T$ is independent of the choice of $D$. For $1 \leq j \leq Nr$ the ray tails $R_j \cap D$ and the intersections $K_j \cap D$ project to closed curves $\tilde{R}_j$ and compact sets $\tilde{K}_j$ in $T$. The connected component of $T \sm \bigcup_{j=1}^{Nr} \tilde{K}_j$ containing $\tilde{R}_j$ is a conformal cylinder $C_j$ of modulus $\Mod(C_j)=\pi/(q \log d)$ having $\tilde{R}_j$ as its core geodesic. All the $C_j$ (hence the $\tilde{R}_j$) are in the same non-trivial isotopy class in $T$.  \vs 

%%%%%%%%%%%%%%%%%%%%%%%%%%%%%%%%%%%%%%%%%%%
\begin{figure}[t]
\centering
\begin{overpic}[width=\textwidth]{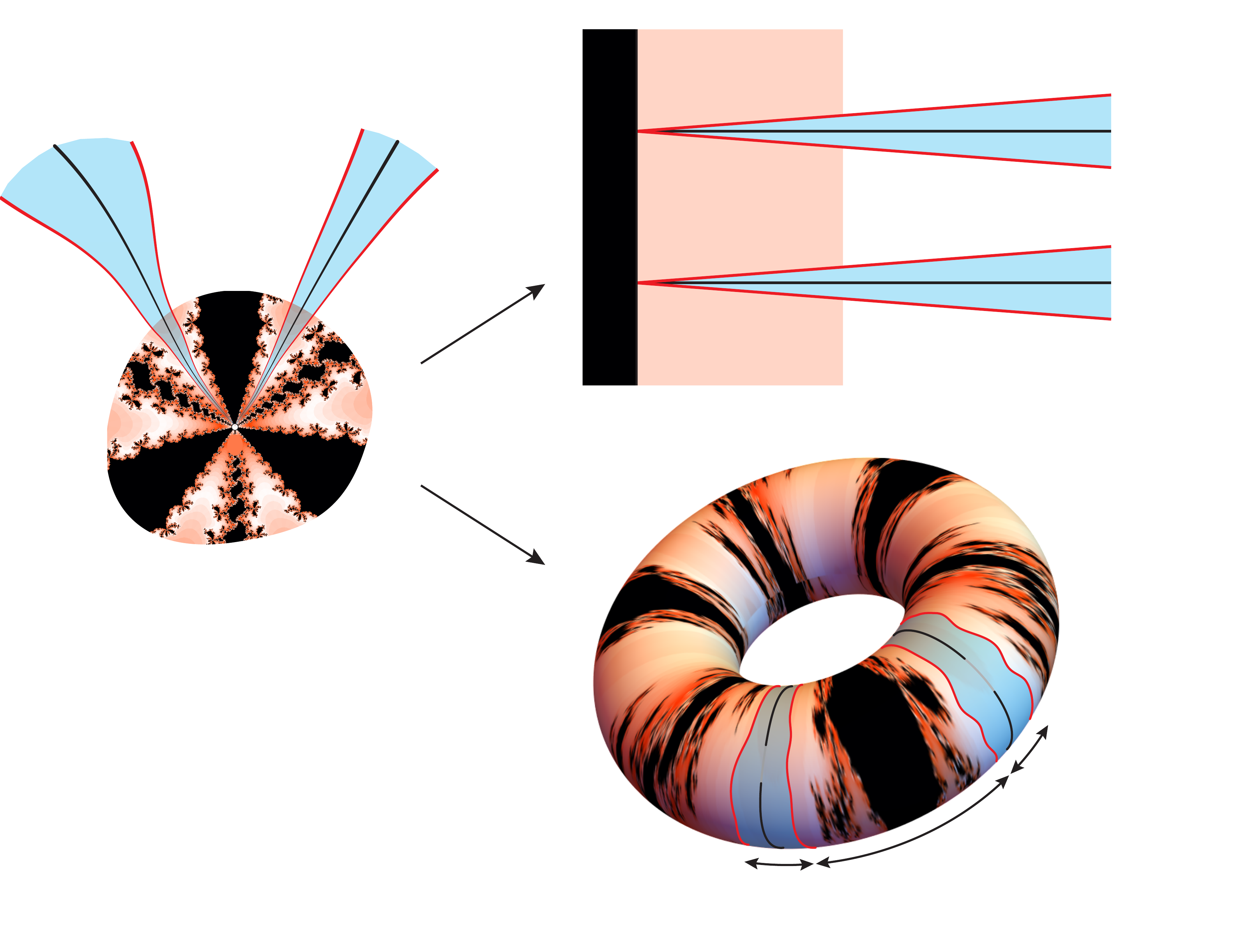}
\put (18.2,40.5) {\tiny $z_0$}
\put (32,66) {\footnotesize $R_j$}
\put (36,63) {\footnotesize {\color{red} $R^-_j$}}
\put (26,67) {\footnotesize {\color{red} $R^+_j$}}
\put (2,66) {\footnotesize $R_{j+1}$}
\put (10,66) {\footnotesize {\color{red} $R^-_{j+1}$}}
\put (-3,59) {\footnotesize {\color{red} $R^+_{j+1}$}}
\put (17,49) {\small {\color{white} $K_j$}}
\put (6,43) {\small $D$}
\put (34,51) {\footnotesize $\log \phi$}
\put (31,32) {\footnotesize quotient}
\put (32.3,29.5) {\footnotesize by $f^{\circ q}$}
\put (90,53.5) {\footnotesize $\eta_j$}
\put (90,50) {\footnotesize {\color{red} $\eta^-_j$}}
\put (90,57) {\footnotesize {\color{red} $\eta^+_j$}}
\put (90,65.5) {\footnotesize $\eta_{j+1}$}
\put (90,62) {\footnotesize {\color{red} $\eta^-_{j+1}$}}
\put (90,69) {\footnotesize {\color{red} $\eta^+_{j+1}$}}
\put (70,17) {\small {\color{white} $\tilde{K}_j$}}
\put (85,31) {\small $T$}
\put (61,4) {\small $H_j$}
\put (73,6.5) {\small $G_j$}
\put (83.5,14) {\small $H_{j+1}$}
\put (60,14.3) {\footnotesize $\tilde{R}_j$}
\put (77,21.5) {\footnotesize $\tilde{R}_{j+1}$}
\end{overpic}
\caption{\small Oblique rays landing at a repelling point and various cylinders in the quotient torus.}  
\label{torus}
\end{figure}
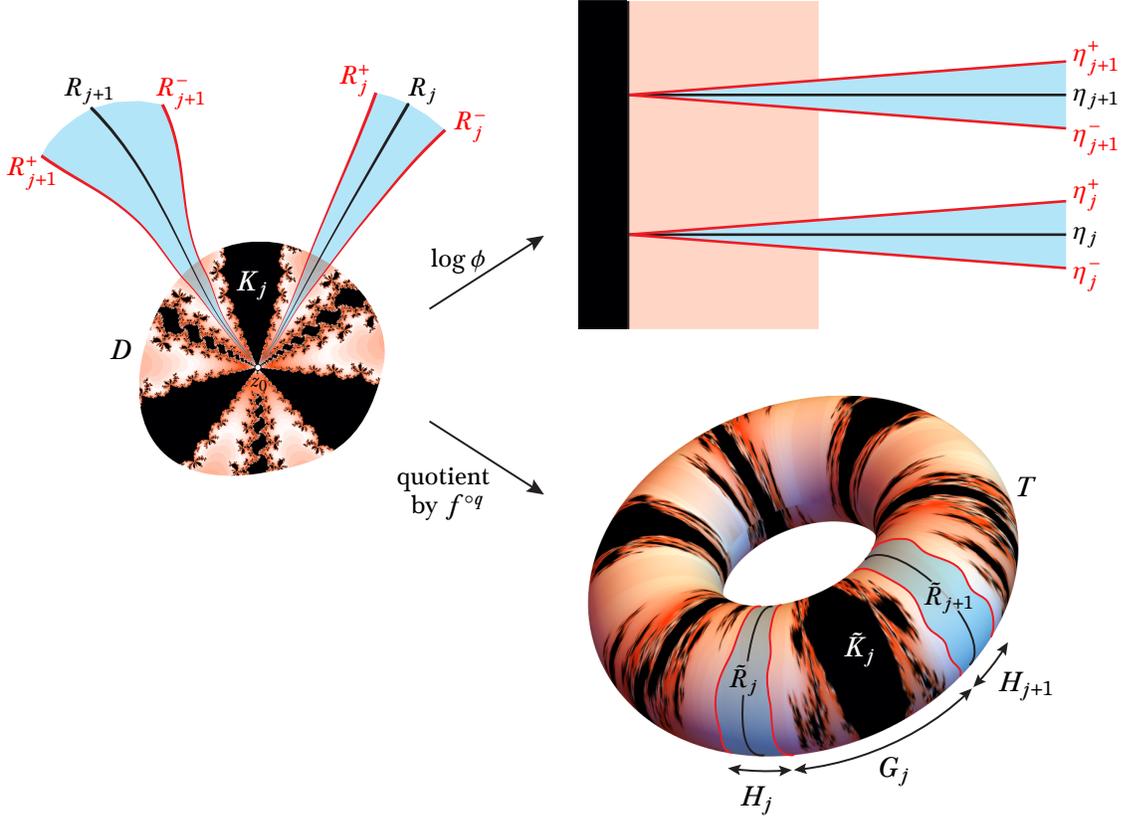
%%%%%%%%%%%%%%%%%%%%%%%%%%%%%%%%%%%%%%%%%%%    

For the purpose of interpolation, it is often useful to replace the $C_j$ by smaller cylinders with smooth (in fact, real analytic) boundary curves as follows. Suppose $0 \leq \theta_j <1$ is the angle of $R_j$, which is an integer multiple of $1/(d^q-1)$. Let $\phi: \C \sm K \to \C \sm \ov{\D}$ be the \Bottcher coordinate of $f$. The map $\log \phi: \C \sm \ov{\D} \to \H^+=\{ z: \re(z)>0 \}$ lifts $R_j$ to the horizontal ray $\eta_j(t) =t+2\pi \ii \, \theta_j$ for $t>0$. The rays 
$$
\eta_j^\pm(t)= t + 2\pi \ii \, (\theta_j \pm \varsigma t) 
$$ 
map back to the arcs $R_j^\pm = (\phi^{-1}\circ \exp)(\eta_j^{\pm})$ which we call the {\bit oblique rays of slope} $\pm \varsigma$ at angle $\theta_j$ (see \figref{torus}). Evidently $R_j^\pm$ land at $z_0$ and $R_j^\pm=f^{\circ q}(R_j^\pm)$. Moreover, the $2Nr$ tails $R_j^\pm \cap D$ for $1 \leq j \leq Nr$ are pairwise disjoint. In fact, $R_j^\pm \cap D_n$ are pairwise disjoint for all large $n>0$, where $D_n$ is the connected component of $f^{-nq}(D)$ containing $z_0$, and $f^{\circ nq}: D_n \to D$ is a biholomorphism. The sector in $D$ delimited by $R_j^\pm$ projects to a cylinder $H_j \subset C_j$ having $\tilde{R}_j$ as its core geodesic. It is easy to see that $\Mod(H_j)$ tends to $0$ as $\varsigma \to 0$ and it tends to $\Mod(C_j)=\pi/(q \log d)$ as $\varsigma \to +\infty$.  \vs  

We are also interested in the complementary sectors in $D$ that meet $K$. Following \cite{EY} we consider the invariant   
$$
\Mod_{z_0}(K_j)=\inf\{ \Mod(C): C \supset \tilde{K}_j \} = \sup \{ \Mod(C): C \subset \tilde{K}_j \}
$$
where the cylinders $C \subset T$ are taken in the same isotopy class as $C_j$. This invariant is not affected by quasiconformal deformations that keep the conformal structure on $K_j$ unchanged. Moreover, $\Mod_{z_0}(K_j)>0$ if and only if $z_0$ is on the boundary of a Fatou component of $f$ in $K_j$. The sector in $D$ delimited by $R^+_j$ and $R^-_{j+1}$ projects to a cylinder $G_j \subset T$ which contains $\tilde{K}_j$. As $\varsigma \to \infty$, $\Mod(G_j) \to \Mod_{z_0}(K_j)$. As $\varsigma \to 0$, $\Mod(G_j)$ tends to the modulus of the cylinder bounded by $\tilde{R}_j$ and $\tilde{R}_{j+1}$, which is at least 
$$
\frac{\pi}{2q \log d}+\Mod_{z_0}(K_j)+\frac{\pi}{2q \log d}=\Mod_{z_0}(K_j)+\frac{\pi}{q \log d}
$$ 
by the \Grotszch inequality. \vs 

We note that all the above constructions can be carried out for polynomial-like maps hybrid equivalent to $f$, or more generally for local holomorphic maps obtained from $f$ by a conjugation near $z_0$.  \vs   

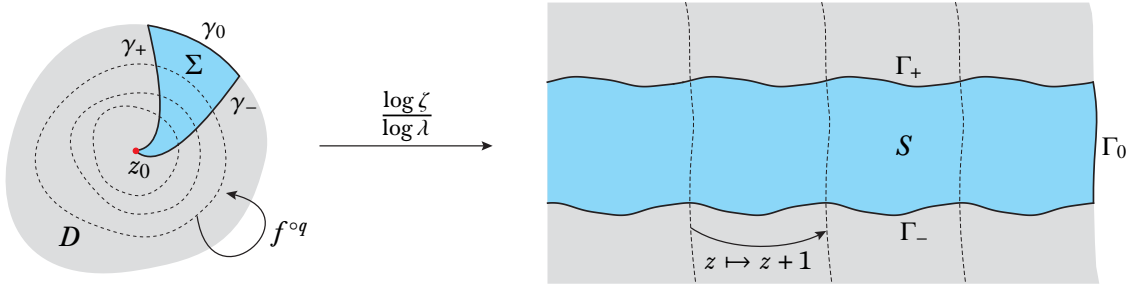
\begin{figure}[t]
	\centering
	\begin{overpic}[width=\textwidth]{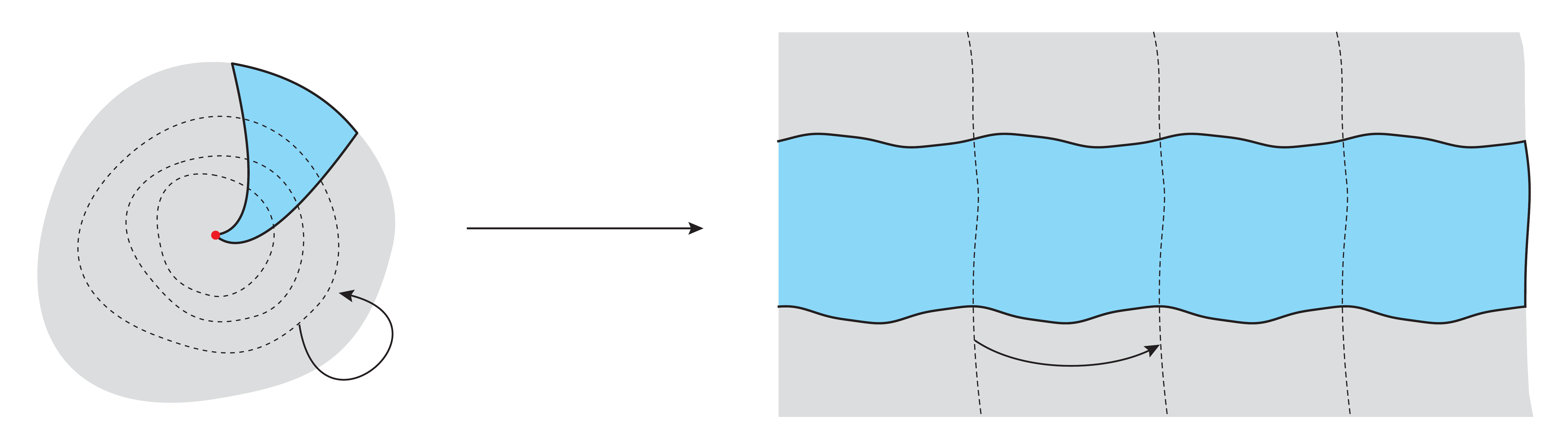}
		\put (18,20) {\small $\Si$}
		\put (13,11.3) {\small $z_0$}
		\put (7,5) {\small $D$}
		\put (12.5,22) {\footnotesize $\ga_+$}
		\put (22,17) {\footnotesize $\ga_-$}
		\put (19.5,23.5) {\footnotesize $\ga_0$}
		\put (35,15.7) {$\frac{\log \zeta}{\log \la}$}
		\put (80,13) {\small $S$}
		\put (80,20) {\footnotesize $\Ga_+$}
		\put (80.4,5.5) {\footnotesize $\Ga_-$}
		\put (98,13) {\footnotesize $\Ga_0$}
      \put (63.3,3) {\footnotesize $z \mapsto z+1$}
      \put (25.5,5.5) {\footnotesize $f^{\circ q}$}
	\end{overpic}
	\caption{\small An invariant sector $\Si$ and its half-strip model $S$.}  
	\label{strip}
\end{figure}

Let us now discuss quasiconformal interpolation on sectors whose idea goes back to the thesis of B. Bielefeld \cite{B}. Consider $f$ as above and let $\Si$ be a sector in $D$ delimited by two smooth arcs $\ga_\pm$ landing at $z_0$. Thus, $\bd \Si$ consists of $\ga_\pm$ and a third arc $\ga_0 \subset \bd D$. We assume $\Si$ is invariant in the sense $f^{\circ q}(\Si) \cap D=\Si$. Denote by $\zeta: (D,z_0) \to (\D,0)$ the (unique up to a rotation) linearizing map that conjugates $f^{\circ q}$ to $z \mapsto \la z$. There is a univalent branch of $z \mapsto \log \zeta(z)/\log \la$ that maps $\Si$ to a half-strip $S$ and conjugates $f^{\circ q}$ to the unit translation $z \mapsto z+1$. We call $S$ the {\bit half-strip model} of $\Si$. Note that $\bd S$ consists of two ``horizontal'' arcs $\Ga_\pm$ that are invariant under $z \mapsto z-1$, joined by a ``vertical'' arc $\Ga_0$ (see \figref{strip}). \vs

%%%%%%%%%%%%%%%%%%%%%%%%%%%%%%%%%%%%%%%%%%%
%%%%%%%%%%%%%%%%%%%%%%%%%%%%%%%%%%%%%%%%%%%    
   
Now suppose we have another holomorphic map $\hat{f}$ with a repelling fixed point $\hat{z}_0$ and invariant sector $\hat{\Si} \subset \hat{D}$ with boundary arcs $\hat{\ga}_\pm, \hat{\ga}_0$ and the corresponding half-strip model $\hat{S}$ with boundary arcs $\hat{\Ga}_\pm, \hat{\Ga}_0$. Consider a homeomorphism $h: \bd \Si \to \bd \hat{\Si}$, which maps each boundary arc of $\Si$ to the corresponding boundary arc of $\hat{\Si}$. We say that the restrictions $h: \ga_\pm \to \hat{\ga}_\pm$ are {\bit near-linear} if they are smooth and the induced maps $H: \Ga_\pm \to \hat{\Ga}_\pm$ have uniform $C^1$ distance from a translation, i.e., there is a constant $M>1$ such that 
$$
|H(z(t))-z(t)|<M \qquad \text{and} \qquad M^{-1} \leq \frac{(H \circ z)'(t)}{z'(t)} \leq M
$$
where $t \mapsto z(t)$ is any $C^1$ parametrization of $\Ga_\pm$. Observe that if $h: \ga_\pm \to \hat{\ga}_\pm$ conjugate $f^{\circ q}$ to $\hat{f}^{\circ q}$, they are automatically near-linear since in this case the induced maps $H: \Ga_\pm \to \hat{\Ga}_\pm$ commute with the unit translation. \vs 

We can now state the following two results that we will use in \S \ref{extf}. For a detailed account we refer the reader to the book \cite{BF}, especially Proposition 2.26 and Lemmas 2.15, 8.4, and 8.5: 

\LEM \label{bf1}
Let $f,\Si$ and $\hat{f},\hat{\Si}$ be as above and $h: \Si \to \hat{\Si}$ be the unique conformal isomorphism which maps each boundary arc of $\Si$ to the corresponding boundary arc of $\hat{\Si}$. If $\Mod(\Si/f^{\circ q})=\Mod(\hat{\Si}/\hat{f}^{\circ q})$, then the restrictions $h: \ga_\pm \to \hat{\ga}_\pm$ are near-linear.
\ENDLEM

\LEM [Quasiconformal interpolation on sectors] \label{bf2}
Let $f,\Si$ and $\hat{f},\hat{\Si}$ be as above and $h: \bd \Si \to \bd \hat{\Si}$ be a homeomorphism which maps each boundary arc of $\Si$ diffeomorphically to the corresponding boundary arc of $\hat{\Si}$. If $h: \ga_\pm \to \hat{\ga}_\pm$ are near-linear, then $h$ extends to a quasiconformal homeomorphism $\Si \to \hat{\Si}$. 
\ENDLEM 

\subsection{Modifying the external class of $Q_c$}\label{mec}

Let us return to the three steps of the proof of \thmref{surject} outlined in the beginning of this section. Step 1 can in fact be performed for any choice of $c \in\Mbrot$ and $P_a \in \LL_0(t)$. As before let $B_a \cong \D$ be the immediate basin of attraction of $0$ for $P_a$ and let $B'_a$ denote the connected component of $P_a^{-1}(B_a)$ other than $B_a$. Let $P_a|_U: U \to U'$ be a central quadratic-like restriction with connected filled Julia set $\ov{B_a}$. The restriction $h_a: \C \sm P_a^{-1}(\ov{B_a}) \to \C \sm \ov{B_a}$ of $P_a$ represents the ``external class'' of $P_a|_U$. Applying the general straightening theorem of Douady and Hubbard (in the interpretation of Lyubich as described in \cite{Ly}), we obtain a holomorphic map $f: \C \sm \ov{\La} \to \C$ with a quadratic-like restriction hybrid equivalent to $Q_c$ and with the same external class $h_a$ as $P_a|_U$. Let us explain the construction of $f$. This will also allow us to setup the notation. \vs

The fundamental annulus $\ov{U'} \sm U$ is contained in the domain of $h_a$. Choose a quadratic-like restriction $Q_c: V \to V'$, say with $V$ and $V'$ bounded by equipotential curves. Choose a diffeomorphism $\eta: \ov{V'}\sm V \to \ov{U'}\sm U$, with $\eta \circ Q_c = P_a \circ \eta$ on $\bd V$. Extend $\eta$ to a quasiconformal conjugacy $\eta: V'\sm K_c \to U' \sm \ov{B_a}$ by iterated lifting under the dynamics. Define the $Q_c$-invariant conformal structure $\mu$ on $V'$ as follows: On $V'\sm K_c$ let $\mu$ be the pull-back under $\eta$ of the standard conformal structure on $U' \sm \ov{B_a}$, and on $K_c$ let $\mu$ be the standard conformal structure. Then $\mu$ has bounded dilatation, so there is a quasiconformal homeomorphism $\zeta: V' \to \D$ which integrates $\mu$ and satisfies $\zeta(0) = 0$. Define a Riemann surface $\YY$ whose underlying set is the disjoint union $V' \sqcup \, (\Chat \sm \ov{B_a})$ modulo the smallest equivalence relation which identifies $z \in V' \sm K_c$ with $\eta(z)\in U' \sm \ov{B_a}$. The topology of $\YY$ is induced by the natural projection and its analytic atlas consists of $\zeta$ on $V'$ and the identity map on $\Chat \sm \ov{B_a}$. Notice that the change of coordinates $\zeta\circ\eta^{-1}$ is biholomorphic by the construction. Evidently $\YY$ is compact and simply connected. Take the uniformization $\Phi: \YY \to \Chat$ normalized by the conditions $\Phi(0)=0, \Phi(\infty)=\infty$ and $\Phi\circ Q_c\circ\Phi^{-1}(z) = \text{const.} + z^2 + O(z^3)$ near the origin. Then the restriction $\Phi: \Chat \sm \ov{B_a} \to \Chat \sm \Phi(K_c)$ is biholomorphic. Finally, set $\La := \Phi(B'_a)$ and define $f: \C \sm \ov{\La} \to \C$ by 
$$
f:= \begin{cases}
\Phi \circ Q_c \circ \Phi^{-1} & \quad \text{on} \ \Phi(V) \\
\Phi \circ h_a \circ \Phi^{-1} & \quad \text{on} \ \Phi(\C\sm P_a^{-1}(\ov{B_a})).
\end{cases}
$$ 
Evidently $f$ is well defined and holomorphic since $\eta$ conjugates $Q_c$ to $h_a=P_a$ in $V \sm K_c$. The restriction $f:\Phi(V) \to \Phi(V')$ is a quadratic-like map hybrid equivalent to $Q_c$ and with external class given by $h_a$. This completes the construction of $f$. \vs

We define the ``filled Julia set'' of $f$ by $K_f := \Phi(K_c) \cup \Phi(K_a \sm \ov{B_a})$. The Green's function in the basin of infinity $\C \sm K_a$ pushes forward under $\Phi$ to yield the ``Green's function'' of $f$ in $\C \sm K_f$. This allows us to freely talk about ``equipotentials'' and ``external rays'' of $f$. We use the notation $R_f(\theta)$ for the external ray $\Phi(R_a(\theta))$.     

\subsection{Basic constructions for $f$} \label{bcon}

Fix a pair $(c,\nu) \in \check{\MM}_t \times \check{\MM}_0$. In the construction of \S \ref{mec}, let $P_a$ be the unique cubic in $\LL_0(t)$ with $\chiper(P_a)=\nu$ (\corref{L0}). As usual, let $0 \leq t_1 < \cdots t_k=t < \cdots < t_q <1$ be the $\m_2$-orbit of $t$ with combinatorics $\sigma \in S_q$ and simulating $\m_3$-orbits $\OO_k=\{ x_1, \ldots, x_q \}, \OO_{k-1}=\{ y_1, \ldots, y_q \}$. As in \S \ref{cpren}, for each $1 \leq i \leq q$ there is a wake $W_a(t_i)$ bounded by the pair $R_a(x_i), R_a(y_i)$ that co-land at the repelling point $z_a(t_i) \in \bd B_a$, and the sub-wake $W'_a(t_k)$ of the critical wake $W_a(t_k)$ bounded by the pair $R_a(y'_k), R_a(x'_k)$ that co-land at $z'_a(t_k) \in \bd B'_a$ with $P_a(z'_a(t_k))=P_a(z_a(t_k))=z_a(t_j)$. Here $j:=\sigma(k)$ and  $x'_k=x_k+1/3, y'_k=y_k-1/3$, so $x_k<y'_k<1/2<x'_k<y_k$. Let $W_i := \Phi(W_a(t_i))$ for $1\leq i \leq q$ and $W'_k := \Phi(W'_a(t_k))$. For each $1 \leq i \leq q$ denote by $z_i$ the common landing point of the pair $R_f(x_i), R_f(y_i)$ which bound $W_i$, and let $z'_k \in \bd{\La}$ be the common landing point of the pair $R_f(y'_k),R_f(x'_k)$ which bound the sub-wake $W'_k \subset W_k$. \vs

At this point a wishful way to proceed is to take the values of $f$ on the rays $R_f(y'_k), R_f(x'_k)$ and use a quasiconformal interpolation to extend $f$ to the wake $W'_k$, keeping it univalent in a neighborhood of the hole $\La$. However, this is not always possible. In order to proceed, we define an auxiliary dynamics $g$ in a neighborhood of $z'_k$ as follows: Choose a sufficiently small simply connected neighborhood $D'_a$ of $z_a(t_j)$ such that $P_a^{-q}(D'_a)$ has a connected component $D_a \Subset D'_a$ which contains $z_a(t_j)$ and such that the restriction $P_a^{\circ q}: D_a \to D'_a$ is a biholomorphism. We may assume $\bd D'_a$ is a smooth Jordan curve which intersects each of the external rays $R_a(x_j),R_a(y_j)$ transversely once and intersects $\bd B_a$ at precisely two points. Let $\De_a, \De'_a$ be the respective connected components of $P_a^{-1}(D_a), P_a^{-1}(D'_a)$ containing $z'_a(t_k)$. We can arrange $\De_a \Subset \De'_a \Subset W_a(t_k)$ and the restriction $P_a: \De'_a \to D'_a$ to be a biholomorphism. Define $g_a: \De_a \to \De'_a$ as the conjugate of $P_a^{\circ q}$ by $P_a$, i.e., $g_a := P_a^{-1} \circ P_a^q \circ P_a$, where $P_a^{-1}$ is the local branch of the inverse mapping $D'_a$ to $\De'_a$. 
Set $\De : = \Phi(\De_a), \De' := \Phi(\De'_a)$ and let $g:\De \to \De'$ be the conjugate map $\Phi \circ g_a \circ \Phi^{-1}$. 
By the construction, $\De \cap \La$, $\De \sm \ov{\La}$, $\De \cap W'_k$ and $\De \sm \ov{W'_k}$ are all connected. Moreover,  
\begin{align*}
g(\De \cap \La) & = \De' \cap \La \\
g(\De \cap W'_k) & = \De' \cap W'_k, \\
f^{\circ q} \circ f & = f \circ g \quad \text{on} \quad \De \sm \ov{\La}, \\ f(\De\sm\ov{W'_k}) & = \Phi(D_a\sm\ov{B_a}) \cap W_j = \Phi(D_a \cap W_a({t_j})).
\end{align*} 
The assumption $(c,\nu) \in \check{\MM}_t \times \check{\MM}_0$ guarantees that the periodic orbit $\{ z_1, \ldots, z_q \}$ of $f$ is repelling. Thus, we can find simply connected neighborhoods $D \Subset D'$ of $z_j$ such that $f^{\circ q}: D \to D'$ is a biholomorphism. By taking $D'_a$ sufficiently small in the beginning of this construction, we can also arrange $\Phi(D'_a \sm B_a) \subseteq D'$. \vs  

Consider the quotient complex tori 
$$
T' := (\De \sm \{ z'_k \})/g \cong (\De_a \sm \{ z'_a(t_k) \})/ g_a \qquad \text{and} \qquad  T:= (D \sm \{ z_j \})/f^{\circ q}.
$$ 
Each of these tori is partitioned (modulo shared boundaries) into two isotopic and essentially embedded cylinders:
\begin{align*}
T' \ \text{partitions into} \quad & (\De \cap W'_k)/g \, \sqcup \, (\De \sm \ov{W'_k})/g \\
T \ \ \text{partitions into} \quad & (D \cap W_j)/f^{\circ q} \, \sqcup \, (D \sm \ov{W_j})/f^{\circ q}.
\end{align*}
The map $f$ induces a biholomorphism between $(\De \sm \ov{W'_k})/g$ and $(D \cap W_j)/f^{\circ q}$. However, we are interested in the two complementary cylinders 
$$
C' := (\De \cap W'_k)/g \qquad \text{and} \qquad C := ( D \sm \ov{W_j})/f^{\circ q},
$$ 
because the quasiconformal extension of $f$ is to be carefully carried out between $\De \cap W'_k$ and $D \sm \ov{W_j}$. \vs 

Consider the cylinder $C_\La:= (\De \cap \La)/g \Subset C'$. By the construction $C_\La$ is conformally isomorphic to the cylinder $(D_a \cap B_a)/P_a^{\circ q}$, which in turn is isomorphic to the standard cylinder $\H^+/\langle z \mapsto 2^q z \rangle$ of modulus $\pi/(q\log 2)$. The two connected components of $\De \sm K_f$ with $z'_k$ on their boundary project to cylinders in $T'$ of common modulus $\pi/(q\log 3)$ and the rays $R_f(x'_k)$ and $R_f(y'_k)$ project to their core geodesics, so half the modulus of each cylinder comes from the projection of $(\De \cap W'_k) \sm K_f$. Similarly, all the external rays landing at $z_j$ (such as $R_f(x_j), R_f(y_j)$ but possibly others due to merging) project to core geodesics of disjoint cylinders in $T$ coming from various components of $D \sm K_f$ each having the modulus $\pi/(q\log 3)$. \vs

We remark that $\Mod(C')$ can be made as large as we wish by changing the conformal structure on $\La$ and its preimages under $f$ (this fact will be used later in our quasiconformal interpolation step). Indeed, since $\Mod(C_\La) < \Mod(C')$, it suffices to show that $\Mod(C_\La)$ can be made arbitrarily large. The idea is of course standard. Let $m>0$ and $C(m)$ be any cylinder of modulus $m$. Take a quasiconformal homeomorphism $\psi: C_\La \to C(m)$. Let $\mu$ be the conformal structure on $\De \cap \La$ obtained by pulling back the standard conformal structure of $C(m)$ under $\psi$ and then lifting under $g$. Extend $\mu$ to the union of the iterated $f$-preimages of $\De \cap \La$ and define it to be the standard conformal structure in the rest of the plane. Evidently $\mu$ has bounded dilatation, so there is a quasiconformal homeomorphism $\phi: \C \to \C$ which integrates $\mu$ and is normalized by the conditions $\phi(0)=0, \phi(\infty)=\infty$ and $(\phi \circ f \circ \phi^{-1})(z)=\text{const.} + z^2 + O(z^3)$ near the origin. The map $f$ can now be replaced with $\phi \circ f \circ \phi^{-1}$ which still has the central and peripheral renormalizations hybrid equivalent to $Q_c$ and $Q_\nu$, but for which the corresponding cylinder $\phi(\De \cap \La)/(\phi \circ g \circ \phi^{-1})$ has modulus equal to $m$. Thus, without loss of generality, we may assume that the cylinder $C_\La$ (and therefore $C'$) for our map $f$ has arbitrarily large modulus. 

\subsection{A quasiregular extension of $f$}\label{extf}

The external rays of $f$ landing on the orbit $\{ z_1, \ldots, z_q \}$ are related to the external rays of $Q_c$ landing on the orbit $\{ w_i :=\Phi^{-1}(z_i) \}$. Since $f$ has a central renormalization hybrid equivalent to $Q_c$, the discussion of \S \ref{qpms} shows that there are three possibilities according as the portrait of $\{ w_1, \ldots, w_q \}$ is trivial, primitive or satellite: (i) $\{ z_1, \ldots, z_q \}$ has period $q$ and each $z_i$ has exactly two rays $R_f(x_i),R_f(y_i)$ landing on it; (ii) $\{ z_1, \ldots, z_q \}$ has period $q$ and each $z_i$ has exactly three rays  $R_f(x_i),R_f(y_i),R_f(u_i)$ with disjoint orbits landing on it; (iii) $\{ z_1, \ldots, z_q \}$ has a smaller period $p=q/r$ and every $z_i$ has exactly $2r$ rays landing on it, $r$ rays in the orbit of $R_f(x_i)$ and $r$ rays in the orbit of $R_f(y_i)$. Whatever the case, the number $n$ of the  external rays landing at each $z_i$ satisfies $2 \leq n \leq 2q$. Denote by $R_1,\ldots, R_n$ the rays landing at $z_j=f(z_k)$, where the labeling is counterclockwise starting with $R_1 = R_f(y_j)$ and ending with $R_n = R_f(x_j)$. \vs
 
Let $K_i$ be the intersection of $K_f$ with the sector bounded by the closures of the rays $R_i$ and $R_{i+1}$, taking the index $i$ modulo $n$. Evidently $K_n \subset W_j$. Set 
$$
m_i := \Mod_{z_j}(K_i) \qquad (1 \leq i \leq n) 
$$
as defined in \S \ref{qcint} and note that the $m_i$ are unaffected by the change of the conformal structure on $\La$ described at the end of the previous section. We shall thus suppose 
\begin{equation}\label{modlowbound}
m = \Mod(C_\La) \geq \sum_{i=1}^{n-1} m_i + (n-1)\ \frac{\pi}{q\log 3}.
\end{equation}
For each $1 \leq i \leq n-1$ replace the rays $R_i, R_{i+1}$ surrounding $K_i$ by the oblique rays $\Ga_{2i}:=R^+_i,\Ga_{2i+1}:=R^-_{i+1}$ of slopes $\varsigma_i,-\varsigma_i$. Setting $\Ga_1:=R_1, \Ga_{2n}:=R_n$, we obtain $2n$ rays $\Ga_1, \ldots, \Ga_{2n}$ arranged counterclockwise around $z_j$ which we may assume do not intersect each other in $D$. The quotient under $f^{\circ q}$ of the sector in $D$ delimited by $\Ga_{2i},\Ga_{2i+1}$ is a sub-cylinder $C_{2i}$ of $C=( D \sm \ov{W_j})/f^{\circ q}$. By the discussion in \S \ref{qcint}, $\Mod(C_{2i})$ tends to $m_i$ as $\varsigma_i \to \infty$ and it tends to a limit $\geq m_i+\pi/(q \log 3)$ as $\varsigma_i \to 0$. Thus, we can choose $\varsigma_i>0$ such that  
\begin{equation}\label{mm1}
\Mod(C_{2i}) = m_i + \frac{\pi}{2q\log 3} \qquad (1 \leq i \leq n-1).
\end{equation}
For $1 \leq i \leq n$ the quotient under $f^{\circ q}$ of the sector in $D$ delimited by $\Ga_{2i-1},\Ga_{2i}$ is a sub-cylinder of $C$ that we denote by $C_{2i-1}$.  \vs

Choose a sufficiently small Green's potential $s>0$ such that 
the rays $\Ga_1, \ldots, \Ga_{2n}$ are disjoint up to level $s$. 
Let $U$ and $U'$ be the topological disks bounded by the equipotentials of $f$ at levels $s/3$ and $s$, respectively. We replace each ray $\Ga_i$ with the arc $\Ga_i \cap \ov{U'}$ and for simplicity continue to denote the trimmed ray by $\Ga_i$. For $1 \leq i \leq 2n-1$ we denote by $\Si_i$ the sector in $U'$ with the vertex at $z_j$ that is delimited by $\Ga_i$ and $\Ga_{i+1}$. Thus, $(\Si_i \cap D)/f^{\circ q} = C_i$. Observe that $\Si_{2i} \cap K_f = K_i$ while $\Si_{2i-1} \cap K_f = \es$. \vs 

Let us now carry out a similar construction near $z'_k$ for the cylinder $C' = (\De \cap W'_k)/g$ bounded by the Jordan curves $\ga_1 := R_f(y'_k)/g$ and $\ga_{2n} := R_f(x'_k)/g$. Supplement these curves with $2n-2$ additional disjoint smooth Jordan curves $\ga_2, \ldots \ga_{2n-1}$ in $C'$, all in the same isotopy class as $\ga_1,\ga_{2n}$ and labeled consecutively so as to obtain pairwise disjoint cylinders $C'_i \subset C'$ bounded by $\ga_i$ and $\ga_{i+1}$ for $1 \leq i \leq 2n-1$. By the choice \eqref{modlowbound} we can choose the $\ga_i$ such that  
\begin{equation}\label{mm2}
\Mod(C'_{2i}) = m_i + \frac{\pi}{2q\log 3} \qquad (1 \leq i \leq n-1) \\
\end{equation}
We can lift and extend $\ga_1, \ldots, \ga_{2n}$ to global arcs $\Ga'_1, \ldots, \Ga'_{2n}$ in $\ov{W_k}$ as follows. Set $\Ga'_1 := R_f(y'_k) \cap \ov{U}, \Ga'_{2n} := R_f(x'_k) \cap \ov{U}$, and for $2 \leq i \leq 2n-1$ choose a smooth embedded arc $\Ga'_i$ connecting $z'_k$ within $W'_k$ to $\bd U$ such that $(\Ga'_i \cap \De)/g = \ga_i$. We may assume that $\Ga'_1, \ldots, \Ga'_{2n}$ are pairwise disjoint except that their common extremity at $z'_k$. For $1 \leq i \leq 2n-1$ denote by $\Si'_i$ the sector in $U$ with the vertex at $z'_k$ that is delimited by $\Ga'_i$ and $\Ga'_{i+1}$, so $(\Si'_i \cap \De)/g = C'_i$. \vs

Let $\zeta_{2i}: \Si'_{2i} \to \Si_{2i}$ be the conformal isomorphism which maps $z'_k$ to $z_j$, $\Ga'_{2i} \cap \bd U$ to $\Ga_{2i} \cap \bd U'$ and $\Ga'_{2i+1} \cap \bd U$ to $\Ga_{2i+1} \cap \bd U'$. By \eqref{mm1} and \eqref{mm2}, 
$$
\Mod((\Si'_{2i} \cap \De)/g)=\Mod(C'_{2i})=\Mod(C_{2i})=\Mod((\Si_{2i} \cap D)/f^{\circ q}).
$$
It follows from \lemref{bf1} that the boundary maps $\zeta_{2i}: \Ga'_{2i} \to \Ga_{2i}$ and $\zeta_{2i} : \Ga'_{2i+1} \to \Ga_{2i+1}$ are near-linear for $1 \leq i \leq n-1$. Moreover, the restrictions $f: \Ga'_1 \to \Ga_1$ and $f: \Ga'_{2n} \to \Ga_{2n}$ are also near-linear since they conjugate $g$ to $f^{\circ q}$. \lemref{bf2} now shows that these boundary maps yield quasiconformal extensions $\zeta_{2i-1}: \Si'_{2i-1} \to \Si_{2i-1}$ on the complementary sectors. We may thus define a proper quasiregular map $F: U \to U'$ of degree $3$ by
$$
F(z) := 
\begin{cases}
f(z) & \qquad z \in U \sm W'_k \\
\zeta_i(z) & \qquad z \in \Si'_i \quad (1 \leq i \leq 2n-1).
\end{cases}
$$

\subsection{Completion of the proof of \thmref{surject}} \label{pfsur}

The set where $\bar{\bd}F \neq 0$ is contained in the union ${E'}_{\!\!\odd} := \bigcup_{i=1}^n \ov{\Si'_{2i-1}} \subset \ov{U \cap W'_k}$. Each forward orbit of $F$ can meet ${E'}_{\!\!\odd}$ at most once because $F({E'}_{\!\!\odd}) = E_\odd := \bigcup_{i=1}^n \ov{\Si_{2i-1}}$, where $E_\odd \cap {E'}_{\!\!\odd}= \es$ and all further iterates of $E_\odd$ under $F$ are disjoint from $W_k$ and contained in the escaping set $\C \sm K_f$. Define a conformal structure $\mu$ on $U'$ by taking the standard conformal structure on $E_\odd$ and extending it to the union $\bigcup_{i \geq 0} F^{-i}(E_\odd)$ by iterated pull-backs. On the rest of $U'$ let $\mu$ be the standard conformal structure. Notice that $\mu$ is well-defined since the branches of $F^{-q}=f^{-q}$ that map $E_\odd$ to itself are conformal. Moreover, the iterated branches $F^{-i}$ can contain at most one non-conformal branch of $F^{-1}$ by the above observation. Thus, $\mu$ is $F$-invariant and has bounded dilatation. \vs

Let $\Psi: U' \to \D$ be a quasiconformal homeomorphism that integrates $\mu$ with say $\Psi(0)=0$. Then $G:= \Psi \circ F \circ \Psi^{-1}: \Psi(U) \to \D$ is a polynomial-like map of degree $3$. Moreover, by the construction of $F=f$ in $U \sm W'_k$, the map $G$ has a quadratic-like restriction around the critical point $(\Psi \circ \Phi)(\om_1)=0$ which is hybrid equivalent to $Q_c$ and $G^{\circ q}$ has a quadratic-like restriction around the critical point $(\Psi \circ \Phi)(\om_2)$ which is hybrid equivalent to $Q_\nu$. A suitably normalized straightening of $G$ will then be the desired cubic $P \in \NN(t)$ with $\Chi(P)=(c,\nu)$. \hfill $\Box$

\section*{Appendix. On the intersection of lemon limbs}

Although distinct lemon limbs are disjoint, there could be plenty of non-trivial intersections between a limb $\LL(t)$ and a rotated limb $\LL^\ast(s)$. This appendix provides some relevant results and examples; they will be mostly independent from the rest of this paper. \vs

Throughout the following discussion we assume $t,s$ are periodic angles under $\m_2$ with periods $q, p \geq 2$ and combinatorics $\sigma \in S_q, \tau \in S_p$, respectively. As usual, let $\{ x_1, \ldots, x_q \}$ and $\{ y_1, \ldots, y_q \}$ be the simulating orbits for $t=t_k$ with simulating angles $x_k,y_k$. Similarly, let $\{ u_1, \ldots, u_p \}$ and $\{ v_1, \ldots, v_p \}$ be the simulating orbits for $s=s_\ell$ with simulating angles $u_\ell, v_\ell$. Set $I_i:=[x_i,y_i]$ for $1 \leq i \leq q$ and $J_j:=[u_j,v_j]$ for $1 \leq j \leq p$. Recall from \S \ref{subsec:sp} that $I_{\sigma(k)}$ and $J_{\tau(\ell)}$ are the shortest intervals in their respective collections with $|I_{\sigma(k)}|=1/(3^q-1)$ and $|J_{\tau(\ell)}|=1/(3^p-1)$. Similarly, $I_k$ and $J_\ell$ are the longest intervals in their respective collections with $|I_k|=3^{q-1}/(3^q-1)$ and $|J_\ell|=3^{p-1}/(3^p-1)$.  

\begin{theorem*}
$\LL(t) \cap \LL(s) = \es$ whenever $t \neq s$. 
\end{theorem*}

\PROOF
Assume by way of contradiction that there is a cubic $P \in \LL(t) \cap \LL(s)$ and $t \neq s$. For $1 \leq i \leq q$ let $z_i$ denote the co-landing point of the rays $R(x_i),R(y_i)$ and for $1 \leq j \leq p$ let $w_j$ denote the co-landing point of $R(u_j),R(v_j)$. Setting $y'_k:=y_k-1/3 < 1/2 < x'_k:=x_k+1/3$ as before, the rays $R(y'_k),R(x'_k)$ co-land at a preimage of $z'_k$ of $z_{\sigma(k)}$ other than $z_k$. Since $w_\ell$ is periodic but $z'_k$ is not, the arcs $R(u_\ell) \cup R(v_\ell) \cup \{ w_\ell \}$ and $R(x'_k) \cup R(y'_k) \cup \{ z'_k \}$ are disjoint. The length comparison  $|J_\ell|>1/3>|[y'_k,x'_k]|$ then shows that 
\begin{equation}\label{uxyv}
u_\ell<y'_k<1/2<x'_k<v_\ell. 
\end{equation}
As the intervals $I_k$ and $J_\ell$ both contain $1/2$, there are only two possibilities: \vs

(i) $I_k$ and $J_\ell$ are not nested. Then, by \eqref{uxyv}, 
\begin{equation}\label{xyuv}
x_k<u_\ell<y'_k<1/2<y_k<v_\ell \qquad \text{or} \qquad  u_\ell<x_k<1/2<x'_k<v_\ell<y_k. 
\end{equation}
In either case the arcs $R(x_k) \cup R(y_k) \cup \{ z_k \}$ and $R(u_\ell) \cup R(v_\ell) \cup \{ w_\ell \}$ necessarily meet, so $z_k=w_\ell$. Since the ray $R(x_k)$ landing at $z_k$ is fixed under $P^{\circ q}$, every ray landing at $z_k$ (such as $R(u_\ell)$) must be fixed under $P^{\circ q}$. This shows $\m_3^{\circ q}(u_\ell)=u_\ell$, so $p|q$. A similar argument shows $q|p$, and we conclude that $p=q$. In particular, every $u_j,v_j$ has period $q$ under $\m_3$. Now \eqref{xyuv} together with the fact that $\m_3$ maps each of $[x_k,y'_k]$ and $[x'_k,y_k]$ homeomorphically onto $I_{\sigma(k)}$ shows that either $u_{\tau(\ell)}$ or $v_{\tau(\ell)}$ lies in the interior of $I_{\sigma(k)}$. This is a contradiction since the endpoints $x_{\sigma(k)}$ and $y_{\sigma(k)}$ of $I_{\sigma(k)}$, being $1/(3^q-1)$ apart, are adjacent period $q$ points of $\m_3$. \vs

(ii) $I_k$ and $J_\ell$ are nested. Assume without loss of generality that $q \leq p$. If $q=p$, then $|I_k|=3^{q-1}/(3^q-1)=|J_\ell|$, so $x_k=u_\ell, y_k=v_\ell$ and therefore $t=s$, contradicting our hypothesis. Let us then assume that $q<p$ so $|I_k|>|J_\ell|$. By \eqref{uxyv},
$$
x_k \leq u_\ell<y'_k<1/2<x'_k<v_\ell \leq y_k. 
$$
Since the iterate $\m_3^{\circ q}$ is an expanding diffeomorphism from each of the intervals $[x_k,y'_k]$ and $[x'_k,y_k]$ onto $I_k$, it follows that $x_k \leq u_\ell \leq u_{\tau^q(\ell)}<y_k$ and $x_k<v_{\tau^q(\ell)} \leq v_\ell \leq y_k$, so $J_\ell$ and $J_{\tau^q(\ell)}$ meet. This is a contradiction since by $q<p$ the intervals $J_\ell$ and $J_{\tau^q(\ell)}$ must be disjoint.  
\ENDPROOF

It follows from the above theorem and symmetry that the rotated limbs $\LL^\ast(t), \LL^\ast(s)$ are also disjoint whenever $t \neq s$. However, the intersection $\LL(t) \cap \LL^\ast(s)$ may well be non-empty in both cases $t=s$ and $t \neq s$. \vs

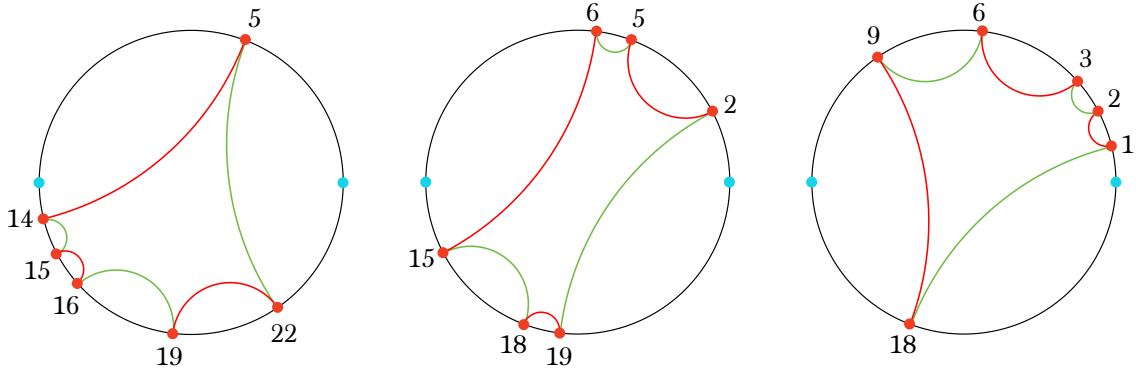
\begin{figure}[t]
	\centering
	\begin{overpic}[width=\textwidth]{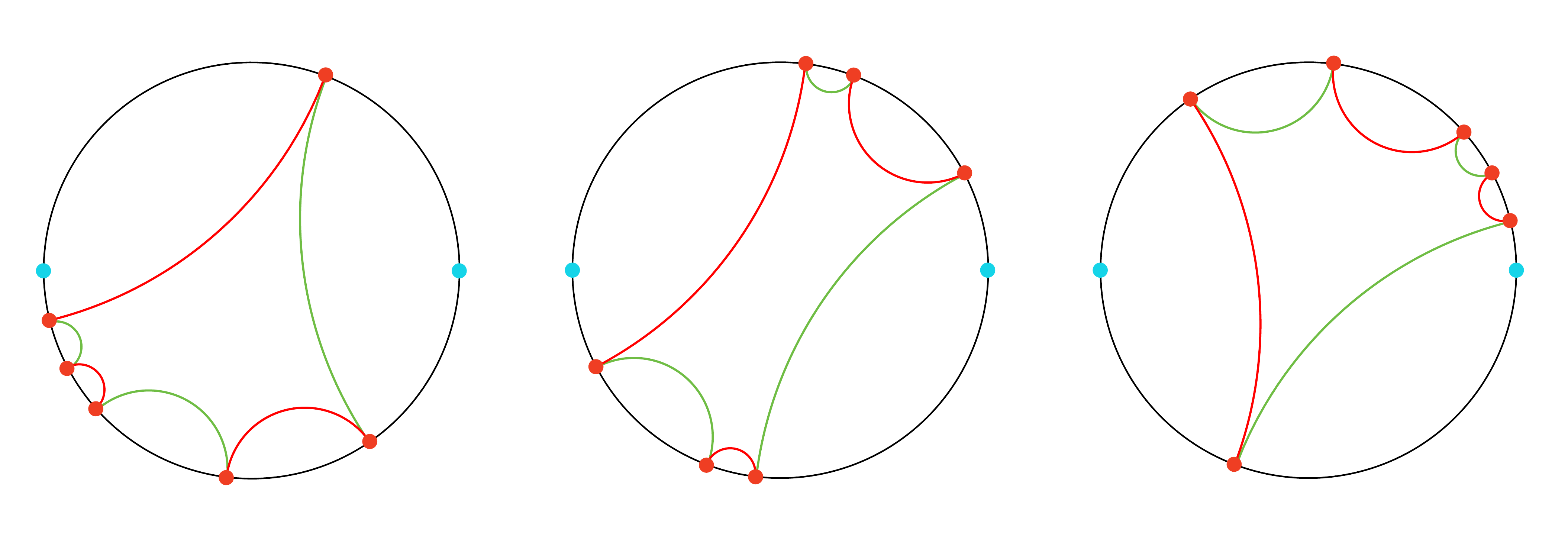}
		\put (21,30.5) {\footnotesize $5$}
		\put (0,13) {\footnotesize $14$}
		\put (1.5,8.5) {\footnotesize $15$}
		\put (4,5.4) {\footnotesize $16$}
		\put (13,1) {\footnotesize $19$}
		\put (23,3) {\footnotesize $22$}
		\put (50.5,31) {\footnotesize $6$}
		\put (35,9.5) {\footnotesize $15$}
		\put (43,2) {\footnotesize $18$}
		\put (47,1) {\footnotesize $19$}
		\put (54.5,30.5) {\footnotesize $5$}
		\put (62.5,23) {\footnotesize $2$}
		\put (97.2,19.5) {\footnotesize $1$}
		\put (96,23.5) {\footnotesize $2$}
		\put (93.5,27) {\footnotesize $3$}
		\put (84.2,31) {\footnotesize $6$}
		\put (75,29) {\footnotesize $9$}
		\put (77,2) {\footnotesize $18$}
	\end{overpic}
	\caption{\small The complementary co-landing patterns in the intersections $\LL(1/7) \cap \LL^\ast(4/7)$ on the left, $\LL(2/7) \cap \LL^\ast(2/7)$ in the middle, and $\LL(4/7) \cap \LL^\ast(1/7)$ on the right (angles are shown in multiplies of $1/26$). All three patterns can be realized by the rays landing on a repelling or degenerate parabolic fixed point of combinatorial rotation number $1/3$.}   
	\label{LL}
\end{figure}

For example, suppose $t,s$ are {\bit complementary angles} in the sense that they belong to the same $\m_2$-orbit $\{ t_1, \ldots, t_q \}$ of rotation type (i.e., having combinatorics of degree $1$) with $t=t_k$ and $s=t_{q-k+1}$ for some $1 \leq k \leq q$. We claim that $\LL(t) \cap \LL^\ast(s) \neq \es$. In fact, if $\OO_0, \ldots, \OO_q$ are the $\m_3$-orbits associated with $\{ t_1, \ldots, t_q \}$ as in \S \ref{subsec:sp}, then the simulating orbits in $\LL(t)$ are $\OO_{k-1}, \OO_k$ while those in $\LL^\ast(s)$ are $\OO^\ast_{q-k}, \OO^\ast_{q-k+1}$, also equal to $\OO_k,\OO_{k-1}$ by \eqref{invol}. In other words, the co-landing rays coming from $\LL(t)$ and $\LL^\ast(s)$ have the same angles but their pairings are complementary. Thus, for any $P \in \LL(t) \cap \LL^\ast(s)$, the $2q$ rays $R(\theta)$ for $\theta \in \OO_{k-1} \cup \OO_k$ land at a common fixed point of $P$, repelling or parabolic with combinatorial rotation number of the form $p/q$. Examples of such cubics can be easily constructed as suitable intertwinings of pairs of quadratics in the $p/q$-limb of the Mandelbrot set (compare \figref{LL}). \vs

%%%%%%%%%%%%%%%%%%%%%%%%%%%%%%%%%%%%%%%%%%%
%%%%%%%%%%%%%%%%%%%%%%%%%%%%%%%%%%%%%%%%%%%   

Beyond the case of complementary angles, many other intersections are possible. For the simplest example with non-rotational combinatorics, take the $\m_2$-orbit  $\{ 1, 2, 3, 4 \}/5$ with $\sigma=(1243)$. One can verify that the six intersections 
\begin{align*}
\LL(1/5) \cap \LL^\ast(1/5) & & \LL(1/5) \cap \LL^\ast(3/5) & & \LL(2/5) \cap \LL^\ast(4/5) \\
\LL(3/5) \cap \LL^\ast(1/5) & & \LL(4/5) \cap \LL^\ast(1/5) & & \LL(4/5) \cap \LL^\ast(4/5)  
\end{align*}
are all non-empty. More dramatically, the following result, which we state without proof, shows that a limb can intersect infinitely many rotated limbs: 

\begin{theorem*}
$\LL(1/(2^{q_0}-1)) \cap \LL^\ast(1/(2^q-1)) \neq \es$ whenever $q \geq q_0 \geq 3$. 
\end{theorem*}

It would be interesting to find an exact characterization of the angles $t,s$ for which $\LL(t) \cap \LL^\ast(s) \neq \es$. All such intersections for non-complementary angles $t,s$ must have a combinatorics that is inherently cubic. Explicitly, if $\{ x_i, y_i \}_{1 \leq i \leq q}$ is the union of the simulating orbits of $t$ and $\{ u_j, v_j \}_{1 \leq j \leq p}$ is the union of the simulating orbits of $s$, then the combinatorics of $\{ x_i, y_i \}_{1 \leq i \leq q} \cup \{ u_j+1/2, v_j+1/2 \}_{1 \leq j \leq p}$ under $\m_3$, as a permutation in $S_{2p+2q}$, has degree $3$. \vs

\noindent
{\bf Acknowledgments.} C.L.P. would like to thank the Danish Council for Independent Research -- Natural Sciences for support via the grant DFF-1026-00267B. S.Z.  acknowledges the partial support of the Research Foundation of The City University of New York via grant TRADB-53-144.

\end{document}